\documentclass[12pt,a4paper,reqno]{amsart}

\usepackage{amssymb}
\usepackage{amscd}
\usepackage{amsfonts}
\usepackage{graphicx}
\usepackage[british]{babel}
\usepackage{enumerate}

\numberwithin{equation}{section}

\theoremstyle{plain}
  \newtheorem{theorem}{Theorem}[section]
  \newtheorem{proposition}[theorem]{Proposition}
  \newtheorem{lemma}[theorem]{Lemma}
  \newtheorem{corollary}[theorem]{Corollary}
  
  \newtheorem{definition}[theorem]{Definition}
  \newtheorem{claim}{Claim}
  \newtheorem*{wredmain}{Proposition \ref{W-red-Main-thm}}

\theoremstyle{remark}
  \newtheorem*{remarks}{Remarks}
  \newtheorem*{remark}{Remark}

\renewcommand{\leq}{\leqslant}
\renewcommand{\geq}{\geqslant}

\newcommand{\stsubsection}[1]{\subsection*[#1]{\sc #1}}

\newcommand\E{\mathbb{E}}
\newcommand\Z{\mathbb{Z}}
\newcommand\R{\mathbb{R}}

\newcommand\C{\mathbb{C}}

\newcommand\A{\mathcal{A}}
\newcommand\N{\mathbb{N}}
\newcommand\W{\overline{W}}

\renewcommand\P{\mathcal{P}}

\newcommand\Q{\mathcal{Q}}
\newcommand\eps{\varepsilon}
\newcommand\vol{\operatorname{vol}}
\newcommand\lcm{\operatorname{lcm}}

\newcommand\pdspan{\langle \mathcal P_D \rangle}

\newcommand\starpdspan{\langle {\mathcal P}^*_D \rangle}
\newcommand\starqdspan{\langle {\mathcal Q}^*_D \rangle}
\newcommand\qdspan{\langle \mathcal Q_D \rangle}
\newcommand{\kron}[2]{\left(\genfrac{}{}{.5pt}{}{#1}{#2}\right)}
\newcommand{\dsum}{\sideset{}{^{\prime}}{\sum}}
\newcommand{\Mod}[1]{\;(\mathrm{mod}\,#1)}
\newcommand\<{\langle}
\renewcommand\>{\rangle}

\begin{document}
\title[Linear correlations amongst numbers weighted by $R_f$]{Linear
correlations amongst numbers represented by positive definite binary
quadratic forms}

\author{Lilian Matthiesen}
\address{School of Mathematics\\
University Walk\\
Bristol, BS8 1TW\\
United Kingdom
}
\email{l.matthiesen@bristol.ac.uk}

\begin{abstract}
Let $f_1, \dots, f_t$ be positive definite binary quadratic forms, and
let $R_{f_i}(n)=|\{(x,y): f_i(x,y)=n\}|$ denote the corresponding
representation functions. 
Employing methods developed by Green and Tao, we deduce asymptotics for
linear correlations of these representation functions.
More precisely, we study the expression
$$\E_{n \in K \cap [-N,N]^d} \prod_{i=1}^t R_{f_i}(\psi_i(n))~,$$
where the $\psi_i$ form a system of affine linear forms no two of which
are affinely related, and where $K$ is a convex body.

The minor arc analysis builds on the observation that polynomial
subsequences of equidistributed nilsequences are still equidistributed,
an observation that could be useful in treating the minor arcs of other
arithmetic questions.

As a very quick application we give asymptotics to the number of
simultaneous zeros of certain systems of quadratic equations in 8
or more variables.
\end{abstract}

\maketitle
\tableofcontents

\addtocontents{toc}{\protect\setcounter{tocdepth}{1}}

\section{Introduction}
The distribution of prime numbers shares many properties with the
distribution of numbers that are representable as a sum of two squares,
an analogy that is occasionally employed to obtain model problems for
questions about the primes.
Let us consider the distributions of the two sets in arithmetic
progressions.
Starting with the average orders, we have on the one hand the Prime
Number Theorem, asserting that $\pi(x) \sim \frac{x}{\log x}$. 
For the set $S$ of sums of two squares on the other hand, Landau
\cite{landau} proved an analogous asymptotic in 1908, namely
$$\sum_{n \leq x} 1_S(n) \sim B \frac{x}{\sqrt{\log x}}~,$$ 
where $B= \frac{1}{\sqrt{2}} \prod_{p\equiv 3 \Mod{4}}(1-p^{-2})^{-1/2}$.

Turning towards more general arithmetic progressions, let $a$ and $q$ be
coprime integers, then the primes congruent to $a \Mod{q}$ satisfy 
$\pi(x;q,a) \sim \frac{x}{\phi(q) \log x}$. 
Building on Landau's result and the analogy to primes, Prachar
\cite{prachar} proved in the 1950s that, when furthermore 
$a\equiv 1 \Mod{\gcd(4,q)}$ holds, sums of two squares show the following
behaviour\footnote{This compact formulation of the result is due to
Iwaniec \cite{iwaniec}.} 
$$\sum_{\substack{n \leq x \\ n \equiv a \Mod{q}}} 1_S(n) 
\sim B_q \frac{x}{\sqrt{\log x}}~,$$
where
$$B_q 
= B q^{-1} \frac{(4,q)}{(2,q)} 
  \prod_{\substack{p\equiv 3 \Mod{4} \\ p|q}} (1+ p^{-1})~.$$
The factor $(4,q)/(2,q)$ describes that the density of sums
of two squares is twice as high in the progression $n\equiv 1 \Mod{4}$ as
it is in $n \equiv 1 \Mod{2}$.
From pairs $(a,q)$ that are not coprime, one needs to remove
those choices from consideration that lead to whole progressions entirely
lying outside the set $S$. 
Examples are integers $n$ such that $n \equiv 3 \Mod{3^2}$, which are
never a sum of two squares, or, 
numbers of the form $(3 \cdot 5) n + 3^2 $, which can only be a sum of two
squares when $3|n$.
When excluding such classes $a \Mod{q}$, the constant $B_q$
only needs to be adapted by restricting the product over $p\equiv 3
\Mod{4}$ to primes dividing $q/\gcd(a,q)$.

Thus, both sets, the primes and the sums of two squares, show some 
uniformity in the distribution in residue classes once one excludes
residue classes that for obvious reason contain too few elements.

It is natural to ask whether this uniformity carries further:
is, for instance, the distribution uniform enough to determine
asymptotically the density of solutions to linear equations within these
sets? 
More precisely, we are interested in an asymptotic for correlations
of the form
$$\sum_{n \in \Z^d \cap K} \prod_{i=1}^t f(\psi_i(n))~,$$
where the $\psi_i : \Z^d \to \Z$ are affine linear forms and where the
arithmetic function $f:\Z \to \R$ is either the characteristic function
$1_S$ of sums of two squares, or it is chosen to be the characteristic
function of primes.

Green and Tao studied these correlations in the case of the primes in
\cite{green-tao-linearprimes}.
They replaced for this purpose the sparse set of primes by a weighted
version of asymptotic density $1$ which is given by the von Mangoldt
function. 
We shall \emph{not} normalise the characteristic function $1_S$ in an
analogous way, but instead consider the intrinsically weighted function
that is given by the representation function of sums of two squares, 
$R(n) = |\{(x,y) \in \Z:x^2+y^2=n\}|$. 
Counting lattice points in a circle of radius $\sqrt{N}$ immediately shows
that the representation function has indeed an asymptotic density given by
$\sum_{n \leq N} R(n) \sim \pi N$.
As we will see, the nilpotent Hardy-Littlewood method which Green and Tao
developed to handle linear correlations among the primes can also be
employed in the case of the representation function $R(n)$.

Instead of restricting attention to representations as sums of two
squares, the slightly more general case of representation by a positive
definite binary quadratic form $f(x,y)= ax^2+bxy+cy^2$ will be
considered. 
The corresponding representation function is then given
by $R_f(n) =|\{(x,y): f(x,y)=n\}|$.

\subsection*{Notation.} Throughout the paper, we write $[N]$ for the set
of numbers $\{1,\dots,N\}$ and $e(x)$ for $\exp(2\pi x)$.
We let $v_p: \N \to \N_0$ denote the $p$-adic valuation.  
If $T$ is a finite set, we use the expectation notation $\E_{t \in T}$
to abbreviate $\frac{1}{|T|}\sum_{t \in T}$.
A linear correlation is defined along a system
$\Psi=(\psi_1,\dots,\psi_t): \Z^d \to \Z^t$ of affine linear forms.
Such a system may be written as $\Psi(n)= \Psi(0) + \dot\Psi(n)$, for a
linear map $\dot\Psi$.
We regard $\dot\Psi$ as fixed, while $\Psi(0)$ may, for instance,
depend on $K$. 
Thus, all implicit constants in asymptotic notation, such as $O()$,
$o()$ and $\ll$, are allowed to depend on the coefficients of the
linear map $\dot \Psi$, the dimensions $d$ and $t$ of the domain and the
target space of $\Psi$, and on the discriminants of the forms 
$f_1,\dots,f_t$.

\subsection*{Methods and results}
The nilpotent Hardy-Littlewood method comprises a scheme that allows
to obtain for any given arithmetic function $h$ with sufficiently
quasi-random behaviour an asymptotic for the expression
$$\sum_{n \in K \cap \Z^d} \prod_{i=1}^t h(\psi_i(n))~,$$
where $K \subset \R^d$ is convex and satisfies $\psi_i \subset [1,N]$ for
each $i \in \{1,\dots,t\}$. 
We proceed to describe the basic set-up for the method. There are two main
requirements on $h$. One is that $h$ has small Gowers-uniformity norms
(see Section \ref{application_section}) and the other is that one can find
a majorant $\nu: \{1,\dots, N\} \to \R_{>0}$ such that 
\begin{enumerate}
 \item[1.] (\emph{majorant property}) the point-wise estimate 
$h(n) \leq C\nu(n)$ is satisfied for an absolute constant $C$
independent of $N$,
 \item[2.] (\emph{density condition}) $h$ has positive relative density in
$\nu$ in the sense that \\ $\E_{n\leq N} h(n) \sim C' \E_{n \leq N}
\nu(n)$,
 \item[3.] (\emph{pseudorandomness}) $\nu$ is a pseudorandom measure.
\end{enumerate}

A \emph{pseudorandom measure} resembles a true probability measure.
Apart from requiring its total mass to be approximately $1$, that is,
$\E_{n \leq N} \nu(n) = 1+o(1)$, there are two further defining
conditions for a pseudorandom measure: 
the linear forms condition and the correlation condition.
Each of them places some independence requirements upon $\nu$.
The linear forms condition for instance requires 
$$\E_{n \in K \cap \Z^d} \prod_{i=1}^t \nu(\psi_i(n)) =1+o(1)$$
to hold for certain systems of affine linear forms.
Once in possession of such a pseudorandom majorant, a number of
tools are available. We will describe them as we encounter them.

Regarding the first condition on $h$, which was the smallness of certain
Gowers-uniformity norms, there is an explicit (and quite strong) necessary
condition that has to be satisfied:
$h$ must be equidistributed in residue classes to small moduli.
The results quoted at the beginning of this introduction show that
neither the characteristic function $1_S(n)$, nor its weighted version
$r(n)$ meet this requirement. 
In such a situation, it may be possible to decompose the function $h$ into
a sum of functions that are more uniformly behaved and consider each of
these functions separately. 
This decomposition is known as $W$-trick and will be carried out in
Section \ref{W-trick}.

In Section \ref{rep-majorant-section} we construct a majorant for the
representation function attached to a primitive form $f$.
This majorant will be slightly modified in accordance to the $W$-trick in
Section \ref{W-trick}.
In Section \ref{linear-forms-and-correlation-conditions} we check that our
majorant is indeed pseudorandom.

In the course of the minor arc analysis, which starts in Section
\ref{non-corr-intro}, we observe that polynomial subsequences of
$\delta$-equidistributed linear nilsequences are still reasonably
equidistributed. 
See Proposition \ref{equid_subsequences} below.
This result will be deduced from the quantitative equidistribution
theory Green and Tao worked out in \cite{green-tao-polynomialorbits}.
In connection with their factorisation theorem \cite[Thm
1.19]{green-tao-polynomialorbits}, it could prove a useful tool
for the minor arc analysis of a wider range of arithmetic problems.

Due to the quite complex foundations of the Green-Tao methods it proved
not feasible to provide a self-contained account of it here. 
This paper therefore strongly depends on \cite{green-tao-linearprimes}.
It furthermore relies on results about the divisor function from
\cite{m-divisorfunction}, which will be used in the construction of the
pseudorandom majorants.

\subsection*{Results}
In \cite{m-divisorfunction} a pseudorandom majorant for the
normalised divisor function 
$\tilde\tau(n) = (\log N)^{-1}\sum_{d|n} 1$ has been constructed. 
Here we shall combine this majorant with a sieving majorant to
obtain a pseudorandom majorant for the function $R_f(n)$ which counts
the number of representations of $n$ by a primitive positive definite
binary quadratic form; results for the non-primitive case are
immediate corollaries.

With this majorant at hand, we obtain, employing the machinery from
\cite{green-tao-longprimeaps,green-tao-linearprimes} in combination
with the inverse theorem for the Gowers-uniformity norms \cite{gtz}, an
asymptotic for the representation function $R_f$ evaluated along systems
of linear equations:
\begin{theorem}\label{main thm}
Let $f_1, \dots f_t$ be primitive positive definite binary quadratic
forms.
Let $\Psi= (\psi_1, \dots ,\psi_t): \Z^d \to \Z^t$ be a system of
affine linear forms such that no two forms $\psi_i$ and $\psi_j$ are
affinely dependent. 
Suppose that the coefficients of the linear part $\dot \Psi$ are bounded
and that $K \subset [-N,N]^d$ is a convex body such that $\Psi(K) \subset
[0,N]^t$.
Then 
\begin{align*}
 \sum_{n\in \Z^d\cap K } R_{f_1}(\psi_1(n)) \dots R_{f_t}(\psi_t(n)) 
= \beta_{\infty} 
   \prod_{p} \beta_p
 +  o(N^d)~,
\end{align*}
where
$$\beta_{\infty}=\vol(K) \prod_{i=1}^t\frac{2\pi}{\sqrt{-D_i}}~,$$
and
$$\beta_p 
= \lim_{m \to \infty}
  \E_{a \in(\Z/p^{m}\Z)^{d}} 
  \prod_{i \in [t]}
  \frac{\rho_{f_i,\psi_i(a)}(p^{m})}{p^{m}}~,$$
with $\rho_{f,A}(q)$ denoting the local number of representations of
$A \Mod{q}$ by $f$, that is,
$$\rho_{f,A}(q)
:= |\{(x,y)\in[q]^2:f(x,y) \equiv A \Mod{q}\}|~.$$
\end{theorem}

Theorem \ref{main thm} extends previous results by Heath-Brown
\cite{heath-brown} and improvements thereof by Browning and de la
Bret{\`e}che \cite{browning-breteche}, where the case of sums of two
squares, $f_i(x,y)=x^2+y^2$, for $i=1, \dots, 4$, together with systems 
$\Psi: \Z^2 \to \Z^4$ was considered. 
We emphasise, however, that, in contrast to the results from
\cite{heath-brown} and \cite{browning-breteche}, we unfortunately do not
obtain explicit error terms in our asymptotic.

The most interesting case of correlation along a system of affine linear
forms is certainly the `infinite complexity' case of
$$\E_{n \leq N} R_f(n+a_1)\dots R_f(n+a_d)~,$$
corresponding to the prime-tuples problem.
In this case the linear forms involved are not independent and thus an
asymptotic would give very strong information on the regularity of
distribution of the function involved. 
Results of this type lie out of reach of the Green-Tao
Hardy-Littlewood method.

It is worth mentioning at this point a recent related
result of Henriot \cite{henriot}, which provides a correct order upper
bound for $\E_{n \leq N} F(|Q_1(n)|,\dots,|Q_t(n)|)$, where 
$F:\N^t \to \R_{\geq0}$
belongs to a family of functions that does include
$F(n_1,\dots,n_t)=\prod_{i=1}^t R_{f_i}(n_i)$, and where the $Q_i$ are
coprime irreducible polynomials.
The bounds in this result are independent of the discriminant of the
polynomial $Q_1 \dots Q_t$.

Theorem \ref{main thm} has some natural arithmetic consequences.
Analysing the frequency of $4$-term arithmetic progressions in sums of
two squares (weighted by the representation function) may be viewed as a
special case of studying the (average) number of simultaneous zeros of a
pair of diagonal quadratic equations, namely solutions to
\begin{align*}
\begin{array}{lrlll}
 x_1^2 + x_2^2  &- 2x_3^2 -2 x_4^2 + x_5^2+ x_6^2 &&= 0 \cr
&x_3^2 + x_4^2 - 2 x_5^2 - 2 x_6^2 &+ ~x_7^2+ x_8^2 &= 0~.
\end{array}
\end{align*}
While the respective system for $3$-term progressions may easily be
handled by the circle method, Heath-Brown mentions in \cite{heath-brown}
that in order to give an asymptotic for the number of $4$-term arithmetic
progressions in sums of two squares  ``it would appear that one would
require a version of the `Kloosterman refinement' for a double
integral''. 
Browning and Munshi \cite{browning-munshi} have succeeded in
showing that the circle method can in fact be employed to study any pair
of quadratic equations in $n\geq9$ variables that takes the form
$$F_1(x_3,\dots,x_n)=-c(x_1^2 + x_2^2), \quad
F_2(x_3,\dots,x_n)=0~.$$
Previously, the classical Hardy-Littlewood method had been successfully
applied to pairs of diagonal quadratic equations in at least $9$
variables:
\begin{theorem}[Cook \cite{cook}]
 Let $F,G: \Z^9 \to \Z$ be integral diagonal quadratic forms such that
for all real $\lambda, \mu$, not both zero, $\lambda F + \mu G$ is an
indefinite form in at least $5$ variables. Then there is some positive
constant $K_0$ such that the number of simultaneous integral zeros of $F$
and $G$ in the box
$$P \leq x_i \leq CP, \qquad i=1, \dots, 9~,$$
is given by
$$\mathcal N(P) = K_0 P^5 + o(P^5) \text{ as } P \to \infty~.$$
\end{theorem}

Our result, which is in fact an analogue of
\cite[Thm.\,1.8]{green-tao-linearprimes}, considers certain highly
singular systems of quadratic equations in $8$ or more variables.

\begin{theorem}\label{thm-application}
Let $t \geq 4$ and let $f_1, \dots, f_t$ be primitive positive definite
binary quadratic forms. For an integer $s \leq t-2$, let 
$A \in M_{s\times t}(\Z)$ be a full rank matrix whose row-span over
$\mathbb Q$ contains no non-trivial element with less than $3$ non-zero
entries.
 
Define a height function $H: \Z^{2t} \to \R_{\geq0}$ by
$$H(x)= \max_{j\in \{1,\dots,t\}} \sqrt{f_j(x_{2j-1},x_{2j})}~.$$
Then the simultaneous
zeros of the system of quadratic forms
$$F_i(x_1,\dots,x_{2t})= \sum_{j=1}^t a_{i,j} f_j(x_{2j-1},x_{2j})~,
\quad i \in \{1, \dots, s\}~,$$
satisfies the following asymptotic:
\begin{align*}
|\{x \in \Z^{2t}: H(x) \leq N,~F_1(x)= \dots= F_s(x)=0 \}| 
= \frac{(2\pi)^t}{\sqrt{|D_1\dots D_t|}}
   \alpha_{\infty}
   \prod_p \alpha_p + o(N^{2t})~,
\end{align*}
where
$$\alpha_p:= 
\lim_{m\to \infty} \frac{
 |\{x \in (\Z/p^m\Z)^{2t} : F_1(x) \equiv \dots \equiv F_s(x) \equiv 0  
   \Mod{p^m} \}|}{(p^m)^{2t-s}}$$
and 
$$\alpha_{\infty}:= |\{ z \in \{1,\dots,N^2\}^t: Az=0 \}|~.$$
\end{theorem}

We conclude this introduction with the fairly short deduction of Theorem
\ref{thm-application}.
\begin{proof}[Proof of Theorem \ref{thm-application} from Theorem
\ref{main thm}]
The number of simultaneous zeros of bounded height of the forms $F_1,
\dots, F_t$ can be reinterpreted in terms of representation functions:
\begin{align}\label{translation-to-corr}
 |\{x \in \Z^{2t}: F_1(x)= \dots =F_t(x) = 0, 
       H(x) \leq N \}| 
= \sum_{\substack{z\in [N^2]^t:Az=0}} \prod_{j=1}^t r_{f_j}(z_j)~.
\end{align}
To turn the latter expression into the form of a linear
correlation, we may follow \cite[\S4]{green-tao-linearprimes}:
Pick a basis for the integer lattice
$$\Gamma := \{z \in \Z^t: Az=0\}~.$$
Since $A$ has full rank, $\Gamma$ has rank $d:= t - s$, and thus there
are linear forms $\psi_1,\dots,\psi_t: \Z^d \to \Z$ such that 
$$\Gamma = \{(\psi_1(n),\dots,\psi_t(n)): n \in \Z^d\}~.$$
This system of forms has finite complexity, as otherwise we would find
$i\not=j$ such that $\alpha_i \psi_i = \alpha_j \psi_j$ for some non-zero
integers $\alpha_i, \alpha_j$. Hence,
$$\Gamma = \{z \in \Z^t: Az=0, \alpha_iz_i-\alpha_jz_j=0\}~,$$ which
implies by the full rank assumption on $A$ that the row-space of $A$
contains a non-trivial element with less than $3$ non-zero entries, a
contradiction. 

Thus, \eqref{translation-to-corr} takes a form to which Theorem \ref{main
thm} applies and we obtain:
\begin{align*}
& |\{x \in \Z^{2t}: F_1(x)= \dots =F_s(x) = 0, 
       H(x) \leq N \}| \\
&= \sum_{\substack{n \in \N^d \cap \Psi^{-1}([1,N^2]^t)}} 
   \prod_{j=1}^t r_{f_j}(\psi_j(n_j)) \\
&= \vol(\R_{\geq 0}^d \cap \Psi^{-1}([0,N^2]^t))
   \frac{(2\pi)^t}{\sqrt{|D_1\dots D_t|}}
   \prod_p \beta_p + o(N^{2d})~.
\end{align*}
Note that 
\begin{align*}
\vol(\R_{\geq 0}^d \cap \Psi^{-1}([0,N^2]^t)) 
&= |\{n \in \Z^{d} : \Psi(n)\in [0,N^2]^t \}| + o(N^{2d})\\
&= |\{z \in \{1,\dots N^2\}^t : Az=0  \}| + o(N^{2d})~,
\end{align*}
which justifies to define 
$\alpha_{\infty}:= |\{z \in \{1,\dots N^2\}^t : Az=0 \}|$.
It remains to interpret the local factors $\beta_p$ in terms of
$F_1,\dots,F_t$.
If $m$ is sufficiently large, the $\mathbb{\Z}$-basis 
$(\psi_j)_{j \in [d]}$ of $\Gamma$ gives rise to a basis of 
$\{z \in (\Z/p^m\Z)^{t}: Az \equiv 0 \Mod{p^m} \}$, whence
\begin{align*}
& \E_{a \in (\Z/p^m\Z)^d} \prod_{j=1}^t
 \frac{\rho_{f_j,\psi_j(a)}(p^m)}{p^{m}} \\
&= p^{-m(t+d)} \sum_{a \in (\Z/p^m\Z)^d} \prod_{j=1}^t
 |\{(x_{2j-1},x_{2j}) \in [p^m]^2: f_j(x_{2j-1},x_{2j}) \equiv
   \psi_j(a) \Mod{p^m} \}| \\
&= \frac{
 |\{x \in (\Z/p^m\Z)^{2t} : F_1(x) \equiv \dots \equiv F_s(x) \equiv 0  
   \Mod{p^m} \}|}{(p^m)^{2t-s}}~,
\end{align*}
which yields $\beta_p = \alpha_p$ for all primes $p$.
\end{proof}

\section{A majorant for the representation function via the Kronecker
sum}
\label{rep-majorant-section}
\subsection*{Preliminaries and notation}
Recall that a binary quadratic form $f(x,y)=ax^2 + bxy + cy^2$ is
primitive when $(a,b,c)=1$ and that its discriminant is given by
$D(f)=b^2-4ac$. 
Throughout this paper all binary quadratic forms will be assumed to be
positive definite.
The number of ways a form $f$ represents an integer $n$ is described
by the \emph{representation function} $R_f: \Z \to \Z$, defined by 
$$R_f(n) := |\{(x,y): f(x,y)=n\}|~.$$
In order to make use of some multiplicative properties of $R_f$, we
introduce the function $r_f: \Z \to \Z$, defined by 
$$r_f(n) := R_f(n)/k(D)~,$$ where $k(D)$ denotes the number of automorphs
of binary forms of discriminant $D$.
We have $k(D)=6,4,2$ according to $D=-3$, $D=-4$ or $D<-4$, respectively.

Closely related to $r_f$ is the function $r_{D(f)}:\Z \to \Z$ 
which counts---up to the factor $k(D)$---the number of ways $n$ is
represented by any equivalence class of forms of discriminant
$D=D(f)$.
We define $r_{D(f)}$ by
$$r_{D(f)}(n) := \sum_{D(f')=D(f)} r_{f'}(n)~,$$ where $f'$ runs through
a complete system of representatives of primitive forms of discriminant
equal to $D(f)$.

The function $r_{D(f)}$ majorises $r_f$ and has some properties that
suggests it may be a good candidate to start the construction of a
pseudorandom majorant with: on the one hand, the number $h(D)$ of
equivalence classes of primitive forms of discriminant $D$ is finite, and
thus the average order of $r_{D(f)}$ is comparable to the average order of
$r_f$; 
on the other hand, $r_{D(f)}$ has an arithmetic representation as a
divisor sum, a structure that proved to be very suited for the
construction of a pseudorandom majorant in both
\cite{green-tao-longprimeaps} and
\cite{m-divisorfunction}.

Let $f$ be a primitive positive definite form of discriminant $D$. 
Then $r_f(n)=0$ for $n<0$ and $r_f(0)=1$.
For positive integers $n$ coprime to $D$, $r_{D}$ has the representation 
$$r_{D(f)}(n) = \sum_{d|n} \kron{D}{d}$$
as a character sum, where the symbol is a Kronecker symbol.
For general $n$, we pick up another factor which depends only on 
$\gcd(n,D)$ and the parities of the $\alpha$ in
$\prod_{p|D, p^{\alpha}\|n} p^{\alpha}$.
We will see in Corollary \ref{r_D-corollary} that
\begin{align} \label{r_D-bound}
r_{D(f)}(n) \ll_D \sum_{d|n} \kron{D}{d}
\end{align}
holds for all $n \in \N$.
 
Recall that the Kronecker symbol is only non-zero when its entries are
coprime and that furthermore the following lemma holds; see for instance
~\cite[Thm. 1.14]{cox}.

\begin{lemma}
If $D \equiv 0,1 \Mod{4}$ is a non-zero integer, then there is a unique
character $\chi_D:(\Z/D\Z)^* \to \{-1,1\}$ such that
$\chi_D\big(p\Mod{D}\big)=\kron{D}{p}$ for odd $p$ coprime to $D$.
\end{lemma}

Let $\mathcal Q_{D}$ denote the set of primes for which
$\chi_D(p)=-1$. Note that this is the union of the primes in a
collection of progressions modulo $D$. By multiplicativity we have
\begin{equation}\label{kronecker}
 \sum_{d|n} \kron{D}{d} 
= \sum_{d|n} \chi_D(d) 
= \prod_{p^a\| n} (1 + \chi_D(p) + \dots + \chi_D(p^a))
= \tau_D(n)
  \prod_{p^{\alpha}\|n, p \in \mathcal Q_D}
  \frac{1}{2} (1+(-1)^{\alpha})~,
\end{equation}
where $$\tau_D(m)=\prod_{p^a\|m, \chi_D(p)=1} (a+1)~.$$

We denote by $\P_D$ the set of primes for which $\chi_D(p)=1$. Thus,
a square-free number $n$ is represented by some form of discriminant
$D(f)$ only if all of its prime factors belong to $\P_D$ or
divide $D(f)$. 

We can say a little more about the sets $\P_D$ and $\mathcal Q_D$:
Since $\chi$ is a non-principal character taking values $\pm 1$, the
fact that $\sum_{a \in (\Z/D \Z)^*} \chi(a) = 0$  implies that both
$\P_D$ and $\mathcal Q_D$ are the union of the primes in exactly
$\frac{\phi(D)}{2}$ progressions modulo $D$. Thus, the square-free
numbers that are coprime to $D$ and representable by some form of
discriminant $D$ are those numbers whose prime factors belong to a set
comprising asymptotically half the prime numbers.

As a final piece of notation, given any set $\mathcal P$ of primes,
let $\langle \mathcal P \rangle $ denote the set of natural numbers all
of whose prime factors belong to $\mathcal P$.
Thus we may write
$$\tau_D(n) = \sum_{d \in \pdspan } 1_{d|n}~.$$

\subsection*{Construction of the majorant}
The key observation for the construction of our majorant for $r_f$ is that
according to \eqref{r_D-bound} and \eqref{kronecker} it suffices to find
two majorants separately: one for a divisor-type function related to
$\tau_D$, and one for the characteristic function of numbers without
$\Q_D$-prime factors.
Writing $\P_D^* = \P_D \cup \{p:p|D\}$, the characteristic function of
interest is $1_{\langle \P_D^* \rangle}$. 
The shifts by square factors of the form 
$\prod_{p \in \mathcal Q_D} p^{2\alpha}$ only influence the asymptotic
density by a constant factor and may be taken care of separately.
If $\nu$ is a majorant for $\tau_D$ and if $\beta$ is a majorant for
$1_{\langle {\P}_D^* \rangle}$, then $r_f(n)$ is majorised by
$$O_D(1)\nu(n) \sum_{m \in \qdspan} \beta(n/m^2) 1_{m^2|n}~.$$
The majorant $\beta$ for $1_{\starpdspan}$ will be chosen
as a sieving majorant.
In fact, the approach via sieve weights in \cite{green-tao-linearprimes}
proves universal enough to apply here too without much change.
Concerning $\nu$, we make use of the results on the divisor function
from \cite{m-divisorfunction}.

Since neither $\tau_D$ nor $1_{\langle {\P}^*_D \rangle}$ has
asymptotic density, we proceed to determine the average order of $\tau_D$
and show that $\E_{n \leq N} \tau_D (n) \asymp (\log N)^{1/2}$.
This suggests to renormalise the factors in the bound on $r_{D(f)}$
as follows
$$
r_{D(f)}(n) \ll_D
\frac{\tau_D(n)}{(\log N)^{1/2}} 
\sum_{\substack{m \in \qdspan \\ m^2|n }}  
1_{\starpdspan}(n/m^2) (\log N)^{1/2} ~.$$
Iwaniec \cite{iwaniec} proves via sieve theory that it is indeed the case
that $1_{\starpdspan}$ is of average order $(\log N)^{-1/2}$.
This bound, however, is not needed here.
\begin{lemma}\label{order-of-tau_D}
$\tau_D$ satisfies the asymptotic bounds 
$$\E_{n \leq N} \tau_D(n) \asymp (\log N)^{1/2}~,$$
where the implicit constants may depend on $D$.
\end{lemma}

\begin{proof}
We have
\begin{align*}
 \E_{n \leq N} \tau_D(n) 
= \frac{1}{N}\sum_{\substack{d \in \pdspan \\ d \leq N}}
  \left[\frac{N}{d}\right]
= \sum_{\substack{d \in \pdspan\\ d\leq N}} \frac{1}{d} + O(1)~.
\end{align*}
To estimate the last sum, observe that on the one hand
\begin{align*}
 \sum_{\substack{d \in \pdspan \\ d \leq N}} \frac{1}{d}
\leq \prod_{\substack{p \in \P_D \\ p \leq N}} (1- p^{-1})^{-1}
\ll (\log N)^{1/2}
\end{align*}
holds, where the last step follows from the prime number theorem in
arithmetic progressions in the form
 $$\sum_{\substack{p\equiv a \Mod{q}\\p\leq N}}
p^{-1} = \frac{1}{\phi(q)} \log \log(N) + O(1)~.$$
The above remains true when replacing $\P_D$ by 
$\starqdspan:=\Q_D \cup \{p:p|D\}$ and $O(1)$ by $O_D(1)$. 
On the other hand the following chain of inequalities allows us to deduce
a matching lower bound
\begin{align*}
\log N + O(1) 
= \sum_{n \leq N} \frac{1}{n} 
 \leq 
 \Bigg( 
   \sum_{\substack{m_1 \in \pdspan \\ m_1 \leq N}} \frac{1}{m_1}
 \Bigg)
 \Bigg( 
   \sum_{\substack{m_2 \in \starqdspan \\ m_2 \leq N}} \frac{1}{m_2}
 \Bigg)
 \ll 
 \Bigg( 
   \sum_{\substack{m_1 \in \pdspan \\ m_1 \leq N}} \frac{1}{m_1}
 \Bigg)
 (\log N)^{1/2}~.
\end{align*}
\end{proof}

\subsection*{The divisor-type majorant}
To start with, we recall the divisor function majorant that was
constructed in \cite{m-divisorfunction} based on Erd\H{o}s's work
\cite{erdos}.
For any $\gamma > 0$ define the truncated divisor function
$\tau_{\gamma}: [N] \to \Z$ by
$$\tau_{\gamma}(n):= \sum_{d \leq N^{\gamma}} 1_{d|n}$$
and the truncated restricted divisor function 
$\tau_{D,\gamma}: [N] \to \Z$ by
$$\tau_{D, \gamma}(n):= 
\sum_{\substack{d \in \pdspan \\ d \leq N^{\gamma}}} 1_{d|n}~.$$

\begin{proposition}[\cite{m-divisorfunction},
Majorant for the divisor function]
\label{chapter3-divisor-majorant}
Let $\xi = 2^{-m}$ for some $m\in \N$.
Let $C_1 > 1$ be a parameter and write $X_0=X_0(C_1,N)$ for the
exceptional set of all $n \leq N$ satisfying either of
the following
\begin{enumerate}
\item  $n$ is excessively ``rough'' in the sense that it is divisible by
some prime power $p^a$, $a \geq 2$, with $p^a > \log^{C_1} N$, or 
\item $n$ is excessively ``smooth'' in the sense that if $n = \prod_{p}
p^a$ then
$$ \prod_{p \leq N^{1/(\log \log N)^3}} p^a \geq N^{\xi/\log \log N}~.
$$
\end{enumerate}
Further, define $U(i,2/\xi):=\{1\}$ for $i= \log_2 (2/\xi) - 2$, and 
$U(i,2/\xi):= \emptyset$ else.  
If $s > 2/\xi$, write $U(i,s)$ for the set of all products of 
$m_0(i,s):= \lceil \xi s(i + 3 - \log_2 s)/100 \rceil$ distinct primes
from the interval $[N^{1/2^{i+1}}, N^{1/2^i}]$.
Define $\tilde \nu_{\xi} : [N] \rightarrow \R_+$ by 
\[ 
 \tilde \nu_{\xi}(n) :=  
 \sum_{s \geq 2/\xi}^{(\log \log N)^3}
 \sum_{i \geq \log_2 s - 2}^{6 \log \log \log N}
 \sum_{u \in U(i,s)} 
 2^s 1_{u|n} \tau_{\xi}(n) +
 1_{n \in X_0} \tau(n)~.
\] 
Then $\tau(n) \leq \tilde \nu_{\xi}(n)$ for all $n\leq N$, provided $N$
is large enough.
\end{proposition}
Note that the main term of $\tilde\nu_{\xi}$ has low complexity in that it
only involves small divisors since all $u \leq N^{\xi}$.
Restricting all occurrences of divisor functions in $\tilde \nu_{\xi}$ to
only count divisors in $\pdspan$, yields a majorant for $\tau_D$ of the
same order of magnitude as $\tau_D$.
We make one further modification and replace the cut-off in the
definition of $\tau_{D,\gamma}$ by a smooth cut-off of the form which
appears in Green and Tao's $\Lambda$-majorant.
This turns out to be advantageous when establishing the linear forms
condition.
Thus, let $\chi:\R \to \R_{\geq 0}$ be a smooth, even function that is
supported on $[-1,1]$ and satisfies the properties $\chi(x)=1$ for 
$x \in [-1/2,1/2]$ and $\int_0^1 |{\chi}'|^2~dx = 1$.
Define $\tau^*_{D,\gamma}: [N] \to \Z$ by
$$\tau^*_{D, \gamma}(n):= 
\sum_{\substack{d \in \pdspan \\ d \leq N^{\gamma}}} 1_{d|n} 
\chi \left( \frac{\log d}{\log N^{\gamma}} \right)~.$$
Then $\tau_{D,\gamma/2} \leq \tau^*_{D, \gamma}(n) 
\leq \tau_{D, \gamma}(n)$ holds.
With this definition we have the following lemma.

\begin{lemma}[A majorant for $\tau_D$]
Let the sets $U(i,s)$ be those which Proposition
\ref{chapter3-divisor-majorant} produces for $\xi=\gamma/2$.
Let $\nu_{D,\gamma}: [N] \to \R$ be defined by
\begin{align*}
C \nu_{D,\gamma} (n) :=  
 \frac{1}{\sqrt{\log N}} \Bigg( 
 \sum_{s \geq 4/\gamma}^{(\log \log N)^3}
 \sum_{i \geq \log_2 s - 2}^{6 \log \log \log N}
 \sum_{u \in U(i,s)}
 2^s
 1_{u|n} 
 \tau^*_{D,\gamma}(n) 
 + 1_{n \in X_0} \tau_{D}(n) \Bigg)~.
\end{align*}
Then $\tau_{D}(n)/(\log N)^{1/2} \leq C \nu_{D,\gamma} (n)$ for all $n
\in [N]$ and there is some constant $C$ bounded independently of $N$ such
that $\E_{n \leq N} \nu_{D,\gamma} (n) = 1 + o(1)$. 
\end{lemma}

\begin{proof}
 We begin by checking the majorisation property. 
For any $n \in [N]$, write $n=n_{\P} m$ where $n_{\P}$ is the largest
factor of $n$ that belongs to $\pdspan$. Then 
$$\tau_{D}(n) = \tau(n_{\P}) \leq \tilde \nu_{\gamma/2} (n_{\P}) 
\leq C (\log N)^{1/2} \nu_{D,\gamma}(n_{\P}) 
= C (\log N)^{1/2} \nu_{D,\gamma}(n)~,$$ as required.
The existence of $C$ follows as in the proof of 
\cite[Prop. 4.2]{m-divisorfunction}, taking into account that 
$\E_{m\leq N} \tau_{D,\gamma/2}(m) \asymp (\log N)^{1/2}$, which is
proved in much the same way as Lemma \ref{order-of-tau_D}.  
\end{proof}

\subsection*{The sieving type majorant}
The next task is to give a majorant $\beta : \N \to \R^+$ for the
characteristic function of the set $\starpdspan$ of
numbers without $\Q_D$-prime factors. 
Adapting the Selberg-sieve majorant for primes from
\cite{green-tao-linearprimes} to the set
$1_{\langle \P^*_D \rangle}$, we aim to remove all integers that have
a prime factors $p$ from $\Q_D$ with $p \leq N^{\gamma}$. 
Let $\chi:\R \to \R$ be a smooth, even function that is supported on
$[-1,1]$ and satisfies the properties $\chi(0)=1$ and
$\int_0^1 |{\chi}'|^2~dx = 1$.
Define in analogy to \cite[App. D]{green-tao-linearprimes}
\begin{equation*}
\beta(n) 
:= \Lambda_{D,\chi, \gamma}(n) 
:= C' (\log N)^{1/2} 
 \left(
 \sum_{d|n, d \in \qdspan} 
 \mu(d)~ 
 \chi \left( \frac{\log d}{\log N^{\gamma}} \right)
 \right)^2~, 
\end{equation*}
for some constant $C'$.
The results from \cite{green-tao-linearprimes} show that $C'$ may be
chosen such that $\E_{n\leq N} \beta(n) = 1 + o(1)$.
This will play a role in Section
\ref{linear-forms-and-correlation-conditions}.
Note that $\beta(m)=C'(\log N)^{1/2}$ at every $\Q_D$-prime-free integer
$m \leq N$, and thus we have the pointwise majorisation
$$1_{\starpdspan}(n) (\log N)^{1/2} \leq C'^{-1} \beta(n),
\quad n \in [N]~.$$

\section{A reduction of the main theorem} \label{bar-reduction-section}

While it is possible to apply the nilpotent Hardy-Littlewood method to
the representation function $r_f$ itself, it is the aim of this section to
show that we can deduce the main theorem from a similar statement about a
smoothed version of $r_f$, that is, a function that agrees with $r_{f}$
everywhere except on a sparse set where the restricted divisor
function $\tau_D$ shows exceptionally irregular behaviour.

First note that the pointwise bound 
$r_f(n) \leq r_{D(f)}(n) \leq \tau(n)$ for $n \in \N$ of the
representation function of any primitive positive definite quadratic
form $f$ by the divisor function gives the following second moment
estimate.

\begin{lemma}[Second moment estimate]
Let
${f_1}, \dots, {f_t}$ be primitive positive definite binary quadratic
forms and let $\Psi=(\psi_1, \dots, \psi_t): \Z^m \to \Z^t$ be a system of
affine-linear forms whose linear coefficients are bounded by $L$.
If $K \subset [-N,N]^d$ is a convex body such that $\Psi(K) \subseteq
[0,N]^t$, then
$$
\E_{n \in \Z^m \cap K} \prod_{i \in [t]} r_{f_i}^2(\psi_i(n)) 
\ll_{t,m,L} (\log N)^{O_t(1)}~.
$$
\end{lemma}
\begin{proof}
Let $K':=\{x \in K: \Psi(x) \in [1,N]^t\}$.
Then H{\"o}lder's inequality yields
\begin{align*}
\E_{n \in \Z^m \cap K} \prod_{i \in [t]} r_{f_i}^2(\psi_i(n)) 
&\leq \prod_{i \in [t]} \left( 1 + 
\E_{n \in \Z^m \cap K'} r_{f_i}^{2t}(\psi_i(n))
\right)^{1/t} \\
&\leq \prod_{i \in [t]} \left( 1 + 
\E_{n \in \Z^m \cap K'} \tau^{2t}(\psi_i(n))
\right)^{1/t}
\end{align*}
The remaining steps are standard; cf.~the proof of 
\cite[Lemma 3.1]{m-divisorfunction} for details.
\end{proof}

The next lemma, which is a combination of some technical lemmas from
\cite{erdos}, describes an exceptional set for the divisor function,
i.e.~a sparse set containing those numbers on which the divisor function
behaves particularly irregularly.

\begin{lemma}\label{exceptional-set}
Let $C_1>1$ be a parameter and write $X_0$ for the set of all positive
$n \leq N$ satisfying either of the following
\begin{enumerate}
\item $n$ is excessively ``rough'' in the sense that it is divisible by
some prime power $p^a$, $a \geq 2$, with $p^a > \log^{C_1} N$, or 
\item $n$ is excessively ``smooth'' in the sense that if $n = \prod_{p}
p^a$ then
\[ \prod_{p \leq N^{1/(\log \log N)^3}} p^a \geq N^{\gamma/\log \log N},\]
\item $n$ has a large square divisor $m^2|n$, $m > N^{\gamma}$.
\end{enumerate}
Then
$$\E_{n \in K\cap \Z^d} \sum_{i=1}^t 1_{\psi_i(n) \in X_0} 
\ll \log^{-C_1/2}N~.$$
\end{lemma}
\begin{proof}
See \cite{erdos} for the original results or \cite[\S3]{m-divisorfunction}
for their adaptation to this situation.
\end{proof}
The previous two lemmas allow us to deduce the main theorem from an
equivalent statement about smoothed versions of the representation
functions $r_{f_i}$.
The particular smoothed functions we shall work with will be chosen in
Section \ref{W-trick}.
\begin{lemma}\label{bar-reduction}
Let ${f_1}, \dots, {f_t}$ be primitive positive definite binary quadratic
forms.
For each $i \in [t]$, let $\bar r_{f_i}:\{ 0,\dots,N\} \to \R$ denote a
function that agrees with $r_{f_i}$ on $[N] \setminus X_0$, that is,
outside the exceptional set of the divisor function, and which further
satisfies
$0\leq \bar r_{f_i}(n) \leq r_{f_i}(n)$ for all $n \in X_0 \cup \{0\}$.
If the parameter $C_1$ of the exceptional set is sufficiently large, then
the main theorem holds if and only if under the same conditions
$$ 
\sum_{n\in \Z^d\cap K } 
\bar r_{f_1}(\psi_1(n)) \dots
\bar r_{f_t}(\psi_t(n)) 
 = \beta_{\infty} \prod_{p} \beta_p + o(N^d) ~.
$$
\end{lemma}
\begin{proof}
This follows by the Cauchy-Schwarz inequality from the previous two
lemmas and the bound
$$\sum_{n \in K \cap \Z^d} \sum_{i=1}^t 1_{\psi_i(n)=0} 
\ll N^{d-1}~.$$
\end{proof}

The above lemma in particular shows that a pseudorandom majorant used in
a proof only needs to majorise the function $r_f$ (or $\bar r_f$) on the
set of positive unexceptional integers. 
We can therefore truncate the summation over dilates of 
$m^2$, $m \in \qdspan$, in the majorant to those $m$ with $m< N^{\gamma}$.
Furthermore, we may restrict attention to the case where 
$\Psi(K) \subseteq [1,N]^t$

\section[Distribution in residue classes]{Distribution in residue classes
}\label{reduction-to-genus-class}

The transference principle from 
\cite{green-tao-longprimeaps, green-tao-linearprimes}, which we shall
employ later, only works with functions $h$ that are sufficiently
quasirandom in the sense that all $U^k$-norms $\| h - \E h \|_{U^k}$ up
to some order $k$, determined by the specific system $\Psi$ one is
working with, are small.
A necessary condition for the uniformity norms to be small, is that the
function $h$ at hand is equidistributed in residue classes to small
moduli. 
This condition is in fact equivalent to requiring that $h$ does not
correlate with periodic nilsequences of short period, cf.~Section
\ref{non-corr-intro}.

As seen at the start of the introduction, the representation function
$r_f$ does not have this property. 
To remove these obstructions to uniformity, one can try to split the
function $r_f$ into a sum of functions each of which does not detect a
difference between residue classes to small moduli.
This strategy is known as $W$-trick.
In order find a suitable decomposition, we shall investigate the
quantities
$$\E_{n \leq N} 1_{n \equiv \beta \Mod{q}}~ r_f(n)$$
for fixed period $q$ and fixed residue class $\beta$.
Define
$$\rho_{f,\beta}(q):=|\{(x,y) \in [q]^2: f(x,y) \equiv \beta \Mod{q}\}|$$
to be the number of representations of $\beta\Mod{q}$, 
and let $K(N) = f^{-1}([0,N]) \subseteq \R^2$. This is the area
enclosed by the ellipse $f(x,y)=N$ and hence a convex set of volume
$$\vol K(N) = \frac{2\pi N}{\sqrt{-D}}~.$$

A volume packing argument,
cf.~\cite[App. A]{green-tao-linearprimes}, yields
$$ \sum_{\substack{n \leq N\\ n \equiv \beta \Mod{q}}} r_{f}(n) k(D)
=\sum_{\substack{(x,y)\in K(N)\cap \Z^2\\f(x,y)\equiv \beta\Mod{q}}} 1 
=\frac{\rho_{f,\beta}(q)}{q^2}\vol(K)+ O(\sqrt{N}q)~,$$
which proves the following lemma.
\begin{lemma} \label{distribution_in_APs}
Let $P:=\{n\leq N:n \equiv \beta\Mod{q}\}$ be an arithmetic progression.
Then the average of the representation function of $f$ along $P$ satisfies
\begin{align*} 
 \E_{n \in P} r_{f}(n)
  &= \frac{2\pi}{k(D) \sqrt{-D}} \frac{\rho_{f,\beta}(q)}{q} 
    + O(|P|^{-1/2} q^{2})~.
\end{align*}
\end{lemma}

In view of this lemma it is not surprising that we will make use of
several further observations on the densities
$\rho_{f,\beta}(q)q^{-1}$, 
which will be established in Section \ref{representation_mod_q}.

\begin{lemma}\label{ind-of-g}
$\rho_{f,\beta}(q)$ only depends on the genus class of $f$.
\end{lemma}
\begin{proof}
 Two forms $f_1$ and $f_2$ belong to the same genus if and only if
they are locally equivalent in the following sense: 
for every non-zero integer $m$ there exists $\sigma_m \in Gl_2(\Z/m\Z)$
such that
$$f_1(x,y)\equiv f_2((x,y)\sigma_m)\mod{m}~.$$ 
Thus, $\rho_{f_1,\beta}(q)=\rho_{f_2,\beta}(q)$ for all positive
integers $q$ and all $\beta \in [q]$.
\end{proof}
The reason this lemma is important to us is that it allows us to consider
instead of $r_f$ the following more regularly behaved function in all
questions regarding the distribution in residue classes.
Let the genus class representation function $r_{g}:\N \to \N$ be
defined by
$$r_g(n) = \E(r_f(n) | f \in g)~,$$
where $f$ runs through a system of representatives of classes in the
genus $g$.
Under the assumptions of Lemma \ref{distribution_in_APs} we then have 
\begin{align}\label{genus-class-fn-in-APs}
\E_{n \in P} r_{f}(n) 
= \E_{n \in P} r_{g(f)}(n) + O(|P|^{-1/2} q^{2})~, 
\end{align}
where $g(f)$ denotes the genus that contains $f$.

\section{Results from the theory of binary quadratic forms}
\label{appendix-on-forms}

The aim of this section is to prove the bound \eqref{r_D-bound} on the
number of representations of a positive integer $n$ by a form of
discriminant $D$, which was used to construct the majorant function in
Section \ref{rep-majorant-section}.

\subsection{Representation by primitive forms of fixed discriminant}

The question of whether or not $m$ is properly representable by a
primitive form of discriminant $D$ is linked to the solubility in $x$ of
the congruence
\begin{equation}\label{criterion}
x^2 \equiv D \Mod{4m}~;
\end{equation}
see \cite[p.506]{IK} or \cite[p.172]{rose}.
If $f$ is a form of discriminant $D$ that represents $m$ properly, then
$f$ is equivalent to $\< m,n,*\> = mX^2 + nXY + *Y^2$, where $D=n^2 - 4mk$
for some integer $k$. 

\begin{claim}
Consider the solutions $x=n$ to \eqref{criterion} that satisfy
$0<n\leq2m$. These form a complete set of incongruent solutions
modulo $2m$.
Those solutions among them for which $\< m,n,(D-n^2)/4m \>$ is primitive
are in one-to-one correspondence with the distinct classes of primitive
forms that represent $m$ properly.
\end{claim}

\begin{proof}
Let $f$ be a primitive form and suppose there are coprime $u$ and $v$ such
that $f(u,v)=m$.
Choose a solution $(z_0,w_0)$ to $1=uz_0-vw_0$.
Then 
$$f'(X,Y)
:= f ( (X,Y)
\left(\begin{array}{cc}
u & v \cr
w_0 & z_0
\end{array}
\right) )
= mX^2 + n XY + \frac{D-n^2}{4m} Y^2
$$
is an equivalent form with leading coefficient $m$.
Choosing different solutions $w = w_0 + w'$ and $z= z_0 + z'$ to 
$1 = uz-vw$, we have $w'=tu$ and $z'=tv$ for some non-zero
integer $t$, which implies that the middle coefficient $n$ is
unique modulo $2m$.
In particular $\<m,n_1,*\> \sim \<m,n_2,*\>$ if and only if 
$n_1 \equiv n_2 \Mod{2m}$.

Observe that in the other direction every solution $x=n$ to 
$x^2 \equiv D \Mod{4m}$ yields an equivalence class $\<m,n,*\>$ of forms
of discriminant $D$ that represents $m$ properly.
\end{proof}

In order to determine the number of classes of forms that represent $m$,
we are interested in two pieces of information:
\begin{enumerate}
 \item[1.] the number of solutions $x$ to $x^2 \equiv D \Mod{4m}$, and
 \item[2.] how many of these solutions yield \emph{primitive} forms
$(m,x,*)$ of discriminant $D$.
\end{enumerate}
A third necessary piece of information regards the number of proper
representations by a fixed class of forms: any two proper
representations of $m$ by a fixed form $f$ are related by an automorph. 
Thus each class $C(f)$ of forms equivalent to $f$ represents $m$ properly
in $k(D)$ different ways, where $k(D)$ is the number of automorphs of
forms of discriminant $D$.

In order to analyse the the number of solutions to \eqref{criterion}, we
introduce the related irreducible quadratic polynomial $P(x)=x^2 - D$,
which has discriminant $4D$.

Let $\rho(a):=|\{k \in [a]: P(k) \equiv 0 \Mod{a} \}|$ denote the number
of zeros modulo $a$. 
The counting function $\rho$ is multiplicative by the Chinese remainder
theorem, which leaves us to determine $\rho$ at prime powers.
If $p \nmid D$, then (cf.~\cite[Thm 12.3.4]{hua})
\begin{align*}
 \rho(p^{\alpha}) = 
\left\{
\begin{array}{cl}
2              & \text{ if } p=2, {\alpha}=2 \cr
2(1+\chi_D(p)) & \text{ if } p=2, {\alpha}>2 \cr
1+\chi_D(p)    & \text{ if } p>2 ~.
\end{array}
\right.
\end{align*}
In the remaining case of primes $p|D$,
Hensel's lemma implies that
$$\rho(p^{\alpha})=\rho(p^{v_p(4D)+1}) \quad \text{ if }
\alpha> v_p(4D)~.$$
For $p|D$ we will show below that, in fact, there are no primitive forms
that properly represent an integer $m$ with $v_p(m) > v_p(D)$ for some
prime $p$.

If $m$ is coprime to $D$, then each solution to $x^2 \equiv D \Mod{4m}$
yields a primitive form, and $\rho(4m)$ is directly linked to the number
$r^*_{D}(m)$ of classes of primitive forms that represent $m$
\emph{properly}: $r^*_{D}(m) = \frac{1}{2} \rho(4m)$.

We turn to the case where $\gcd(D,m)>1$.
If there is a prime $p$ dividing $\gcd(D,4m)$ to an odd power, then
solutions to $D=n^2-4mk$ yield primitive forms if and only if each such
$p$ divides both $D$ and $4m$ to the same power.

Considering the set of forms arising from solutions to \eqref{criterion},
we can, if $\gcd(m,4D)>1$, retrieve the number of primitive
forms among them via an inclusion-exclusion argument. 
Indeed, when $d=\gcd(m,n,k)$, then $m/d$ is properly represented by the
form $\<\frac{m}{d},\frac{n}{d},\frac{k}{d}\>$ of
discriminant $Dd^{-2}$.
Note that automorphs of forms of the first kind are also automorphs of
forms of the second kind and vice versa.

Let $p^{\alpha}\|4m$ and suppose that $p^{\sigma}\|D$, $\sigma>1$.

We begin by analysing the largest range for 
$\alpha$, $\alpha > \sigma >1$.
When $\sigma$ is odd, then there are, as seen above, no primitive forms
that represent $p^{\alpha}$ properly.
Suppose next that $\sigma$ is even and define
$$\rho'(p^{\alpha})
:=|\{x: x^2 \equiv D \Mod{p^{\alpha}} \}| 
 -|\{x: x^2 \equiv Dp^{-2} \Mod{p^{\alpha-1}}\}|~.$$
This quantity counts the number of solutions to 
$x^2 \equiv D \Mod{p^{\alpha}}$ for which 
$x^2 = D + k p^{\alpha}$ for some $k$ not divisible by $p$. 
The expression for $\rho'(p^{\alpha})$ simplifies to 
\begin{align*}
 \rho'(p^{\alpha}) 
 = |\{x: x^2 \equiv Dp^{-\sigma} \Mod{p^{\alpha-\sigma}} \}| 
  -|\{x: x^2 \equiv Dp^{-\sigma} \Mod{p^{\alpha-1-(\sigma-2)}}\}|~,
\end{align*}
which is seen to be $0$ by Hensel's lemma 
(note that $p \nmid Dp^{-\sigma}$).
Thus, no power $p^{\alpha} \nmid D$ of a discriminant-prime with $p^2|D$
is properly representable by a primitive form.

What remains are even powers $p^{\alpha}|D$, $\alpha < \sigma$ and the
case $p^{\alpha}\|D$.
In the former case, any solution to $D=n^2+4mk$ with $p^{\alpha/2}\|n$
satisfies $p \nmid k$. 
Hence there are $p^{\alpha/2}(1-p^{-1})$ choices for $n \Mod{p^{\alpha}}$.
In the latter case, $p \nmid k$ holds if and only if
$p^{\lceil \alpha/2 \rceil}|n$, hence there are 
$p^{\lfloor \alpha/2 \rfloor}$ choices in this case.

In total, the number $r^*(m)$ of primitive forms properly representing $m$
is given by
$$\frac{1}{2}(1+1_{2\nmid D})
\prod_{\substack{p|m, p \nmid D}} (1 + \chi_D(p))
\prod_{\substack{q|D \\ q^{\alpha} \| 4m,\\ q^{\sigma} \| 4D}}
\Big(q^{\alpha/2}(1-q^{-1}) 1_{\alpha < \sigma} 1_{\alpha \text{ even}}
 + q^{\lfloor \alpha/2 \rfloor} 1_{\alpha=\sigma}\Big)~,
$$
where $p$ and $q$ run over primes, and where the factor $1/2$ takes
account of the fact that for every solution $x \in [4m]$, $x+2m$ is the
unique other solution determining the same class of forms.

Collecting everything together, we obtain the following explicit
expression for $r_{D}$:
\begin{corollary}\label{r_D-corollary}
The total number of representations (proper and improper ones) of an
integer $m$ by classes of primitive forms of discriminant $D$ satisfies
\begin{align*}
&r_{D}(m) =
 \sum_{\delta^2|m} r^*_D(m/\delta^2) \\
&=
 \frac{1+1_{2\nmid D}}{2} \sum_{\substack{\delta^2|m \\ (\delta,D)=1}}
  \prod_{\substack{p \nmid D \\ p^{\alpha} \| m \delta^{-2} }} 
    (1 + \chi_D(p))
  \prod_{\substack{q | (D,m) \\ q^{\alpha}\| 4m \\ q^{\sigma} \| 4D }}
 \bigg(
 q^{\lfloor \min(\alpha, \sigma-1)/2 \rfloor}
 1_{\alpha \, \text{\em even}} 
 + 
 q^{\lfloor \sigma/2 \rfloor}
 1_{\alpha \equiv \sigma \Mod{2}}
 1_{\alpha\geq\sigma}
 \bigg) \\
&\ll \sqrt{D} \sum_{\substack{d|m}} \chi_D(d)~,
\end{align*}
where $p,q$ run over primes. 
\end{corollary}

\subsection{Representation by genera}
Recall that the representation function $r_g: \N \to \N$ of a
genus class $g$ was defined to be
$r_g(n) = \E (r_{f'}(n)\mid f' \in g)$,
where $f'$ runs through a system of representatives.
This function is of interest since by Lemma \ref{distribution_in_APs} and
Lemma \ref{ind-of-g}, it has the same distribution in residue classes as
any function $r_f$ with $f \in g$.
We aim to reduce the problem of determining the number of representations
of an integer $n$ by a specific genus class to that of counting
certain representations of the factor $n'$ of $n$ that is coprime to $D$.

This is advantageous for the following reason. 
The values in $(\Z/D\Z)^*$ that are represented by a form $f$ with
$D(f)=D$ form a \emph{coset} of the subgroup in $(\Z/D\Z)^*$ that is
generated by the values the principal form represents, 
c.f.~\cite[Lemma 2.24]{cox}. 
Thus, different genera represent disjoint sets of values in $(\Z/D\Z)^*$.
This means that the character sum expression of the function
$r_{D(f)}$ which counts representations of all classes in $h(D)$ yields
an arithmetic expression for the function $r_{g}$ which just considers
those classes of genus $g$.
Indeed, let $\mathcal R_g$ denote the non-zero residues modulo $D$ that
are represented by forms in $g$.
Then for $n'$ \emph{coprime} to $D$ we have
$$
  r_{g}(n')
= \frac{1}{|g|} 
  \sum_{b \in \mathcal R_g} 
  1_{n' \equiv b \Mod{D}} 
  \sum_{d|n'} \chi_D(d)~.
$$

For an arbitrary positive integer $n$, let $n= n_D \tilde n^2 n'$ be the
factorisation for which $n'$ is coprime to $D$ and $n_D$ is the
largest divisor $n_D|(n,D)$ such that $\frac{n}{n' n_D}=\tilde n^2$ is a
square. 
This factorisation is chosen in such a way that Corollary
\ref{r_D-corollary} implies $r_{D}(n)=r_{D}(n'n_D)$, which is of interest
because in $n'n_D$ the factor that is not coprime to $D$ is bounded.

Let $\< n'n_D,b,c \>$ be a primitive form properly representing $n'n_D$. 
Then, since $n_D | D$, we have $(n_D,b)>1$ and hence $(n_D,c)=1$ by
primitiveness of the form. Since further $(n_D,n')=1$, we have
\begin{align*}
\< n'n_D,b,c\>
&\simeq \<c,-b,n_Dn'\>\\
&\simeq \<c,-b,n_Dn'\>*\<n_D,-b,n'c\>*\<n'c,-b,n_D\> \\
&\simeq \<c n_D,-b,n'\>*\<n'c,-b,n_D\> \\
&\simeq \<n',b,c n_D\>*\<n_D,b,n' c\>~.
\end{align*}
Note that all forms involved are primitive.

Thus, we can decompose the representation into separate ones for the
coprime factors $n'$ and $n_D$. 
We aim to use this multiplicative property of representation by primitive
forms of fixed discriminant in conjunction with the following lemma.

\begin{lemma} \label{cosets}
 The principle genus $\mathcal{G}_0$ is a subgroup of the class group (a
finite and Abelian group). The genera form cosets of $\mathcal{G}_0$ in
the class group.
\end{lemma}

\begin{proof}
 See e.g.~\cite[p.197, Thm 2.8]{rose}.
\end{proof}

With the help of this lemma we have
$$r_{g}(n) = |g| \sum_{g'} r_{g*g'^{-1}}(n_D) r_{g'}(n')~.$$
If the residue $n' \Mod{D}$ is representable by a form of discriminant
$D$, then let $g_{n'}$ denote the unique genus class that represents 
$n' \Mod{D}$.
We may use the arithmetic representation of $r_{g_{n'}}$ to obtain
the following lemma.
\begin{lemma}\label{general-r_g}
 Given $n = n_D n' \tilde n^2$ as above and a genus class $g$, then
$$r_{g}(n) = r_{g*g_{n'}^{-1}}(n_D) \sum_{d|n'} \chi_D(d)~.$$
\end{lemma}

\section{Representation in $\Z/q\Z$} \label{representation_mod_q}
This section contains several results on the densities
$\rho_{f_i,\beta}(p^{\alpha})p^{-\alpha}$, which will be established using
results from the previous section and the following proposition.

\begin{proposition}
\label{r'_D(f)-major_arc_estimate}
Let $P = \{q_0 m + \beta_0: m \leq M \}$ be a progression
such that $D|q_0$ and $\beta_0 \not\equiv 0 \Mod{p^{\alpha}}$ for any
$p^{\alpha}\|q_0$. 
Then
\begin{align*}
  \E_{n \in P} \sum_{d|n} \chi_D(d)
&= C \prod_{p|q_0} (1-\chi_D(p)p^{-1}) 
   \sum_{\alpha \geq 0} 1_{p^{\alpha}|\beta_0} \chi_D(p^{\alpha}) +
   O\Big(\frac{q_0^{1/2-\eps}}{M^{1/2}}\Big),
\end{align*}
where 
$C = (1 + \chi_D(\frac{\beta_0}{(\beta_0,q_0)})) L(1,\chi_D) = O(1)$.
\end{proposition}

We defer the proof to the end of the section.
The following lemma is a rather immediate consequence.

\begin{lemma}\label{rho-approximate}
Let $q$ be a positive integer that is divisible by $D$ and 
let $\beta \in [q]$ be such that $\beta \not\equiv 0 \mod{p^{\alpha}}$
for any $p^{\alpha}\|q$. Then
\begin{align*}
 \frac{\rho_{f,\beta}(q)}{q} 
=  C'
   \prod_{p|q} (1-\chi_D(p)p^{-1}) 
   \sum_{\alpha \geq 0} 1_{p^{\alpha}|\beta} \chi_D(p^{\alpha})~,
\end{align*}
where 
$C' 
= r_{g*g_{\beta'}^{-1}}(\beta_D)
  (1+\chi_D(\frac{\beta}{(\beta,q)})) h(D)
= O(1)$.
\end{lemma}
\begin{proof}
By Lemma \ref{distribution_in_APs} and Lemma \ref{general-r_g} we have
for 
$P(M)=\{m \equiv \beta \Mod{q}: m \leq M\}$
\begin{align*}
 \frac{\rho_{f,\beta}(q)}{q} \frac{2 \pi}{k(D) \sqrt{-D}}
&= \lim_{M\to \infty}\E_{n \in P(M)} r_f(n) 
= \lim_{M\to \infty}\E_{n \in P(M)} r_{g(f)}(n) \\
&= r_{g*g_{\beta'}^{-1}}(\beta_D) 
   \lim_{M\to \infty}\E_{n \in P(M)}\sum_{d|n} \chi_D(d)~.
\end{align*}
By Proposition \ref{r'_D(f)-major_arc_estimate} and the class
number formula the result follows.
\end{proof}

With the help of the previous lemma and a result of Stewart
\cite{stewart}, we obtain the following more explicit information on the
densities $\rho_{f_i,\beta}(p^{\alpha})p^{-\alpha}$.
\begin{lemma} \label{rho-bounds}
\begin{enumerate}[\upshape(a)]
 \item Let $p_0$ be a prime that divides $D$ and suppose that
$\beta \not\equiv 0 \Mod{p_0^{\alpha}}$. Then
 $$\rho_{f,\beta}(p_0^{\alpha})p_0^{-\alpha} = O(1)$$
as $\beta$ and $\alpha$ vary. 
If $\alpha \geq v_{p_0}(D)$ and $\beta \not\equiv 0 \Mod{p_0^{\alpha}}$,
then
$$\rho_{f,\beta}(p_0^{\alpha})p_0^{-\alpha} =
\rho_{f,\beta + kp_0^{\alpha}}(p_0^{\alpha+1})p_0^{-(\alpha+1)}$$
for any $k \in \Z/p_0\Z$
 \item If $p_0\nmid D$ then we have for $\beta \not\equiv 0
\Mod{p_0^{\alpha}}$ 
$$\rho_{f,\beta}(p_0^{\alpha})p_0^{-\alpha} 
= (1-\chi_D(p_0)p_0^{-1})
\sum_{j \geq 0} 1_{p_0^{j}|\beta} \chi_D(p_0^{j})~.$$
 \item Let $p$ be any prime. Then
$$\rho_{f,0}(p^{\alpha}) p^{-\alpha} \ll \alpha$$
holds.
\end{enumerate}
\end{lemma}
\begin{proof}
(a) 
We may assume $\alpha > v_{p_0}(D)$. 
Let $\beta_0 \in (\Z/D\Z)^*$ be a residue representable by $f$ and let
$\beta_1$ be such that 
$\beta_1 \equiv \beta \Mod{p_0^\alpha}$ and
$\beta_1 \equiv \beta_0 \Mod{p^{v_{p}(D)}}$ for any prime
divisor $p \not= p_0$ of $D$.
By choice of $\beta_0$ we have $\rho_{f,\beta_1}(p^{v_p(D)}) \geq 1$
for $p \not= p_0$.
The previous lemma yields
\begin{align*}
 \frac{\rho_{f,\beta_1}
    (p_0^{\alpha} \prod_{p|D, p \not= p_0} p^{v_{p}(D)})}
  {p_0^{\alpha} \prod_{p|D, p \not= p_0} p^{v_{p}(D)}}
 = O(1)~,
\end{align*}
whence the first part of (a) follows by multiplicativity of $\rho$.
Define $\beta_2 \in \Z$ to be such that 
$\beta_2 \equiv \beta_0 \Mod{p^{v_p(D)}}$ for any prime divisor 
$p \neq p_0$ of $D$ and 
$\beta_2 \equiv \beta + k p_0^{\alpha} \Mod{p_0^{\alpha + 1}}$.
Then, by Lemma \ref{rho-approximate},
\begin{align*}
 \frac{\rho_{f,\beta}(p_0^{\alpha})}{p_0^{\alpha}}
 &= \frac{\rho_{f,\beta_2}(p_0^{\alpha})}{p_0^{\alpha}}
 = \frac{\rho_{f,\beta_2}
    (p_0^{\alpha} \prod_{p|D, p \not= p_0} p^{v_{p}(D)})}
  {p_0^{\alpha} \prod_{p|D, p \not= p_0} p^{v_{p}(D)}}
  \prod_{p|D, p \not= p_0} 
  \bigg( 
  \frac{\rho_{f,\beta_0}(p^{v_{p}(D)})}{p^{v_{p}(D)}}
  \bigg)^{-1} \\
 &= r_{g*g_{\beta_2'}^{-1}}((\beta_2)_D)
  (1+\chi_D(\beta_2')) h(D)
  \prod_{p|D, p \not= p_0} \bigg( 
  \frac{\rho_{f,\beta_0}(p^{v_{p}(D)})}{p^{v_{p}(D)}}
  \bigg)^{-1} \\
 &= \frac{\rho_{f,\beta_2}(p_0^{\alpha+1})}{p_0^{\alpha+1}}
  = \frac{\rho_{f,\beta + kp_0^{\alpha}}(p_0^{\alpha+1})}
         {p_0^{\alpha+1}}~.
\end{align*}

The proof of part (b) is almost identical.
Let $\beta_0 \in (\Z/D\Z)^*$ be a residue representable by $f$ and let
$\beta_1$ be such that 
$\beta_1 \equiv \beta \Mod{p_0^\alpha}$ and
$\beta_1 \equiv \beta_0 \Mod{p^{v_{p}(D)}}$ for any prime $p|D$. 
Then $(\beta_1,D)=1$ and $g(f)$ is the unique genus class
representing $\beta_1 \Mod{D}$.
Hence $r_{g(f)*g_{\beta_1'}^{-1}}(1) > 1$, as the principal genus
represents $1$. 
Since $\beta_1$ is representable by $f$, there is some $m$ such that 
$\sum_{d|mD + \beta_1} \chi_D(d) > 0$, hence, in particular 
$\chi_D(mD + \beta_1) = \chi_D(\beta_1) = 1$.
Two applications of Lemma \ref{rho-approximate} yield
\begin{align*}
 \frac{\rho_{f,\beta}(p_0^{\alpha})}{p_0^{\alpha}}
=  \frac{\rho_{f,\beta_1}(p_0^{\alpha})}{p_0^{\alpha}}
=  \frac{\rho_{f,\beta_1}(p_0^{\alpha}D)}{p_0^{\alpha}D}
   \bigg(\frac{\rho_{f,\beta_1} (D)}{D}\bigg)^{-1}
= (1-\chi_D(p_0)p_0^{-1})
  \sum_{j \geq 0} 1_{p_0^{j}|\beta} \chi_D(p_0^{j}) ~.
\end{align*}

Part (c) follows from \cite[Corollary 2]{stewart}, which implies, as shown
in \cite[Lemma 31]{browning-munshi}, that any quadratic polynomial
$P(x)=a_1x^2+a_2x+a_3$ of discriminant $D_P = a_2^2-4a_1a_3$ satisfies
$$|\{x \in \Z/p^{k}\Z: P(x) \equiv 0 \Mod{p^k}\}| \leq 2
p^{v_p(D_P)/2}~.$$
Consider for fixed $y$ the polynomial $P_y(x)= f(x,y) =ax^2 + bxy +cy^2$
of discriminant $y^2 D(f)$. There are less than $p^{\alpha - k}$
values of $y \in \Z/p^{\alpha}\Z$ for which $p^k\|y$. Thus
$$\rho_{f,0}(p^{\alpha}) 
\leq 2 \sum_{k=0}^{\alpha-1} p^{\alpha - k} p^{k + v_p(D(f))/2}
\ll \alpha p^{\alpha}~.$$
\end{proof}

An immediate corollary, which will be essential for the
$W$-trick, states that the $\rho$-densities are constant for lifts of
non-zero residues $\beta \Mod{p^{\alpha}}$ to higher powers of $p$:

\begin{corollary}\label{lifted-densities}
 Let $p$ be a prime and suppose that $\alpha \geq v_p(D)$ and that 
$\beta \not\equiv 0 \Mod{p^{\alpha}}$.
Then
$$\rho_{f,\beta}(p^{\alpha})p^{-\alpha} =
\rho_{f,\beta + kp^{\alpha}}(p^{\alpha+1})p^{-(\alpha+1)}$$
for all $k \in \Z/p\Z$.
\end{corollary}
\begin{proof}
 This follows from part (a) and (b) of Lemma \ref{rho-bounds}.
\end{proof}

\begin{proof}[Proof of Proposition \ref{r'_D(f)-major_arc_estimate}]
Multiplicativity and the assumption on $\beta_0$ yield
\begin{align*}
  \E_{n \in P} \sum_{d|n} \chi_D(d)
&= \bigg(\prod_{p|q_0}  
   \sum_{\alpha \geq 0} 1_{p^{\alpha}|\beta_0} \chi_D(p^{\alpha})
   \bigg)
   \E_{m \leq M} \sum_{\substack{d|qm+\beta}} \chi_D(d)~,
\end{align*}
where $(q,\beta)=1$, and $q$ and $q_0$ have the same prime divisors.
We will estimate the mean value of $\sum_{\substack{d|qm+\beta}}
\chi_D(d)$
by the hyperbola method.
Recall that $\chi_D$ is a character to the modulus $\prod_{p|D} p$
(cf.~\cite[Ch.9.3]{MV}) and let $\chi_D^*$ be the character to the
modulus 
$\prod_{p|q} p$ that is induced by $\chi_D$. (Note that $q$ is divisible
by $\prod_{p|D} p$.)
Thus $\chi_D^*(n)$ is only non-zero when $n$ is coprime to $q$. Then
\begin{align*}
  \E_{0 \leq m \leq M} \sum_{\substack{d|qm+\beta}} \chi_D(d)
= \E_{0 \leq m < M} \Bigg(
  \sum_{\substack{d|(q m + \beta)\\ d \leq T}} \chi_D^*(d)
 + \sum_{\substack{d|(q m + \beta)\\ d > T}} \chi_D^*(d) \Bigg)~,
\end{align*}
where the cut-off $T$ will be chosen as $T=\sqrt{qM}$.
We begin with the \emph{large divisors}.

Writing $G:=(\Z/q \Z)^*$ and denoting its dual group by $\hat G$, 
we have
\begin{align*}
\E_{0 \leq m < M} \sum_{\substack{d|(q m + \beta)\\ d > T}} \chi_D^*(d)
& = \frac{1}{M} \sum_{n \leq Mq}
  \frac{1}{|\hat G|} 
  \sum_{\chi \in \hat G} 
  \overline{\chi}(\beta) \chi(n)
  \sum_{\substack{d|n \\ d >T}}
  \chi_D^*(d)~.
\end{align*}
Since $T^2 \geq qM$, this equals
\begin{align*}
  \frac{1}{M} 
  \frac{1}{|\hat G|} 
  \sum_{\chi \in \hat G} 
  \overline{\chi}(\beta) 
  \sum_{m \leq T}
  \chi(m)
  \sum_{T< d \leq Mq/m}
  (\chi\chi_D^*)(d)~.
\end{align*}
The character $\chi\chi_D^*$ is a non-principal character to the modulus
$q$ unless $\chi$ is the character $\chi_D^{*}$
induced by $\chi_D$. 
We consider the cases $\chi = \chi_D^{*}$ and 
$\chi \not= \chi_D^{*}$ separately. 
For $\chi \not= \chi_D^{*}$, we have
\begin{align*}
  \frac{1}{M} 
  \frac{1}{|\hat G|} 
  \sum_{\substack{\chi \in \hat G \\ \chi \not= \chi_D^{*}}} 
  \overline{\chi}(\beta) 
  \sum_{m \leq T}
  \chi(m)
  \sum_{T< d \leq Mq/m}
  (\chi\chi_D)(d) 
= O(qT/M) ~.
\end{align*}
If $\chi =\chi_D^{*}$, we have
\begin{align*}
& \frac{1}{M} 
  \frac{1}{|\hat G|}
  \overline{\chi_D^{*}}(\beta) 
  \sum_{m \leq T}
  \chi_D^{*}(m)
  \sum_{T< d \leq \frac{Mq}{m}}
  (\chi_D^{*}\chi_D)(d) \\
&=\chi_D(\beta)
  \frac{1}{M} 
  \frac{1}{|\hat G|}
  \sum_{\substack{n \leq Mq \\ (n,q)=1}}
  \sum_{m: m^2 < n} 1_{m|n} \chi_D(m)\\
&=\chi_D(\beta)
  \frac{1}{\phi(q)M}
  \sum_{\substack{m \leq \sqrt{qM} \\ (m,q)=1 }}
  \chi_D(m) \Big(\frac{Mq - m^2}{m} \frac{\phi(q)}{q} + O(q)\Big)\\
&=\chi_D(\beta)
  \frac{1}{qM}
  \sum_{\substack{m \leq \sqrt{qM} \\ (m,q)=1 }}
  \chi_D(m) \frac{Mq - m^2}{m} 
  + O ( \log q \sqrt{q/M} ) \\
&=\chi_D(\beta)
  \sum_{\substack{m \leq \sqrt{qM} \\ (m,W)=1 }}
  \frac{\chi_D(m)}{m} 
  - 
  \frac{\chi_D(\beta)}{qM}
  \sum_{\substack{m \leq \sqrt{qM} \\ (m,W)=1 }}
  \chi_D(m) m 
  + O (  \log q \sqrt{q/M} )~.
\end{align*}
The second term is seen to be small, that is $O(\sqrt{q/M})$, by partial
summation.
The first sum, 
$\sum_{m \leq \sqrt{qM}} \frac{\chi_D^*(m)}{m}$,
is a partial sum of the convergent series
$$\sum_{m \geq 1} \frac{\chi_D^*(m)}{m} 
= \prod_{p|q}\Big(1 - \frac{\chi_D(p)}{p}\Big)
  L(1,\chi_D)
= \prod_{p|q}\Big(1 - \frac{\chi_D(p)}{p}\Big)
  \frac{2\pi h(D)}{k(D)\sqrt{-D}} ~.$$
Bounding their difference sum by partial summation, we obtain
\begin{align*}
\sum_{m \leq \sqrt{qM}} \frac{\chi_D^*(m)}{m}
= \prod_{p\leq w(N)}\Big(1 - \frac{\chi_D(p)}{p}\Big)
  L(1,\chi_D) + O( \sqrt{q/M} )~.
\end{align*}
Hence, the large divisors satisfy
\begin{align*}
\E_{m \leq M} \sum_{\substack{d|(q m + \beta)\\ d > T}} \chi_D^*(d)
 = \chi_D(\beta) \prod_{p|q}\Big(1 - \frac{\chi_D(p)}{p}\Big)
  L(1,\chi_D) + O(\log q \sqrt{q/M} )~.
\end{align*}

Concerning the \emph{small divisors} sum, we obtain
\begin{align*}
 \E_{m \leq M}
  \sum_{\substack{d|(q m + \beta); \\ d \leq T}}
  \chi_D^*(d)
&= \frac{1}{M} \sum_{\substack{d \leq T}}
   \Big(\chi_D^*(d) \frac{M}{d} + O(1)\Big)
 = \sum_{\substack{d \leq T}} 
   \frac{\chi_D^*(d)}{d} + O\Big(\frac{T}{M}\Big) \\
&= L(1,\chi_D)
   \prod_{p|q}\Big(1 - \frac{\chi_D(p)}{p}\Big)
  + O\Big(\frac{q}{T}+ \frac{T}{M}\Big)\\
&=  L(1,\chi_D)
   \prod_{p|q}\Big(1 - \frac{\chi_D(p)}{p}\Big)
  + O(\sqrt{q/M})~.
\end{align*}
Putting things together, we obtain the estimate
\begin{align*}
 \E_{0\leq m \leq M} \sum_{d|qm+\beta} \chi_D(d)  
=(1 + \chi_D(\beta)) L(1,\chi_D)
 \prod_{p|q} \Big(1 - \frac{\chi_D(p)}{p} \Big)
 + O(\log q \sqrt{q/M})~,
\end{align*}
which proves the result. 
\end{proof}

\section{$W$-trick}\label{W-trick}

The aim of this section is to find a decomposition of the function $r_f$
into a sum of functions that are equidistributed in residue classes to
small moduli.

In the case of primes, see \cite{green-tao-longprimeaps}, this was
achieved by defining $W = \prod_{p \leq w(N)} p$ to be the product of
primes up to $w(N)$, where $w:\N \to \R$ is a slowly growing function.
For $n$ with $\gcd(n,W)=1$ the von Mangoldt function then splits as
$$\Lambda(n) 
= \sum_{a \in (\Z/W\Z)^*} \Lambda(n)1_{n \equiv a \Mod{W}}~,$$
and it suffices to consider the functions 
$n \mapsto \Lambda(Wn+a)$, $a\in (\Z/W\Z)^*$, which are equidistributed
in residue classes to small moduli.

In the case of the divisor function, the most natural decomposition makes
use of the restricted divisor function that only counts divisors coprime
to $W$ (and is thus likely to be a quasirandom function): define
$$\tau'(n) := \sum_{d: (d,W)=1} 1_{d|n} ~.$$
Then
$$\tau (n) = \tau'(n) \sum_{w} 1_{w|n} ~,$$
where $w$ runs over all integers entirely composed of primes $\leq w(N)$.
The second factor, $\sum_{w} 1_{w|n}$, is almost periodic.
Indeed, let $\alpha(p)$ be such that
$$p^{\alpha(p)-1} < \log^{C_1+1} N \leq p^{\alpha(p)}~.$$
Then any number $n$ that is divisible by some $w$ as above with
$p^{\alpha(p)} | w$ for some $p \leq w(N)$ belongs to the exceptional set
$X_0$ from Lemma \ref{exceptional-set}.
Choosing 
$$
\W := \prod_{p \leq w(n)}p^{\alpha(p)}~,
$$
one can achieve that the second factor is a periodic function of period
$\W$, when adjusting the values of $\tau$ at exceptional integers. 
This way, it suffices to consider the functions of the form
$n \mapsto \tau'(\W n+a)$ for non-zero residues $a\in [\W]$.
In fact, observing that $\tau(\W n+a)=\tau'(\W n+a)\sum_{w|\W}1_{w|a}$ for
unexceptional values of $a$, we essentially consider functions of the
form $n \mapsto \tau(\W n+a)$.

In the case of representation functions a very similar $W$-trick works.
We use the same choice of $\W$ as in the divisor function case above.
\begin{definition} 
Let $\A$ be the set of residues $a \Mod{\W}$ such that
$$\rho_{f,a}(\W) > 0$$
and such that $a \not\equiv 0 \Mod{p^{\alpha(p)}}$~.
\end{definition}
Thus $\A$ contains only residue classes that are representable by $f$,
and every $n \in [N]$ which fails to satisfy the second condition, that
is, for which $n \equiv 0 \Mod{p^{\alpha(p)}}$ holds, belongs to the
exceptional set $X_0$ from Lemma \ref{exceptional-set}.
\begin{definition}[Normalised and $W$-tricked representation function] 
 Let $\beta \in \A$ and define $r'_{f,\beta}:[N/\W] \to \R$ by
\begin{align*}
r'_{f,\beta}(m) 
&= \frac{k(D)\sqrt{-D}}{2 \pi} r_f(\W m + \beta)
\left(\frac{\rho_{f,\beta}(\W)}{\W}\right)^{-1} \\
&= \frac{k(D)\sqrt{-D}}{2 \pi} r_f(\W m + \beta) \prod_{p \leq w(N)}
\left(\frac{\rho_{f,\beta}(p^{\alpha(p)})}{p^{\alpha(p)}}\right)^{-1}~.
\end{align*}
\end{definition}
Thus, by Lemma \ref{distribution_in_APs}
$$\E_{n \leq M}r'_{f,\beta}(m) = 1 + O(\W^3 M^{-1/2})~.$$

\subsection{The major arc estimate}\label{major-arcs}

Our next aim is to give a major arc estimate for the $W$-tricked function
$r'_{f,b}$: 
we show that this function has, up to a small error, a
constant average on arithmetic progressions whose common difference is
small in the sense that it is $w(N)$-smooth.

\begin{definition}
 An integer is called $k$-\emph{smooth}, when each of its prime divisors
is
at most $k$.
\end{definition}

\begin{proposition}[Major arc analysis for $r'_{f,\beta}$]
\label{major_arc_estimate}
Let $P \subseteq [N/\W]$ be a progression of $w(N)$-smooth common
difference $q_1$ and let $\beta \in \A$. 
If $P = \{q_1 m + q_0: 0 \leq m < M \}$ has length $M$, then
$$
  \E_{n \in P} r'_{f,\beta}(n) 
= \E_{0 \leq m < M} r'_{f,\beta}(q_1m+q_0)
= 1 + O\Big(\frac{\W(\W q_1)^2}{M^{1/2}}\Big)~.
$$
\end{proposition}
\begin{proof}
Corollary \ref{lifted-densities} implies
$$\frac{\rho_{f,\beta}(\W)}{\W}
= \frac{\rho_{f,\W q_0+\beta}(\W q_1)}{\W q_1}~.$$
Hence the result follows from Lemma \ref{distribution_in_APs}.
\end{proof}

\subsection{$W$-tricked majorant}
Finally, we need to slightly adapt our majorant function for $r_f$ to its
$W$-tricked version. Let $\beta \in \A$.
Then Lemma \ref{ind-of-g} and Lemma \ref{general-r_g} yield the pointwise
majorisation
\begin{align*}
r'_{f,\beta}(n) 
\leq \sum_{f' \sim_g f} r'_{f',\beta}(n) 
&= O(1) (\rho_{f,\beta}(\W)\W^{-1})^{-1} 
   \sum_{d|\W n +\beta} \chi_D(d) \\
&= O(1) \prod_{\substack{p < w(N) }} 
   \Big(1 - \frac{\chi_D(p)}{p}\Big)^{-1}
   \sum_{\substack{d|\W n + \beta \\ p|d \Rightarrow p> w(N)}}
   \chi_D(d)~,
\end{align*}
where the last step uses Lemma \ref{rho-approximate}.
Since each function $r'_{f',\beta}(n)$ has average order $1+o(1)$, the
last expression is of bounded average order.
Thus, the function
$$r'_{D(f)}(n) 
:= \sum_{\substack{d|n \\ p|d \Rightarrow p>w(n) }} \chi_D(d)$$
may be used in place of $r_{D(f)}$ to run through the construction of
the majorant as in Section \ref{rep-majorant-section}.
In view of the results from that section and the remarks at the end
of Section \ref{bar-reduction-section} we find
$$ r'_{D(f)}(n)
\leq  \beta'_{D,\gamma}(n) \nu'_{D,\gamma}(n)~,$$
where
\begin{align*}
\nu'_{D,\gamma} (n) = 
 \sum_{s \geq 2/\gamma}^{(\log \log N)^3}
 \sum_{i \geq \log_2 s - 2}^{6 \log \log \log N} 
 \sum_{u \in U(i,s)}
 2^s
 1_{u|n} 
 \tau'_{D,\gamma}(n)~,
\end{align*}
with
$$
\tau'_{D,\gamma}
:=
 \sum_{\substack{d \in \pdspan \\ 
       p|d \Rightarrow p > w(N) }} 
 1_{d|n} 
 \chi\left(\frac{\log d}{\log N^{2\gamma}}\right)~,
$$
and
\begin{equation*}
\beta'_{D,\gamma}(n) 
:= 
 \sum_{\substack{m \in \qdspan \\ p| m \Rightarrow p>w(N)
       \\ m < N^{\gamma}}}
 \bigg(
  \sum_{\substack{d \in \qdspan \\ p| d \Rightarrow p>w(N)}} 
  1_{m^2d|n}~\mu(d)~
  \chi \Big(\frac{\log d}{\log N^{\gamma}}\Big)
 \bigg)^2. 
\end{equation*}
For the two factors $\beta'_{D,\gamma}$ and $\nu'_{D,\gamma}$ one shows in
the same way as for the original majorants that 
$$C(\beta'_{D,\gamma})
:= \sqrt{\log N}
\prod_{\substack{q \in \Q_D \\ q \leq w(N)}} (1 + q^{-1})^{-1}
\E_{n \leq N}  
\beta'_{D,\gamma}(n)$$ and 
$$C(\nu'_{D,\gamma})
:=
\frac{1}{\sqrt{\log N}}
\prod_{\substack{p \in \P_D\\ p < w(N)}} (1-p^{-1})^{-1}
\E_{n \leq N} \nu'_{D,\gamma}(n)$$
are bounded independently of $N$. 
Since $\beta'_{D,\gamma}$ and $\nu'_{D,\gamma}$ are given by short
divisor sums running over coprime sets of divisors, the average order of
their product satisfies
$$\E_{n \leq N} \beta'_{D,\gamma}(n)\nu'_{D,\gamma}(n)
= \E_{n \leq N} \beta'_{D,\gamma}(n) 
  \E_{m \leq N}  \nu'_{D,\gamma}(m) + N^{O(\gamma) - 1}~.$$
Indeed, for coprime integers $y_1,y_2 < N^{\gamma}$, we have
$$\E_{n\leq N} 1_{y_1 y_2 | n} = \frac{1}{y_1 y_2} + O(N^{2\gamma - 1})
= \E_{n\leq N} 1_{y_1 | n}\E_{n\leq N} 1_{y_2 | n} + O(N^{2\gamma -
1})~,$$
and since the total number of divisors in the sum $\E_{n\leq N}
\beta'_{D,\gamma}(n)\nu'_{D,\gamma}(n)$ is $N^{O(\gamma)}$, the statement
follows.

Since $\prod_{p \leq w(N)} (1 - \chi_D(p) p^{-1})^{-1} = C + o(1)$ for
some constant $C$, we have proved the following lemma.
\begin{lemma}[$W$-tricked majorant]
Let $\beta \in \A$, then
$$r'_{f,\beta}(m) 
\leq \beta'_{D,\gamma}(\W m + \beta) \nu'_{D,\gamma}(\W m + \beta)$$
for all $m \leq N/\W$.
Furthermore, there is a positive real number $C_{D,\gamma} = O(1)$ such
that 
$$\E_{n \leq N}
\frac{\beta'_{D,\gamma}(n)\nu'_{D,\gamma}(n)}{C_{D,\gamma}}
= 1 + o(1)~.$$
\end{lemma}

\section[Local factors and a $W$-tricked version of the main theorem]
{Local factors and the reduction of the main theorem to a
$W$-tricked version}\label{local-factors}

Define the smoothed representation function
$\bar r_f: [N] \to \R$ by
\begin{align*}
\bar r_f(n) := r_f(n) 1_{n \Mod{\W} \in \A}~.
\end{align*}
According to the definition of $\A$, this function satisfies the
conditions of Lemma \ref{bar-reduction}.
Thus it suffices to study correlations of functions $\bar r_f$ in order to
prove the main theorem.
As the main theorem will show, the asymptotic behaviour of these
correlations, 
\begin{align} \label{local-factors-1}
\sum_{ n \in K \cap \Z^d} 
\bar r_{f_1}(\psi_1(n)) \dots \bar r_{f_t}(\psi_t(n))~,
\end{align}
is determined by the local behaviour of the affine-linear system $\Psi$
modulo small primes.

By splitting the summation range into progressions of common difference
$\W$, we reduce the task of estimating \eqref{local-factors-1} to an
assertion, Proposition \ref{W-red-Main-thm} below, about the uniformity of
the $W$-tricked representation functions.
Local factors measuring irregularities of the system $\Psi$ modulo small
primes will appear in this process.

Define for fixed quadratic forms $f_1, \dots, f_t$ and for an
affine-linear system $\Psi: \Z^d \to \Z^t$ the set of residues
\begin{align*}
\A_{\Psi} :=& 
\{ a \in [\W]^d:
 \psi_i(a) \in \A_{f_i} \text{ for all } i \in [t]\}\\ 
=& \{ a \in [\W]^d: \prod_{i=1}^{t}
 \rho_{f_i, \psi_i(a)}(\W) > 0 \text{ and }
 \prod_{i=1}^{t}\psi_i(a) \not\equiv 0 \Mod{p^{v_p(\W)}} \} ~.
\end{align*}
Notice that any $n$ with non-zero contribution to
\eqref{local-factors-1} is congruent modulo $\W$ to an element of this
set.
For a fixed element $a \in \A_{\Psi}$ let 
$\tilde \Psi = (\tilde \psi_1, \dots, \tilde \psi_t):\Z^d \to \Z^t$ be
the affine-linear system satisfying
$$\psi_i(\W m + a) = \W \tilde \psi_i(m) + c_i(a)$$
with $c_i(a) \in [\W]$.
Thus, $\psi_i(a) \equiv c_i(a) \Mod{\W}$, and $\psi_i$ and $\tilde \psi_i$
only differ in the constant term.

The main result will be deduced from the following Proposition.

\begin{proposition}\label{W-red-Main-thm} 
Let $\Psi:\Z^d \to \Z^t$ be a finite complexity system of forms, let
$a \in \A_{\Psi}$, and let $\tilde \Psi:\Z^d \to \Z^t$ be defined as
above.
Then 
$$\sum_{m \in \Z^d \cap K'} \prod_{i=1}^t r'_{f_i,c_i(a)}(\tilde
\psi_i(m)) 
= \vol(K') +  o\Big({(N/\W)^d}\Big)~,$$
where $K' \subseteq [-N/\W, N/\W]^d$ is a convex body such that 
$\W \tilde \Psi (K') + c(a) \subseteq [1,N]^t$.
\end{proposition}

For every $a \in \A_{\Psi}$, define the convex body 
$$K_a:= \{x \in \R^d : \W x + a \in K\}$$
and note that $\vol(K_a) = \vol(K)/\W^d$.
Then we can rewrite \eqref{local-factors-1} by means of
Proposition \ref{W-red-Main-thm} as follows
\begin{align} \nonumber
& \sum_{ n \in K \cap \Z^d} 
  \bar r_{f_1}(\psi_1(n)) \dots \bar r_{f_t}(\psi_t(n)) \\
\nonumber
&=\sum_{a \in \A_{\Psi}}
  \sum_{m \in K_a \cap \Z^d} 
  \prod_{i=1}^t
  r_{f_i}(\psi_i(\W m + a)) \\ 
\nonumber
&=\sum_{a \in \A_{\Psi}}
  \sum_{m \in K_a \cap \Z^d} 
  \prod_{i=1}^t
  r'_{f_i,c_i(a)}(\tilde \psi_i(m)) 
  \frac{\rho_{f_i,\psi_i(a)}(\W)}{\W} \frac{2\pi}{k(D_i)\sqrt{-D_i}}\\
\nonumber
&=\frac{\vol(K) + o(N^d)}{\W^d} 
  \sum_{a \in \A_{\Psi}}
  \prod_{i=1}^t
  \frac{\rho_{f_i,\psi_i(a)}(\W)}{\W} \frac{2\pi}{k(D_i)\sqrt{-D_i}} \\
\label{lf}
&=\Big(\vol(K)\prod_{j=1}^t \frac{2\pi}{k(D_j)\sqrt{-D_j}}  + o(N^d)\Big)
  \E_{a\in [\W]^d}
  1_{a \in \A_{\Psi}}
  \prod_{i=1}^t
  \frac{\rho_{f_i,\psi_i(a)}(\W)}{\W}~.
\end{align}
By the Chinese remainder theorem, the above expectation is in fact a
product over local densities, that is
\begin{align}\label{CRT} \nonumber
& \E_{a\in [\W]^d}
  1_{a \in \A_{\Psi}}
  \prod_{i=1}^t
  \frac{\rho_{f_i,\psi_i(a)}(\W)}{\W} \\
&= \prod_{p<w(N)} \E_{a \in (\Z/p^{\alpha(p)}\Z)^d}
  \prod_{i=1}^t
  \frac{\rho_{f_i,\psi_i(a)}(p^{\alpha(p)})}{p^{\alpha(p)}}
  1_{\psi_i(a) \not\equiv 0 \Mod{p^{\alpha(p)}}}~,
\end{align}
where $\alpha(p)=v_p(\W)$.
To complete the proof that \eqref{lf} and \eqref{CRT} indeed imply the
main theorem, two further lemmas are required.
The first shows that the above factors at primes are essentially local
factors:
\begin{lemma}[Local factors] \label{localfactorslemma}
Let $p$ be a prime. Then
\begin{align*}
 \E_{a \in (\Z/p^{\alpha(p)}\Z)^d}
  \prod_{i=1}^t
  \frac{\rho_{f_i,\psi_i(a)}(p^{\alpha(p)})}{p^{\alpha(p)}}
  1_{\psi_i(a) \not\equiv 0 \Mod{p^{\alpha(p)}}}
= \beta_p + O\big(\log^{-C_1/5} N\big)~,
\end{align*}
where 
\begin{align*}
\beta_p
:= \lim_{m \to \infty} 
  \E_{a \in (\Z/p^{m}\Z)^d}
  \prod_{i=1}^t
  \frac{\rho_{f_i,\psi_i(a)}(p^{m})}{p^{m}}
\end{align*}
is the local factor at $p$.
\end{lemma}

The second lemma is an estimate of the local factors.
\begin{lemma} \label{beta_p-bound-lemma}
Let $\Psi=(\psi_1,\dots,\psi_t):\Z^d \to \Z^d$ be a system of
affine-linear forms for which no two forms $\psi_i$ and $\psi_j$ are
affinely dependent, and all of whose linear coefficients are bounded by
$L$. Then
\begin{equation*}
\beta_p = 1 + O_{t,d,L}(p^{-2})~.
\end{equation*}
\end{lemma}
Thus,
$$\prod_{\substack{p\leq w(N)}} \beta_p = 
\Big(1 + O_{t,d,L}\Big(\frac{1}{w(N)}\Big)\Big)\prod_{p} \beta_p ~.$$
A second consequence of this lemma is that
$\beta_p + O(\log^{-C_1/5} N) = \beta_p(1+O(\log^{-C_1/5} N))$ for all 
$p\gg 1$. 
For the remaining $p \ll 1$, we require an upper bound on $\beta_p$.
Since Lemma \ref{rho-bounds} implies 
$\rho_{f_i,A}(p^{\alpha(p)})p^{-\alpha(p)} 
\ll \alpha(p) \ll \log \log N$ for any 
$A \in \Z/p^{\alpha(p)}\Z$, we may deduce from Lemma
\ref{localfactorslemma} the very crude bound
$\beta_p \ll (\log \log N)^t$.
Thus, by \eqref{lf}, \eqref{CRT} and the two lemmas stated above, we
obtain
\begin{align*}
& \sum_{ n \in K \cap \Z^d} 
  \bar r_{f_1}(\psi_1(n)) \dots \bar r_{f_t}(\psi_t(n))\\ 
&= (\beta_{\infty} + o(N^d))
  \prod_{p<w(N)} \Big(\beta_p + O(\log^{-C_1/5} N)\Big) \\ 
&= (\beta_{\infty} + o(N^d)) \Big(1 + O(\log^{-C_1/5} N)\Big)^{\pi(w(N))}
  \Big(\prod_{p<w(N)} \beta_p 
  + O\Big(\frac{(\log \log N)^{O(t)}}{\log^{C_1/5} N}\Big) \Big)
   \\ 
&= \beta_{\infty} \prod_{p<w(N)} \beta_p + o(N^d) \\
&= \beta_{\infty} \prod_{p} \beta_p + o(N^d)~,
\end{align*}
where we used that $w(N)= \log \log N$.
Apart form the proof of the two lemmas, we have reduced the task of
establishing the main theorem to that of proving Proposition
\ref{W-red-Main-thm}.

We conclude this section with the proofs of the lemmas, for the
purpose of which the following notion is introduced.
\begin{definition}[Local divisor densities]
For a given system $\Psi=(\psi_1,\dots,\psi_t)$ of affine-linear forms,
positive integers $d_1, \dots, d_t$ and their least common multiple 
$m:= \lcm(d_1, \dots, d_t)$ define \emph{local divisor densities} by 
$$
\alpha_{\Psi}(d_1, \dots, d_t) := \E_{n \in (\Z/m\Z)^d} \prod_{i \in [t]}
1_{\psi_i(n) \equiv 0 \Mod{d_i}}~.
$$
\end{definition}

\begin{proof}[Proof of Lemma \ref{localfactorslemma}]
We shall show more precisely that $\beta_p$ satisfies
\begin{align}\label{term1}
 \beta_p
=&\E_{a \in (\Z/p^{\alpha(p)}\Z)^d}
  \prod_{i=1}^t
  \frac{\rho_{f_i,\psi_i(a)}(p^{\alpha(p)})}{p^{\alpha(p)}}
  1_{\psi_i(a) \not\equiv 0 \Mod{p^{\alpha(p)}}}  \\ \label{term2}
& + O\Big((\alpha(p))^t\Big) 
 \sum_{
 \substack{ a_1, \dots, a_t: \\ M:=\max_i a_i \\ \geq \alpha(p)}} 
 \E_{a \in (\Z/p^{M}\Z)^d} 
 \prod_{i=1}^t 1_{\psi_i(a) \equiv 0 \Mod{p^{a_i}}}~.
\end{align}
Suppose $m > \alpha(p)$.
We split the sum
$\E_{a \in (\Z/p^{m}\Z)^d}
  \prod_{i=1}^t
  \rho_{f_i,\psi_i(a)}(p^{m}) p^{-m}$
over residues $a$ into two parts according to whether
\begin{align*}
 \prod_{i=1}^t 1_{\psi_i(a) \not\equiv 0 \Mod{p^{\alpha(p)}}} = 1 
\text{ or }0~.
\end{align*}
First note that for any $a$ with $\psi_i(a) \not\equiv 0
\Mod{p^{j}}$ for all $i \in [t]$, any lift $\Psi(a+kp^{j})$, $k \in
[p]^{d}$ is
component-wise divisible to the same powers of $p$ as $\Psi(a)$.
Hence, Corollary \ref{lifted-densities} implies
\begin{align*}
& \E_{a \in (\Z/p^{m}\Z)^d}
  \prod_{i=1}^t
  \frac{\rho_{f_i, \psi_i(a)}(p^{m})}{p^{m}}
  1_{\psi_i(a) \not\equiv 0 \Mod{p^{\alpha(p)}}} \\
&=\E_{a \in (\Z/p^{\alpha(p)}\Z)^d}
  \prod_{i=1}^t
  \frac{\rho_{f_i, \psi_i(a)}(p^{\alpha(p)})}{p^{\alpha(p)}}
  1_{\psi_i(a) \not\equiv 0 \Mod{p^{\alpha(p)}}}~.
\end{align*}
Thus, the terms of the first type give rise to \eqref{term1}.
Combining part (a), (b) and (c) of Lemma \ref{rho-bounds} yields the
general bound
\begin{align*}
\frac{\rho_{f_i, \psi_i(a)}(p^{m})}{p^{m}}
\ll \sum_{k=0}^m 1_{\psi_i(a) \equiv 0 \Mod{k}}~,
\end{align*}
which shows that terms of the second type are bounded by
\begin{align*}
 O\Big((\alpha(p))^t\Big)
 \sum_{
 \substack{0\leq a_1, \dots, a_t \leq m: 
\\ M:=\max_i a_i \\ \geq \alpha(p)}} 
 \E_{a \in (\Z/p^{M}\Z)^d} 
 \prod_{i=1}^t 1_{\psi_i(a) \equiv 0 \Mod{p^{a_i}}}~.
\end{align*}
This proves the above expression for $\beta_p$.
In order to establish the lemma, it thus remains to bound \eqref{term2},
that is, the sum
over divisor densities 
$$\delta_p:= 
  \sum_{\substack{ a_1, \dots, a_t \\ M:=\max_i a_i \geq \alpha(p)}}
  \alpha_{\Psi}(p^{a_1}, \dots, p^{a_t})~.$$
Since the coefficients of $\dot\Psi$ are bounded, we have
$$\alpha_{\Psi}(p^{a_1}, \dots, p^{a_t})=
  \E_{n\in (\Z/p^{\max_i a_i}\Z)^d} 
  \prod_{i=1}^t 
  1_{\psi_i(n) \equiv 0 \Mod{p^{a_i}}}
\ll p^{-\max_i a_i}~,$$
which yields
\begin{align*}
\delta_p \ll 
 \sum_{\substack{a_1,\dots, a_t \\ \max_i a_i \geq \alpha(p)}}
  p^{-\max_i a_i}~.
\end{align*}
Recall that
$$\alpha(p) = v_p(\W) = (C_1+1)\frac{\log \log N}{\log p} + O(1)$$
for some sufficiently large integer $C_1$. 
Estimating the number of tuples $(a_1, \dots, a_t)$ with $\max_i a_i=j$
crudely by $(j+1)^{t}$, we conclude that for $p \leq w(n)= \log \log N$
\begin{align*}
\delta_p 
&\ll \sum_{j \geq C_1(\log\log N)/2\log p}  p^{-j} j^{t} \\
&\ll \sum_{j \geq C_1(\log\log N)/2\log p}  p^{-j/2} \\
&\ll (\log N)^{-C_1/4}~.
\end{align*}
Hence, $(\alpha(p))^t \delta_p \ll (\log N)^{-C_1/5}$, which proves the
result.
\end{proof}

\begin{proof}[Proof of Lemma \ref{beta_p-bound-lemma}]
We may assume that $p$ is large enough so that $p \nmid D_1 \dots D_t$.
For such primes Lemma \ref{rho-bounds}(c) yields
\begin{align*}
 \beta_p 
 &= \lim_{m \to \infty} \E_{a \in (\Z/p^m\Z)^d}
    \prod_{j \in [t]} \frac{\rho_{f_i, \psi_i(a)}(p^m)}{p^m} \\
 &= \lim_{m \to \infty} \E_{a \in (\Z/p^m\Z)^d}
    \prod_{i=1}^t \Big(1-\chi_{D_i}(p)p^{-1}\Big) 
    \sum_{j\geq 0} 1_{p^j|\psi_i(a)}\chi_{D_i}(p^j) \\
 &= \sum_{a_1, \dots, a_t} \alpha(p^{a_1}, \dots, p^{a_t})
   \prod_{j \in [t]} \Big(1-\chi_{D_j}(p)p^{-1}\Big) 
   \chi_{D_j}(p^{a_j})~.
\end{align*}
By splitting the sum $\sum_{a_1,\dots,a_t}$ into terms according to
whether no $a_i$ is non-zero, exactly one $a_i$ is non-zero, or at least
two $a_i$ are non-zero, we obtain for $\beta_p$ the following.
\begin{align*}
\beta_p 
&= \sum_{a_1, \dots, a_t} \alpha(p^{a_1}, \dots, p^{a_t})
   \prod_{j \in [t]} \Big(1-\chi_{D_j}(p)p^{-1}\Big) 
   \chi_{D_j}(p^{a_j})  \\
&= \prod_{j \in [t]} \Big(1-\chi_{D_j}(p)p^{-1}\Big) 
   \bigg\{ 1 + 
   \sum_{i=1}^t \sum_{a_i>0} 
   \chi_{D_i}(p^{a_i}) p^{-a_i} \bigg\}
 + O\bigg( \sum_{\substack{a_1,\dots a_t: \\ 
   \text{at least}\\ \text {two } a_i>0}}
   \alpha(p^{a_1},\dots p^{a_t})   
   \bigg)~.
\end{align*}
Here we used the fact that, for sufficiently large $p$ with respect to
$t,d$ and $L$, we have 
$\alpha(p^{a_1},\dots,p^{a_t}) = p^{-a_i}$
whenever $a_i$ is the only non-zero exponent.

It is easy to see that the main term equals $1 + O_t(p^{-2})$.
Concerning the error term, we employ the fact that we are dealing with a
finite complexity system of forms. That is, since no two forms are
affinely related, we have for every $p$ which is sufficiently large with
respect to $t,d,L$ that
$$
 \alpha(p^{a_1},\dots,p^{a_t}) 
\leq p^{- \max_{i\not=j}( a_i + a_j)}
\leq p^{-1 -\max_{i}a_i}$$
whenever at least two $a_i$ are non-zero.
There are at most $tj^{t-1}$ choices of coefficients $a_1, \dots, a_t$
that satisfy $\max_i a_i = j$, and thus the contribution of the error term
to the value of $\beta_p$ may be bounded by
\begin{align*}
 O\Big(\sum_{j \geq 1} tj^{t-1} p^{-j-1}\Big) = O_t(p^{-2})~.
\end{align*}
This proves the lemma.
\end{proof}

\subsection*{Simultaneous majorant}
To summarise, we reduced the task of proving the main theorem to that of
proving Proposition \ref{W-red-Main-thm}.
This will carried out by the nilpotent Hardy-Littlewood method in the
remainder of this paper.
In order to apply the method, specifically Proposition \ref{v.neumann}
below, we require for every occurring collection of
$\{r'_{f_i, c_i(a)}: i=1,\dots, t\}$, $a \in \A_{\Psi}$, a pseudorandom
majorant that simultaneously majorises all $r'_{f_i, c_i(a)}$.
The following function has the required majorant property:
\begin{align} \label{sim-majorant}
 \sigma_{(f_i),a}: [N/\W] \to \R^+~, 
\quad 
 \sigma_{(f_i),a}(m) 
 := \E_{i\in[t]} 
    \frac{\beta'_{D_i,\gamma}(\W m + b_i(a))
          \nu'_{D_i,\gamma}(\W m + b_i(a))}
         {C_{D_i,\gamma}}~.
\end{align}

\section{Linear forms and correlation conditions}
\label{linear-forms-and-correlation-conditions}

In this section we check that the majorant $\sigma_{(f_i),a}$ defined in
\eqref{sim-majorant} for a collection of $W$-tricked representation
functions $r'_{f_1,c_1(a)}, \dots, r'_{f_t,c_t(a)}$ is (after a minor
technical modification) indeed a pseudorandom measure, that is, satisfies
the linear forms and correlation conditions. 

Write $M = N/ \W$, let $M'$ be a prime satisfying 
$M < M' \leq O_{t,d,L}(M)$, and define
$\sigma^{*}_{(f_i),a}:[M'] \to \R^+$ by
\begin{align*}
 \sigma^*_{(f_i),a}(n)= 
 \left\{
 \begin{array}{ll}
 \frac{1}{2}(1+\sigma_{(f_i),a}(n)) &\text{ if } n \leq M \cr
 1 & \text{ if } M < n \leq M'~.
 \end{array}
 \right.
\end{align*}
As is seen in \cite[App.D]{green-tao-linearprimes},
$\sigma^{*}_{(f_i),a}$ is $D$-pseudorandom if the following
two propositions, which are technical reductions of the linear forms and
correlation conditions from \cite{green-tao-linearprimes}, hold true. 

\begin{proposition}[$D$-Linear forms estimate]
\label{linear-forms-estimate}
Let $1 \leq d,t \leq D$ and 
let $(i_1, \dots, i_t) \in [t]^t$ be an arbitrary collection of indices.
For any finite complexity system
$\Psi : \Z^d \to \Z^t$ with bounded coefficients $\|\Psi\|_{N} \leq D$
and every convex body $K \subseteq [0,N]^d$ such that $\Psi(K) \subseteq
[1,N/\W]^t$, the estimate
\begin{align}\label{lfc-correlation}
  \E_{n\in \Z^d \cap K } 
& \prod_{j \in [t]} 
  \nu'_{D_{i_j},\gamma}(\W \psi_j(n) + b_{i_j})
  \beta'_{D_{i_j},\gamma}(\W \psi_j(n) + b_{i_j}) \\
\nonumber
&= \bigg( 1 + O_D\Big(\frac{N^{d - 1 +O_D(\gamma)}}{\vol(K)}\Big) + o_D(1)
   \bigg) 
   \prod_{j = 1}^t 
   C_{D_{i_j},\gamma}
\end{align}
holds, provided $\gamma$ was small enough.
\end{proposition}

\begin{proposition}[Correlation estimate]
\label{verification-of-C-Condition}
For every $1 < m_0 \leq D$ there exists a function $\sigma_{m_0} : \Z_{M'}
\to \R^+$ with bounded moments 
$\E_{n \in \Z_{M'}} \sigma_{m_0}^q(n) \ll_{m,q} 1$
such that for every interval $I \subset \Z_{M'}$, every $1 \leq m \leq
m_0$ and every $m$-tuple $(i_1, \dots, i_{m}) \in [t]^{m}$ and every
choice of (not necessarily distinct) $h_1,\dots,h_{m} \in \Z_{M'}$ we have
$$\E_{n \in I} \prod_{j \in [m]} 
  \nu'_{D_{i_j},\gamma}(\W(n+h_j)+b_{i_j}) 
  \beta'_{D_{i_j},\gamma}(\W(n+h_j)+b_{i_j})
\leq \sum_{1 \leq i < j \leq m} \sigma_{m_0}(h_i - h_j)~,$$ 
provided $\gamma$ was small enough.
\end{proposition}

Recall that the $W$-tricked majorant $\nu'_{D_j,\gamma}(\W m + b_j(a))
\beta'_{D_j,\gamma}(\W m + b_j(a))$ for $r'_{f_j,a}$ has divisor sum
structure:
\begin{align}\label{div-sum-structure} \nonumber
&\nu'_{D_j,\gamma}(n) \beta'_{D_j,\gamma}(n) \\ \nonumber
&=  
\Bigg( 
 \sum_{s = 2 /\gamma}^{(\log \log N)^3} 
 \sum_{i = \log_2 s-2}^{6 \log \log \log N}
 \sum_{u \in U(i,s)} 
 \sum_{\substack{d \in \langle \P_{D_j} \rangle \\ (d,uW)=1 }}
 \sum_{v|u}
 2^s 1_{d|n} 1_{u|n} \chi \left( \frac{\log d}{\log N^{2\gamma}} \right)
\Bigg)\times \\
&\qquad\times 
\Bigg(
 \sum_{\substack{m_j \in \langle \Q_{D_j} \rangle \\ (m_j,W)=1}}
 \chi \left( \frac{\log m_j}{\log N^{2\gamma}} \right)
 1_{m_j^2|n}
 \Bigg(
  \sum_{\substack{\eps \in \langle Q_{D_i}\rangle \\ (\eps,W)=1 }} 
  1_{\eps m_j^2|n} \mu(\eps) 
  \chi \left( \frac{\log \eps}{\log N^{\gamma}} \right)
 \Bigg)^2
\Bigg)
~.
\end{align}
The function $\chi$ above is a cut-off. As no characters appear in this
section, there is no danger of confusion.

Our strategy to prove the linear forms estimate is as follows.
The first step is to show that in order to asymptotically evaluate
\eqref{lfc-correlation} we may ignore all terms that arise from divisor
densities of \emph{dependent} divisibility events, that is, events
$\{n:\prod_{i \in [t]} 1_{a_i|\psi_i(n)}\}$ where $(a_1, \dots, a_t)$ are
not pairwise coprime.
The second step is the observation that the densities of
independent divisibility events are, up to a small error, independent of
the system $\Psi$ of forms, which will finally allow us to reduce the
verification of the linear forms condition to the task of verifying it
separately for each of the two factors of each of the majorants in the
case where the $\Psi:\Z \to \Z$ is the identity function.
The same strategy was used in \cite[\S6]{m-divisorfunction}

The main tool to exploit the divisor sum structure of our majorants is
the following simple lemma (see \cite[App.A]{green-tao-linearprimes} for
a proof).

\begin{lemma}[Volume packing argument] Let $K \subseteq [-B,B]^d$ be a
convex body and $\Psi$ a system of affine-linear forms. Then 
$$
 \sum_{n\in \Z^d \cap K} \prod_{i\in [t]} 1_{d_i|\psi_i(n)} 
 = \vol(K) \alpha (d_1, \dots, d_t) + O(B^{d-1} \lcm(d_1,\dots, d_t))~.
$$
\end{lemma}

In order to remove the above mentioned dependent divisibility events, we
need to replace $\chi$ by a multiplicative function.
A way to achieve this has been found by Goldston and Y{\i}ld{\i}r{\i}m and
was employed and modified by Green and Tao \cite{green-tao-linearprimes}
to check the linear forms condition for their majorant function for
$W$-tricked primes.
In this respect, the proof of Proposition \ref{linear-forms-estimate}
below builds on \cite[App.D]{green-tao-linearprimes}.
In particular, we shall employ many of the small technical arguments
from there.

Recall that the cut-off $\chi$ was chosen to be a smooth, compactly
supported function satisfying $\int_0^{\infty} |\chi'(x)|^2 dx = 1$.
Let $\vartheta$ be the modified Fourier transform of $\chi$, defined via
$$e^x \chi(x) = \int_{\R} \vartheta(\xi)e^{-ix\xi} d\xi~.$$
Fourier inversion, compact support and smoothness of $\chi$, and partial
integration yield the bound
$$\vartheta(\xi) \ll_A (1 + |\xi|)^{-A}$$ for all $A > 0$.
Green and Tao make use of this rapid decay to truncate the integral
representation of $\chi$ as follows.
Let $I=\{\xi \in \R : |\xi|\leq \log^{1/2} N^{\gamma}\}$, then for any 
$A > 0$
\begin{align} \label{integral-rep-for-chi}
\nonumber
\chi(\frac{\log m}{\log N^{\gamma}}) 
&= \int_{\R}
   {m}^{-\frac{1+i\xi}{\log N^{\gamma}}} \vartheta(\xi)~d\xi \\
&= \int_{I} m^{-\frac{1+i\xi}{\log N^{\gamma}}} \vartheta(\xi)~d\xi
 + O_A(m^{-1/\log N^{\gamma}} \log^{-A} N^{\gamma})~.
\end{align}
This truncation will later-on simplify the process of swapping integrals
and summations.
We proceed to check the linear forms estimate.

\subsection*{Proof of Proposition \ref{linear-forms-estimate}}
Define the system $\Phi=(\varphi_j)_{j\in[t]}: \Z^d \to \Z^t$ by 
$\varphi_j(n):= \W \psi_j(n) + b_{i_j}$.
A prime $p$ is called exceptional for $\Phi$ if the reduction of $\Phi$
modulo $p$ has affinely depended forms.
For the system defined here, all exceptional primes are bounded by
$w(N)+O(D)$.
All information we will use about $\Phi$ are the bound on exceptional
primes and the fact that it has finite complexity. 
Consider an arbitrary cross term that appears on the left hand side of
\eqref{lfc-correlation} when inserting the definition
\eqref{div-sum-structure} and fixing the parameters $s_j,i_j,u_j$ for each
factor. That is, we consider 
\begin{align*}
&\E_{n\in \Z^d \cap K} \prod_{j \in [t]}
\Bigg( 
 \sum_{\substack{d_j \in \langle \P_{D_j} \rangle \\ (d_j,u_jW)=1}}
 \sum_{v_j|u_j}
 2^{s_j} 1_{d_ju_j|\varphi_j(n)} 
 \chi\left(\frac{\log d_j}{\log N^{\gamma}}\right)
 \Bigg) \times \\
& \qquad \qquad \qquad
 \sum_{\substack{m_j \in  \langle \Q_{D_i} \rangle \\ (m_j,W)=1}} 
 \chi\left(\frac{\log m_j}{\log N^{\gamma}}\right)
 \Bigg(
 \sum_{e_j \in \langle \Q_{D_i} \rangle} 
 1_{e_j m_j^2|\varphi_j(n)} \mu(e_j) 
 \chi\left(\frac{\log e_j}{\log N^{\gamma}}\right)
\Bigg)^2  \\
&= \sum_{\substack{\mathbf{d},\mathbf{m}, \mathbf{e},\mathbf{e'} }}
 \Bigg(
 \prod_{j \in [t]} 2^{s_j}
 \mu(e_j) \mu(e'_j) \tau(u_j)
 \prod_{\substack{x \in \\ \{d_j,m_j,e_j,e'_j\}}}
 \chi\left(\frac{\log x}{\log N^{\gamma}}\right)
 \Bigg)
 \E_{n\in \Z^d \cap K}
 \prod_{i \in [t]}
 1_{u_id_im_i^2\eps_i|\varphi_i(n)}
\end{align*}
where $\eps_i=\lcm(e_i,e'_i)$ and where we denote by bold letters such as
$\mathbf{d}$ any $t$-tuple of positive $w(N)$-smooth integers which we
shall implicitly assume to satisfy the correct multiplicative
restrictions, e.g.~$d_i \in \<\P_i\>$ and $(d_i,v_i W)=1$ in this case.
 
Note that $u_jd_jm_j^2\eps_j = N^{O(\gamma)}$ for all summands with
non-zero contribution. 
Indeed, $d_j,e_j,e'_j, m_j \leq N^{\gamma}$ by definition of the cut-off. 
We have $u_j < N^{\gamma}$ by construction of the divisor majorant, as
the $u_j$ arise as divisors of certain numbers bounded by $N^{\gamma}$
(c.f.~also
the remarks following Proposition 4.2 of \cite{m-divisorfunction}). 
Therefore, the volume packing lemma implies
\begin{align*}
& \E_{n\in \Z^d \cap K} \prod_{i \in [t]}
 2^{s_i} \tau(u_i) 1_{u_id_im_i^2\eps_i|\varphi_i(n)} \\
&= \alpha_{\Phi}(u_1d_1m_1^2\eps_1, \dots, u_td_tm_t^2\eps_t)
   \prod_{i \in [t]} 2^{s_i} \tau(u_i)
 + O\Big(M^{d-1 + O(\gamma)}/\vol(K)\Big)~,
\end{align*}
where the bound $2^{s_j} \leq 2^{(\log\log N)^3} \ll M^{\gamma}$
allowed to hide the factors $2^{s_j}$ in the error term.

Since $u_jd_jm_j^2\eps_j = N^{O(\gamma)}$, there are only $N^{O(\gamma)}$
terms all together in all sums of the majorant, including those over
$s_j$, $i_j$ and $u_j$.
This and the boundedness of $\chi$ imply that the volume packing error
term has a total contribution of $O(M^{d-1 + O(\gamma)}/\vol(K))$ towards
\eqref{lfc-correlation}, and we are left to deal with the main term, that
is
\begin{align*}
 \sum_{\mathbf{d},\mathbf{m},\mathbf{e},\mathbf{e'}}
 \alpha_{\Phi}((u_id_im_i^2\eps_i)_{i\in[t]}) 
 \prod_{j \in [t]} 2^{s_j}
 \mu(e_j) \mu(e'_j) \tau(u_j)
 \prod_{\substack{x \in \\ \{d_j,m_j,e_j,e'_j\}}}
 \chi\!\left(\frac{\log x}{\log N^{\gamma}}\right)~.
\end{align*}
Next we show that we may assume that each $u_i$ is coprime to
$u_jd_jm_j^2\eps_j$ for all $j\not=i$ and that $(u_i,d_im_i^2\eps_i)=1$.
These properties yield
$$\alpha_{\Phi}(u_1d_1m_1^2\eps_1, \dots, u_td_tm_t^2\eps_t)
= \alpha_{\Phi}(d_1m_1^2\eps_1, \dots, d_tm_t^2\eps_t) 
\frac{1}{u_1 \dots u_t}~.$$
We shall also abbreviate $\mathbf{u}=(u_1,\dots,u_t)$, implicitly assuming
that the conditions $u_j \in U(i_j,s_j)$ on these tuples still apply.
\begin{claim} \label{claim-1}
 For all choices of $(s_j)_{j\in[t]}$ and $(i_j)_{j\in[t]}$ we have
\begin{align*}
 &\sum_{\substack{\mathbf{d},\mathbf{m},\mathbf{e},\mathbf{e'}}} 
 \sum_{\mathbf{u}} 
 \alpha_{\Phi}((u_id_im_i^2\eps_i)_{i\in[t]}) 
 \prod_{j \in [t]}  2^{s_j}
 \mu(e_j) \mu(e'_j) \tau(u_j)
 \prod_{\substack{x \in \\ \{d_j,m_j,e_j,e'_j\}}}
 \chi\!\left(\frac{\log x}{\log N^{\gamma}}\right)\\
&= \sum_{\substack{\mathbf{d},\mathbf{m},\mathbf{e},\mathbf{e'}}} 
 \alpha_{\Phi}((d_im_i^2\eps_i)_{i\in[t]})  
 \dsum_{u_1, \dots, u_t}
 \prod_{j \in [t]} 
 \frac{2^{s_j}\tau(u_j)}{u_j}
 \mu(e_j) \mu(e'_j) 
 \prod_{\substack{x \in \\ \{d_j,m_j,e_j,e'_j\}}}
 \chi\!\left(\frac{\log x}{\log N^{\gamma}}\right) \\
 & \qquad + O_D(N^{-(\log \log N)^{-4}})~,
\end{align*}
where $\dsum$ indicates that the sum is extended only over choices 
$(u_1, \dots, u_t)$ satisfying the coprimality conditions
$(u_i,u_jd_jm_j^2\eps_j)=1$ whenever $i\not=j$ and
$(u_i,d_im_i^2\eps_i)=1$ for $i\in[t]$.
\end{claim}
\begin{proof}
 We have to bound the contribution from excluded choices of 
$(u_1, \dots, u_t)$. Any prime divisor of any $u_i$ is at least as large
as $N^{1/(\log \log N)^3}$ by construction.
Thus, whenever the coprimality conditions fail, the divisibility events
we
are considering are included in $\{n:p^2|\prod_{i\in[t]} \phi_i(n)\}$ for
some $p > N^{1/(\log \log N)^3}$.
By finite complexity and the bounds on exceptional primes of $\Phi$, we
have 
$$
\sum_{N^{(\log \log N)^{-3}}<p<N^{\gamma}} 
\E_{n\in \Z^d \cap K} 1_{p^2|\prod_{i} \phi_i(n)}
\ll_t \sum_{N^{(\log \log N)^{-3}}<p<N^{\gamma}} p^{-2}
=O_t(N^{-(\log \log N)^{-3}})~.
$$
We will make use of this with the help of Cauchy-Schwarz.
Since $2^{s_j} \leq 2^{(\log \log N)^3}$ and since $\chi^2$ is at
most $1$, we can crudely bound the following second moment
\begin{align*}
& \E_{n \in \Z^d \cap K} \prod_{i\in[t]}
 \Bigg(
 \sum_{\substack{\mathbf{d},\mathbf{m},\mathbf{u},\mathbf{e},\mathbf{e'}}}
 1_{u_id_im_i^2\eps_i|\psi_i(n)} 
 2^{s_i} \tau(u_i) 
 \prod_{\substack{x \in \\ \{d_j,m_j,e_j,e'_j\}}}
 \chi^2\!\left(\frac{\log x}{\log N^{\gamma}}\right) \Bigg)^2 \\
&\ll 2^{2t(\log \log N)^3} 
 \prod_{i\in[t]}
 \Bigg( \E_{n \in \Z^d \cap K}
 \Bigg(
 \sum_{\substack{\mathbf{d},\mathbf{m},\mathbf{u},\mathbf{e},\mathbf{e'}
        \\ \in [N^{\gamma}]^t}}
 1_{u_id_im_i^2\eps_i|\psi_i(n)} 
 \tau(u_i) \Bigg)^{2t}
 \Bigg)^{1/t} \\
&\ll (\log N)^{O(t)}  2^{2t(\log \log N)^3}~.
\end{align*}
The combination of these two bounds proves the claim.
\end{proof}
Note that the same argument furthermore shows that the main term from
Claim \ref{claim-1} equals
\begin{align} \label{claim-1-main} \nonumber
&\sum_{\substack{\mathbf{d}, \mathbf{m}, \mathbf{e},\mathbf{e'}}} 
 \alpha_{\Phi}((d_im_i^2\eps_i)_{i\in[t]})
 \prod_{j \in [t]} 
 \sum_{u_j}
 \frac{2^{s_j}\tau(u_j)}{u_j}
 \mu(e_j) \mu(e'_j) 
 \prod_{\substack{x \in \\ \{d_j,m_j,e_j,e'_j\}}}
 \chi\!\left(\frac{\log x}{\log N^{\gamma}}\right) \\ 
& \qquad + O_D(N^{-(\log \log N)^{-4}})~.
\end{align}
Thus, we are left to deal with the main term in \eqref{claim-1-main}. 
We proceed by inserting the integral representation
\eqref{integral-rep-for-chi} of each of the $4t$ factors involving
$\chi$.
Multiplying out this product we obtain a main term and number error terms.
Since
$\chi(\frac{\log m}{\log N^{\gamma}}) \ll m^{-1/\log N^{\gamma}}$,
all these error terms may be seen to be of the same form, which allows
us to combine them into one error term. 
Writing $z_{j,k}= (1 + i \xi_{j,k})/\log N^{\gamma}$ for
$j\in[t]$,$k\in[4]$ and noting that
$|z_{j,k}| \ll (\log N^{\gamma})^{-1/2}$, the main term from
\eqref{claim-1-main} is seen to equal
\begin{align}\label{insert-ints}
 \sum_{\mathbf{d},\mathbf{m},\mathbf{e},\mathbf{e'}}
 \bigg( 
 \prod_{i \in [t]} 
 \sum_{u_i} \frac{2^{s_i}\tau(u_i)}{u_i} 
 \bigg)
 \alpha_{\Phi}(d_1m_1^2\eps_1, \dots, d_t&m_t^2\eps_t) \times \\ 
\nonumber
 \Bigg\{ \int_{I} \dots \int_I
 \prod_{j \in [t]} 
 \mu(e_j) \mu(e'_j)
 e_j^{-z_{j,1}}
 {e'_j}^{-z_{j,2}} &
 d_j^{-z_{j,3}} 
 m_j^{-z_{j,4}}
 \prod_{k\in[4]}
 \vartheta(\xi_{j,k})
 ~d\xi_{j,k}
 \\   
\nonumber
+ &O_A \Big( \log^{-A} N^{\gamma} 
           \prod_{j \in [t]} (e_je'_jd_jm_j)^{-1/\log N^{\gamma}}
 \Big)
 \Bigg\}.
\end{align}
The error term here indeed has small contribution: 
On the one hand, we have
\begin{align*}
 \sum_{s_1, \dots, s_t} \sum_{i_1, \dots, i_t} 
 \prod_{j \in [t]}
 \sum_{\substack{u_j \in U(i_j,s_j)}}
 \frac{2^{s_{j}}\tau(u_j)}{u_j}
= O(1)~.
\end{align*}
See the proof of \cite[Proposition 4.2]{m-divisorfunction} for
details.
On the other hand, the divisor sum is bounded:
\begin{align*}
& \sum_{\mathbf{d},\mathbf{m},\mathbf{e},\mathbf{e'}} 
 \alpha_{\Phi}(d_1m_1^2\eps_1, \dots, d_tm_t^2\eps_t)
 \prod_{j \in [t]} (e_je'_jd_jm_j)^{-1/\log N^{\gamma}} \\
&= \sum_{\mathbf{d},\mathbf{m},\mathbf{e},\mathbf{e'}} 
 \prod_{\substack{p>w(N)\\ p^{a_i}\|d_im_i^2\eps_i}}
 \alpha_{\Phi}(p^{a_1}, \dots, p^{a_t}) 
 \prod_{\substack{j \in [t]\\ p^{a'_j}\|e_je'_jd_jm_j }} 
 (p^{a_j + a'_j})^{-1/\log N^{\gamma}} \\
& \ll \prod_{p>w(N)} (1 + p^{-(1+1/\log N^{\gamma})})^{-O(t)}
\ll \log^{O(t)}N~,
\end{align*}
Here, we crudely bounded the number of occurring $t$-tuples
$(a_1, \dots, a_t)$ that satisfy
$\max_i a_i = k$ by $k^{O(t)}$ and apply to each of these tuples the
bound $\alpha_{\Phi}(p^{a_t}, \dots, p^{a_t}) \ll p^{-k}$.

Thus, when choosing $A$ in \eqref{insert-ints}
sufficiently large, the error term above makes
a total contribution of $\ll_A \log^{-A/2}N^{\gamma}$.

It remains to estimate the main term from above. 
Changing the order of summation and integration leads to an absolutely
convergent sum in the integrand. 
Since the range of integration is compact this change is permitted and,
hence, the main term is equal to
\begin{align} \label{whole_integral} 
 \int_{I} \dots \int_I 
 \Bigg(
 \sum_{\mathbf{d},\mathbf{m},\mathbf{e},\mathbf{e'}}
 \bigg(
 \prod_{i \in [t]}
 \sum_{u_i} \frac{2^{s_i}\tau(u_i)}{u_i} 
 \bigg)
 \alpha_{\Phi}((d_km_k^2\eps_k)_{k\in[t]}) \\ \nonumber
 \times
 \prod_{j \in [t]} 
 \mu(e_j) \mu(e'_j)
 e_j^{-z_{j,1}}
 {e'_j}^{-z_{j,2}} &
 d_j^{-z_{j,3}} 
 m_j^{-z_{j,4}}
 \Bigg) 
 \prod_{k\in[4]}
 \vartheta(\xi_{j,k})
 ~d\xi_{j,k}~.
\end{align}
Our next aim is to show that all relevant terms in the integrand are in
fact the independent terms, that is, they are those terms for which the
$t$ products $u_id_im_i^2\eps_i$, $i\in[t]$ are pairwise coprime.
This will eventually allow us to swap the sums with the product while only
introducing a small error. 
For the $u_i$ we have just done this.

Since each entry of $\mathbf{d}$, $\mathbf{m}$, $\mathbf{e}$, and
$\mathbf{e'}$ is completely composed of primes $\geq w(N)$, the following
claim holds.
\begin{claim} \label{claim-2}
 We have
\begin{align*}
&\sum_{\mathbf{d}, \mathbf{m}, \mathbf{e},\mathbf{e'}}
 \alpha_{\Phi}((d_im_i^2\eps_i)_{i\in[t]}) 
 \prod_{j \in [t]} 
 \mu(e_j) \mu(e'_j)
 e_j^{-z_{j,1}}
 {e'_j}^{-z_{j,2}} 
 d_j^{-z_{j,3}} 
 m_j^{-z_{j,4}} \\
&= (1 + O_D(w(N)^{-1}))
 \dsum_{\mathbf{d}, \mathbf{m}, \mathbf{e},\mathbf{e'}}
 \alpha_{\Phi}((d_im_i^2\eps_i)_{i\in[t]}) 
 \prod_{j \in [t]} 
 \mu(e_j) \mu(e'_j)
 e_j^{-z_{j,1}}
 {e'_j}^{-z_{j,2}} 
 d_j^{-z_{j,3}} 
 m_j^{-z_{j,4}} ~,
\end{align*}
where $\dsum$ indicates that the summation is extended only over choices
of $t$-tuples that satisfy the coprimality condition
$(d_im_i\eps_i,d_{i'}m_{i'}\eps_{i'})=1$ for any $i \not= i'$.
\end{claim}
\begin{proof}
 Note that the summand is multiplicative and may be written as a product
over primes $p>w(N)$.
Any summand 
$$
 \alpha_{\Phi}((d_im_i^2\eps_i)_{i\in[t]}) 
 \prod_{j \in [t]} 
 \mu(e_j) \mu(e'_j)
 e_j^{-z_{j,1}}
 {e'_j}^{-z_{j,2}} 
 d_j^{-z_{j,3}} 
 m_j^{-z_{j,4}}
$$
with entries failing coprimality may be factorised into a product of one
factor of the same form that satisfies coprimality and one factor for
which every prime $p$ that appears as a divisor of some $d_im_i\eps_i$
divides at least another $d_{i'}m_{i'}\eps_{i'}$, $i'\not=i$.
For a fixed tuple $(k_1,\dots, k_t)$ of the latter type (that is, $p|k_i$
implies $p|\prod_{i'\not=i}k_{i'}$), the contribution may be bounded
as follows employing the triangle inequality:
\begin{align*}
&\alpha(k_1,\dots,k_t)
\bigg| \dsum_{\substack{\mathbf{d}, \mathbf{m}, \mathbf{e},\mathbf{e'}
 \\(d_im_i\eps_i,~ k_1 \dots k_t)=1}}
 \alpha_{\Phi}((d_im_i^2\eps_i)_{i\in[t]}) 
 \prod_{j \in [t]} 
 \mu(e_j) \mu(e'_j)
 e_j^{-z_{j,1}}
 {e'_j}^{-z_{j,2}} 
 d_j^{-z_{j,3}} 
 m_j^{-z_{j,4}}
 \bigg| \\
&= \alpha(k_1,\dots,k_t) \\
& \qquad \bigg| \prod_{p|k_1 \dots k_t} \big(1 + O(p^{-1})\big)
 \dsum_{\substack{\mathbf{d}, \mathbf{m}, \mathbf{e},\mathbf{e'}}}
 \alpha_{\Phi}((d_im_i^2\eps_i)_{i\in[t]}) 
 \prod_{j \in [t]} 
 \mu(e_j) \mu(e'_j)
 e_j^{-z_{j,1}}
 {e'_j}^{-z_{j,2}} 
 d_j^{-z_{j,3}} 
 m_j^{-z_{j,4}} \bigg| \\
& \leq \alpha(k_1,\dots,k_t) \\
& \qquad \prod_{p|k_1 \dots k_t} \big(1 + O(p^{-1})\big)
\bigg| \dsum_{\substack{\mathbf{d}, \mathbf{m}, \mathbf{e},\mathbf{e'}}}
 \alpha_{\Phi}((d_im_i^2\eps_i)_{i\in[t]}) 
 \prod_{j \in [t]} 
 \mu(e_j) \mu(e'_j)
 e_j^{-z_{j,1}}
 {e'_j}^{-z_{j,2}} 
 d_j^{-z_{j,3}} 
 m_j^{-z_{j,4}} \bigg|~.
\end{align*}
Next, we bound the sum over all occurring terms 
$\alpha(k_1,\dots,k_t) \prod_{p|k_1 \dots k_t} (1 + O(p^{-1}))$.
Written as a product over primes, a crude bound for this quantity is
given by
\begin{align*}
 \prod_{p>w(N)}
 \Big\{1 + 
   \sum_{\substack{a_1, \dots, a_t : \\ \text{at least two} \\ a_i>0}}
   O(a_1^4 + \dots + a_t^4) 
   \alpha_{\Phi}(p^{a_1}, \dots, p^{a_t}) \big(1+ O(p^{-1})\big)
 \Big\} - 1~,
\end{align*}
where we used the very crude bound $\tau_5(p^{a_i}) \ll a_i^4$ on the
generalised divisor function $\tau_5$.
The five factors correspond to $d_i$, $m_i^2$, $\eps_i / e_i$,
$\eps_i/e'_i$ and $e_ie'_i/\eps_i$.
To further bound the above expression, we observe that the number of
tuples $(a_1,\dots,a_t)$ with $\max_i a_i = k$ is at most $t(k+1)^{t-1}$.
For such choices of $(a_1,\dots,a_t)$, we have $\sum_i a_i^4 \leq tk^4$
and $\alpha_{\Phi}(p^{a_1}, \dots, p^{a_t}) \leq p^{-k-1}$, since $\Phi$
has finite complexity and at least two of the $a_i$ are non-zero.
Further, for large enough $p$, we have
$p^{-k}t^2k^{t+3}(1+O(p^{-1})) < p^{-3k/4}$ for all $k\geq1$.
We certainly may assume that $N$ is large enough for $p>w(N)$ to satisfy
this condition.
Thus 
\begin{align*}
 \sum_{\substack{a_1, \dots, a_t : \\ \text{at least two} \\ a_i>0}}
   O(a_1^4 + \dots + a_t^4) 
   \alpha_{\Phi}(p^{a_1}, \dots, p^{a_t})\big(1+ O(p^{-1})\big)
\leq \sum_{k \geq 1} p^{-3k/4 -1}
\leq (p^{-1 - 1/2})~.
\end{align*}
Since 
$$\prod_{p>w(N)}( 1 + p^{-3/2}) -1 \leq \sum_{n>w(N)}n^{-3/2} \ll
w(N)^{-1/2}~,$$
the result follows.
\end{proof}
Note that in the above claim 
\begin{align}\label{rewritten_term} \nonumber
&\dsum_{\mathbf{d}, \mathbf{m}, \mathbf{e},\mathbf{e'}}
 \alpha_{\Phi}((d_im_i^2\eps_i)_{i\in[t]}) 
 \prod_{j \in [t]} 
 \mu(e_j) \mu(e'_j)
 e_j^{-z_{j,1}}
 {e'_j}^{-z_{j,2}} 
 d_j^{-z_{j,3}} 
 m_j^{-z_{j,4}} \\
&= 
 \dsum_{\mathbf{d}, \mathbf{m}, \mathbf{e},\mathbf{e'}}
 \prod_{j \in [t]} 
 \frac{\mu(e_j) \mu(e'_j)}{\eps_j}
 e_j^{-z_{j,1}}
 {e'_j}^{-z_{j,2}} 
 d_j^{-1-z_{j,3}} 
 m_j^{-2-z_{j,4}}
\end{align}
holds. 
The last step of the rearrangement is to show that we may swap the
inner product and sum in the integrand.
\begin{claim}\label{claim-4}
The sum and product in \eqref{rewritten_term} may be interchanged:
\begin{align*}
& \dsum_{\mathbf{d}, \mathbf{m}, \mathbf{e},\mathbf{e'}}
 \prod_{j \in [t]} 
 \frac{\mu(e_j) \mu(e'_j)}{\eps_j}
 e_j^{-z_{j,1}}
 {e'_j}^{-z_{j,2}} 
 d_j^{-1-z_{j,3}} 
 m_j^{-2-z_{j,4}} \\
&= (1 + O(w(N)^{-1/2}))
 \prod_{j \in [t]} 
 \sum_{d_j, m_j, e_j, e'_j}
 \frac{\mu(e_j) \mu(e'_j)}{\eps_j}
 e_j^{-z_{j,1}}
 {e'_j}^{-z_{j,2}} 
 d_j^{-1-z_{j,3}} 
 m_j^{-2-z_{j,4}} 
~.
\end{align*}
\end{claim}
\begin{proof}
 The proof of this claim is essentially the same as the one of the
previous claim.
\end{proof}
The next claim will imply that the integral \eqref{whole_integral}
equals, up to a small error, the integral of the main term from Claim
\ref{claim-4}.
\begin{claim}
 \begin{align*}
 \int_{I} \dots &\int_I 
 \Bigg|
 \prod_{j \in [t]} 
 \sum_{d_j, m_j, e_j, e'_j}
 \frac{\mu(e_j) \mu(e'_j)}{\eps_j}
 e_j^{-z_{j,1}}
 {e'_j}^{-z_{j,2}} 
 d_j^{-1-z_{j,3}} 
 m_j^{-2-z_{j,4}}
 \prod_{k\in[4]} \vartheta(\xi_{j,k})
 \Bigg|
 \prod_{\substack{(j',k') \\ \in[t]\times[4]}}
 ~d\xi_{j',k'} \\
&= O(1)~.
 \end{align*}
\end{claim}
\begin{proof}
(Cf.~\cite{green-tao-linearprimes}, equation (D.23) and the proof
thereof.)
We begin by writing the integrand as a product over primes
\begin{align*}
&\Bigg|
 \prod_{j \in [t]} 
 \sum_{d_j, m_j, e_j, e'_j}
 \frac{\mu(e_j) \mu(e'_j)}{\eps_j}
 e_j^{-z_{j,1}}
 {e'_j}^{-z_{j,2}} 
 d_j^{-1-z_{j,3}} 
 m_j^{-2-z_{j,4}}
 \prod_{k\in[4]} \vartheta(\xi_{j,k})
 \Bigg| \\
&\ll_A
\prod_{j \in [t]}
\prod_{q\in\Q_j}
\big(1 - q^{-1-z_{j,1}} - q^{-1-z_{j,2}} + q^{-1-z_{j,1}-z_{j,2}}\big) 
\prod_{\substack{p \in\P_j }}(1-p^{-1-z_{j,3}})^{-1}
\prod_{k\in[4]} (1+|\xi_{j,k}|)^{-A} 
\end{align*}
By the prime number theorem in arithmetic progressions, we have for 
$\Re s>0$
$$\sum_{p\in\Q_i} p^{-1-s} = \frac{1}{2}\log \frac{1}{s} + O_{D_i}(1)~.$$
This also holds for $\P_i$ in place of $\Q_i$.
Thus, choosing $A$ sufficiently large the above is seen to be bounded by
\begin{align*}
&\prod_{j \in [t]}
\prod_{k\in[4]} (1+|\xi_{j,k}|)^{-A}
\prod_{q\in\Q_j}
\big(1 - q^{-1-z_{j,1}} - q^{-1-z_{j,2}} + q^{-1-z_{j,1}-z_{j,2}}\big) 
\prod_{\substack{p \in\P_j }}(1-p^{-1-z_{j,3}})^{-1} \\
&\ll \prod_{j \in [t]}
\prod_{k\in[4]} (1+|\xi_{j,k}|)^{-A}
|z_{j,1}|^{1/2}|z'_{j,2}|^{1/2}|z_{j,1}+z'_{j,2}|^{-1/2}|z_{j,3}|^{-1/2}\\
&\ll \log^{t}N^{\gamma} \log^{-t}N^{\gamma}
\prod_{j \in [t]} (1+|\xi_{j,1}|)^{1/2}(1+|\xi_{j,2}|)^{1/2}
\prod_{k\in[4]} (1+|\xi_{j,k}|)^{-A} \\
&\ll \prod_{j \in [t]}
\prod_{k\in[4]}
(1+|\xi_{j,k}|)^{-A/2}~.
\end{align*}
For any $A>2$ the integral of the final expression is $O(1)$.
\end{proof}
Together with Claim \ref{claim-2}, equation \eqref{rewritten_term} and
Claim \ref{claim-4}, the above Claim implies that the integral
\eqref{whole_integral} is given by
\begin{align*}
\int_{I} \dots \int_I
 \prod_{j \in [t]} 
 \sum_{d_j, m_j, e_j, e'_j}
 \frac{\mu(e_j) \mu(e'_j)}{\eps_j}
 e_j^{-z_{j,1}}
 {e'_j}^{-z_{j,2}} 
 & d_j^{-1-z_{j,3}} 
 m_j^{-2-z_{j,4}}
 \prod_{k\in[4]} \vartheta(\xi_{j,k})
 \prod_{\substack{(j',k') \\ \in[t]\times[4]}}
 ~d\xi_{j',k'} \\
& \times \bigg(
 \prod_{i \in [t]}
 \sum_{u_i} \frac{2^{s_i}\tau(u_i)}{u_i} 
 \bigg)
 + o(1).
\end{align*}
Removing the truncation of the integral again, the latter expression is
seen to equal:
\begin{align*}
 \bigg(
 \prod_{i \in [t]}
 \sum_{u_i} \frac{2^{s_i}\tau(u_i)}{u_i} 
 \bigg) 
  \prod_{j \in [t]} 
  \sum_{d_j, m_j, e_j, e'_j}
  \frac{\mu(e_j) \mu(e'_j)}{d_j m_j^2 \eps_j}
  \prod_{\substack{x \in \\ \{e_j,e'_j,m_j,d_j\}}}
  \chi(\frac{\log x}{\log N^{\gamma}})
  + o(1)~.
\end{align*}
Putting everything together, we have shown that
\begin{align*}
& \E_{n\in \Z^d \cap K } 
  \prod_{j \in [t]} 
  \nu'_{D_{i_j},\gamma}(\varphi_j(n))
  \beta'_{D_{i_j},\gamma}(\varphi_j(n))\\
&= (1 + O_d(w(N)^{-1/2})) \\
&\qquad \prod_{j \in [t]} \Bigg(
 \sum_{s_j}
 \sum_{i_j}
 \sum_{u_j}
 \sum_{d_j}
 \sum_{m_j, e_j, e'_j}
 \frac{2^{s_j}\tau(u_i)}{u_j}
 \frac{\mu(e_j) \mu(e'_j)}{d_j m_j^2 \eps_j}
 \prod_{\substack{x \in \\ \{e_j,e'_j,m_j,d_j\}}}
 \chi(\frac{\log x}{\log N^{\gamma}})
  + o(1)\Bigg)~.
\end{align*}
The last expression now is independent of $\Phi$.
Applying the asymptotic in each of the known one-dimensional cases 
\begin{align*}
 \E_{n \leq N} \nu'_{D_{i_j},\gamma}(n) \beta'_{D_{i_j},\gamma}(n)
&= C_{D_{i_j},\gamma} + o(1)~,
\end{align*}
where $\Phi: \Z \to \Z$ is given by the identity,
implies that each of the factors above is of the correct form. 
This completes proof of the Proposition.

\stsubsection{Proof of Proposition \ref{verification-of-C-Condition}}
The proof of the correlation estimate follows in a very similar manner
to those of the corresponding estimates for the divisor function majorant
in \cite[\S7]{m-divisorfunction} and the von Mangoldt function majorant
from \cite[App.D]{green-tao-linearprimes}.
We restrict attention to the case of pairwise distinct $h_i$; the
remaining case follows, as before, by choosing $\sigma_{m_0}(0)$
sufficiently large.
Employing the volume packing lemma, we may show as in
\cite[\S7]{m-divisorfunction}, that 
 \begin{align*}
\E_{n \in I} \prod_{j \in [m]} 
  \nu'_{D_{i_j},\gamma}(\W(n+h_j)+b_{i_j}) 
  \beta'_{D_{i_j},\gamma}(\W(n+h_j)+b_{i_j})
&\ll \prod_{\substack{p|\Delta\\p>w(N)}} \sum_{a_1,\dots, a_m} 
 \alpha(p^{a_1},\dots,p^{a_m})~,
\end{align*}
where
$\Delta := \prod_{j\not=j'}(\W(h_j-h_{j'})+b_{i_j}-b_{i_{j'}})$.
This estimate allows us to proceed as in \cite[\S7]{m-divisorfunction}.

\section{Application of the transference principle}
\label{application_section}

This section provides a quick overview of the results around the von
Neumann theorem and the inverse theorem for the Gowers norms.
We apply these results in the end of the section to reduce Proposition
\ref{W-red-Main-thm} to a non-correlation estimate. 

In the dense setting, that is, if $g: \Z \to \R$ is a bounded function
with asymptotic density, the Gowers uniformity norms, defined as
$$
\| g \|_{U^s[N]} 
 := \Bigg( 
    \E_{x \in [N]} 
    \E_{h\in[N]^s} 
    \prod_{\omega \in \{0,1\}^s} 
    g(x + \omega \cdot h ) 
    \Bigg)^{1/2^{s}}~,
$$
capture all information on the correlations of $g$ with respect to finite
complexity systems.
This generalises as follows.
\begin{proposition}[Green-Tao \cite{green-tao-linearprimes}, generalised
von Neumann theorem]\label{v.neumann}
Let $t,d,L$ be positive integer parameters.
Then there are constants $C_1$ and $D$, depending on $t,d$ and $L$, such
that the following is true.
Let $C$, $C_1 \leq C \leq O_{t,d,L}(1)$, be arbitrary and
suppose that $N' \in [CN,2CN]$ is a prime. 
Let $\nu: \Z_{N'} \to \R^{+}$ be a $D$-pseudorandom measure, and suppose
that $f_1,\dots,f_t : [N] \to \R$ are functions with 
$|f_i(x)| \leq \nu(x)$ for all $i \in [t]$ and $x \in [N]$. 
Suppose that $\Psi=(\psi_1,\dots,\psi_t)$
is a finite complexity system of affine-linear forms whose
linear coefficients are bounded
by $L$.
Let $K\subset [-N,N]^d$ be a convex body such that 
$\Psi(K) \subset [N]^t$.
Suppose also that 
\begin{equation}\label{uniformity-condition;v.N.Thm}
 \min_{1\leq j \leq t} \|f_j\|_{U^{t-1}[N]} = o(1)~.
\end{equation}
Then we have
$$\sum_{n \in K} \prod_{i \in [t]} f_i(\psi_i(n)) 
  = o(N^d)~.$$ 
\end{proposition}

Establishing the Gowers-uniformity condition
\eqref{uniformity-condition;v.N.Thm} itself is a task that is
conceptually equivalent to that of finding an asymptotic for 
$\sum_{n \in K} \prod_{i\in [t]} f(\psi_i(n))$ directly, and should
therefore not be any easier.
The specific system of affine-linear forms that appears in the definition
of the uniformity norms, however, allows an alternative characterisation
of Gowers-uniform functions.

\subsection*{A characterisation of Gowers-uniform functions}

Whether or not a function $f$ is Gowers-uniform, is characterised by the
non-existence or existence of a polynomial nilsequence\footnote{For
definitions of nilmanifolds and nilsequences, see, for instance,
\cite{green-tao-polynomialorbits}.}
that correlates with $f$.
On the one hand, correlation with a nilsequence obstructs uniformity: 

\begin{proposition}[Green-Tao \cite{green-tao-linearprimes}, Cor. 11.6]
\label{nilsequences-obstruct-uniformity}
Let $s \geq 1$ be an integer and let $\delta \in (0,1)$ be real. 
Let $G/\Gamma = (G/\Gamma, d_{G/\Gamma})$ be an $s$-step nilmanifold with
some fixed smooth metric $d_{G/\Gamma}$ , and let $(F(g(n)\Gamma))_{n \in
\N}$ be a bounded $s$-step nilsequence with Lipschitz constant at most
$L$.
Let $f : [N] \to \R$ be a function that is bounded in the $L_1$-norm,
that is, assume $ \|f\|_{L_1} = \E_{n \in [N]} |f(n)| \leq 1$.
If furthermore
$$ \E_{n \in [N]} f(n) F(g(n)\Gamma) \geq \delta $$
then we have
$$ \|f\|_{U^{s+1}[N]} \gg_{s,\delta,L,G/\Gamma} 1~. $$
\end{proposition}

An inverse result to this statement has been known as the Inverse
Conjecture for the Gowers norms for some time and has recently been
resolved, see \cite{gtz}. 
The inverse conjectures are stated for bounded functions.
With our application to the normalised divisor function in mind, we only
recall the transferred statement,
c.f.~\cite[Prop. 10.1]{green-tao-linearprimes}, here.

\begin{proposition}[Green-Tao-Ziegler,
Relative inverse theorem for the Gowers norms]
\label{inverse theorem} 
For any $0 < \delta \leq 1$ and any $C\geq 20$, there exists a finite
collection $\mathcal M_{s,\delta,C}$ of $s$-step nilmanifolds $G/\Gamma$,
each equipped with a metric $d_{G/\Gamma}$, such that the following
holds. 
Given any $N\geq1$, suppose that $N'\in[CN,2CN]$ is prime, that 
$\nu:[N'] \to \R^+$ is an $(s+2)2^{s+1}$-pseudorandom measure, suppose
that $f:[N] \to \R$ is any arithmetic function with $|f(n)|\leq \nu(n)$
for all $n \in [N]$ and such that $$\|f\|_{U^{s+1}[N]} \geq \delta~.$$
Then there is a nilmanifold $G/\Gamma \in \mathcal M_{s,\delta,C}$ in the
collection and a $1$-bounded $s$-step nilsequence 
$(F(g(n)\Gamma))_{n \in \N}$ on it that has Lipschitz constant 
$O_{s,\delta,C}(1)$, such that we have the correlation estimate
$$|\E_{n\in [N]} f(n)F(g(n)\Gamma)|\gg_{s,\delta,C} 1~.$$
\end{proposition}

This inverse theorem now reduces the required uniformity-norm estimate
\eqref{uniformity-condition;v.N.Thm} to the potentially easier task of
proving that the centralised version of $f$ does not correlate with
polynomial nilsequences.

\subsection{Reduction of the main theorem to a non-correlation
estimate}

We already reduced the main theorem to the $W$-tricked version given in
Proposition \ref{W-red-Main-thm}, which we now restate:

\begin{wredmain}
Let $\Psi:\Z^d \to \Z^t$ be a finite complexity system of forms, let
$a \in \A_{\Psi}$, and let $\tilde \Psi:\Z^d \to \Z^t$ be the translate
of $\Psi$ defined as in Section \ref{local-factors}.
Then 
$$\E_{m \in \Z^d \cap K'} \prod_{i=1}^t r'_{f_i,c_i(a)}(\tilde \psi(m)) 
= 1 + o_{t,d,L}(1)~,$$
where $K' \subseteq [-N/\W,N/\W]^d$ is a convex body such that 
$\W \tilde \Psi (K') + c(a) \subseteq [1,N]^t$.
\end{wredmain}

Writing
$$
\E_{m \in \Z^d \cap K'} 
\prod_{i=1}^t r'_{f_i,c_i(a)}(\tilde \psi(m))
= \E_{m \in \Z^d \cap K'} \prod_{i=1}^t
\Big(\Big(r'_{f_i,c_i(a)}(\tilde \psi(m))-1 \Big)+1\Big)
$$
and multiplying out, we obtain a constant term $1$ and all other terms are
of a form the generalised von Neumann theorem applies to, provided we can
show that
$$\|r'_{f_i,c_i(a)} - 1 \|_{U^{t-1}} = o(1)$$ 
for all $i \in [t]$.
By the inverse theorem, it suffices to show that
$$|\E_{n\in[N/\W]} (r'_{f_i,c_i(a)}(n) - 1) F(g(n)\Gamma)| 
= o_{G/\Gamma,t}(1)$$ 
for all $(t-2)$-step nilsequences $(g(n)\Gamma)_{n \leq N/\W}$ and
$1$-bounded Lipschitz functions $F$.
This task will be carried out in the sections 
\ref{reduction_to_torus}--\ref{conclusion-section}.

\section{Non-correlation with nilsequences} \label{non-corr-intro}

The so far standard line of attack to obtain a result of the form
`the function $h$ does not correlate with $k$-step nilsequences' is to
employ the Green-Tao factorisation theorem
\cite[1.19]{green-tao-polynomialorbits}, which allows us to reduce this
task to the case where the nilsequence is close to being equidistributed. 
A separate estimate which shows that $h$ does not correlate with
periodic (nil)sequences allows us to further assume that the Lipschitz
function involved has zero mean, that is, $\int_{G/\Gamma} F=0$.
Periodic sequences are regarded as major arcs.
We have already deduced a major arc estimate in Section \ref{major-arcs}.
The remaining case with the strong assumption that the nilsequence
behaves in a very equidistributed way corresponds to the minor arc
analysis of the classical Hardy-Littlewood method, cf.~the discussion in
\cite[\S4]{green-tao-nilmobius}.
The procedure of passing to the equidistributed (minor arc) case is fairly
independent of the individual problem and is completely described in \S2
of \cite{green-tao-nilmobius}. 
Thus, we restrict our attention here to providing the necessary major and
minor arc estimates specific to our problem and only summarise the
procedures from \cite{green-tao-nilmobius} we employ.
 
Our approach to the minor arc estimate is modelled on a strategy one
might choose in the classical setting:
If $\theta$ is a rational that belongs to a suitably chosen notion of
`minor arc', then one obtains an upper bound for the expression 
$$\E_{n \leq N} r_f(n)e(\theta n) =
\frac{1}{N}\sum_{\substack{x,y~:\\f(x,y)\leq N}}
 e(\theta f(x,y))= \frac{1}{N}\sum_{\substack{x,y~:\\f(x,y)\leq N}}
e(\theta (ax^2+bxy+cy^2))$$ 
by splitting into suitable summation ranges, fixing either $x$ or $y$, and
applying Weyl's inequality\footnote{See the next Section for more
details.}.
Thus, in our case, we aim to employ the quadratic structure of the form
$f$ by means of Weyl's inequality in order to deduce the estimate
$$\E_{n\leq N} r_f(n)F(g(n)\Gamma) = o(1)$$
for sufficiently equidistributed sequences $(g(n)\Gamma)_{n\leq N}$.
When working with a sequence $(F(g(n)\Gamma))_{n\in [N]}$ directly, Weyl's
differencing trick may only be employed locally on so called generalised
Bohr neighbourhoods, where one can make the locally polynomial structure
of a nilsequence explicit, cf.~the approach in
\cite{green-tao-quadraticmöbius}.

The crucial fact that makes Weyl's  differencing trick work for
exponential sums is the fact that the exponential function is a group
homomorphism.
Since $F$ is a Lipschitz function, one expects it to have a good,
i.e.~short, Fourier approximation. 
In general, elements of a Fourier basis in the non-abelian case arise from
characters, i.e.~homomorphisms. Thus there is a good chance that it is
possible to employ Weyl's inequality globally for elements of the Fourier
basis and hence for a short Fourier approximation of a Lipschitz
function.

In our case, the situation is considerably simplified by the availability 
of a complete quantitative equidistribution theory for polynomial orbits
on nilmanifolds, which has been worked out by Green and Tao in
\cite{green-tao-polynomialorbits}. 
In particular, their generalisation of Leon Green's theorem (`Quantitative
Leibman theorem' \cite[Thm.1.16]{green-tao-polynomialorbits}) asserts that
any polynomial sequences on a nilmanifold $G/\Gamma$ is
$\delta$-equidistributed\footnote{The quantitative notion of
equidistribution is recalled in Section \ref{reduction_to_torus}.}
\emph{if and only if} its projection on the horizontal torus is
$\delta'$-equidistributed, where the dependence is polynomial. 
The horizontal torus bears the advantage of being isomorphic to an
ordinary torus $\R^{d_{\mathrm{ab}}}/\Z^{d_{\mathrm{ab}}}$.
Consequently, we need not consider the representation theory on nilpotent
Lie groups and their homogeneous spaces;
analysing the projected sequence on the horizontal torus by standard
Fourier analysis, or even the quantitative version of Weyl's
equidistribution theory, is sufficient.
(The latter theory will actually reduce matters to looking at sequences
$\Z \to \R/\Z$ arising from horizontal characters.)

Our strategy, after reducing to the equidistributed case, is the
following: 
Let $P$ denote a polynomial of degree $d$. 
Then equidistribution of $(g(n)\Gamma)_{n\leq N}$ on $G/\Gamma$ implies
that $(\pi \circ g(n))_{n \leq N}$ is equidistributed on the horizontal
torus, which implies, as a consequence of Weyl's equidistribution theory,
that $(\pi \circ g(P(n)))_{n \leq N^{1/d}}$ is equidistributed on the
horizontal torus, which implies that $(g(P(n))\Gamma)_{n \leq N^{1/d}}$ is
equidistributed on $G/\Gamma$. 
The distribution of polynomial subsequences was not considered in
\cite{green-tao-polynomialorbits}, but will follow from results of that
paper.
These results will be proved in sections \ref{reduction_to_torus} and
\ref{poly_subs_section}.

For the above strategy to work, a strong major arc analysis is required,
because the $W$-trick introduces very large coefficients into the
quadratic forms under consideration.
For the major arc analysis, we rely on the observation that all of these
large coefficients turn out to be entirely composed of small prime
factors.
We briefly describe in the next section how this information is used to
choose major and minor arcs in the classical setting.
The general case will be carried out in Section \ref{poly_subs_section}
(especially Corollary \ref{modified_abelian_equid} and Proposition
\ref{modified_general_equid} which deal with polynomial subsequences that
have large but smooth coefficients) and Section \ref{factorisation},
which provides a factorisation of polynomial sequences into major and
minor arcs.

\section{A special choice of major and minor arcs is necessary}
In this section we describe briefly and solely for motivational purposes
how the major and minor arcs are chosen in the model case of
correlation with linear phase functions $e(\theta n)$ instead of general
nilsequences.
Here the task is to show that
$$\E_{n\leq N} (r'_{f,\beta}(n) - 1)e(\theta n) = o(1)~.$$
In Section \ref{major-arcs}, we saw that $(r'_{f,\beta} - 1)$ does not
correlate with any $q$-periodic function of $w(N)$-smooth period $q$,
provided $N/q$ is still quite large.
It is therefore possible to choose the major arcs to consist of all
rationals $\theta \in [0,1)$ that are close to a rational with
$w(N)$-smooth denominator: in that case $e(n\theta)$ is close to a
periodic function with $w(N)$-smooth period.
The minor arcs then comprise all $\theta$ that are not close to
rationals with $w(N)$-smooth denominators.
For such a `minor arc' $\theta$, we automatically have 
$\E_{n\leq N} e(\theta n) = o(1)$.

Thus, we define the major arcs to be
$$\mathfrak{M}:= \bigcup_{q \in \mathfrak Q} \mathfrak{M}_q~,$$
where $\mathfrak Q$ is the following set of all not too large
$w(N)$-smooth denominators 
$$\mathfrak Q:=\{1\leq q \leq N^{\eps} : p| q \implies p \leq w(N) \}$$
and where $\mathfrak{M}_q$ is the set of real numbers that are well
approximated by some rational with $w(N)$-smooth denominator:
$$\mathfrak{M}_q := \Big \{ \theta: \Big|\theta - \frac{\alpha}{q} \Big|
\leq \frac{1}{qN^{1-\eps}} \text{ for some } (\alpha,q)=1 \Big \}~.$$
 
The reason behind this choice of major arc is the following.
When we pass to $W$-tricked versions of the representation function, which
are up to normalisation of the form $n \mapsto r_f(\W n + \beta)$, then
this restriction to a \emph{linear} substructure cannot directly be
expressed by the quadratic form $f$. 
For the minor arcs treatment, we, however, hope to work with the quadratic
form directly.
We will therefore consider all choices $(x',y')\in [\W]^2$ such that 
$f(x',y') \equiv \beta \Mod{\W}$ and consider for each choice the
quadratic form $f(\W x + x', \W y + y')$ in $x,y$.
Fixing either $x$ or $y$, we hope to apply Weyl's inequality when 
$\theta \not\in \mathfrak{M}$ to estimate
\begin{align*}
&\sum_{n \leq (N-\beta)/\W} 
  r(\W n + \beta) e(\theta n) \\
&= \sum_{\substack{x',y'\in [\W] :\\ f(x',y') \equiv \beta \Mod{\W}}}
  \sum_{\substack{x,y \\ f(\W x + x', \W y + y') \leq N}} 
  e\Big(\frac{\theta (f(\W x + x', \W y + y') - \beta)}{\W}\Big)~.
\end{align*}
Here we obtain for fixed $x',y'$ and either fixed $x$ or fixed $y$ a
quadratic inside the exponential with leading coefficient $\theta \W a$
or $\theta \W c$ where $a$ and $c$ are coefficients of $f$. 
For the application of Weyl's inequality, we require that this leading
coefficient is close to a rational with large denominator.

Since $ac\W \ll N^{o(1)}$, the choice of major and minor arcs guarantees
that, when $\theta \not\in \mathfrak M$, i.e.
$$\Big|\theta - \frac{\alpha}{q} \Big| \leq \frac{1}{qN^{1-\eps}}$$
for some $q$ that has a prime factor $>w(N)$, or satisfies $q>N^{\eps}$,
then 
$$\Big|a\W \theta - \frac{\alpha'}{q'} \Big| \leq
\frac{1}{qN^{1-\eps-o(1)}}~,$$
where $q'$ has a prime factor $>w(N)$, or satisfies $q'>N^{\eps-o(1)}$.
Thus, $a\W \theta$ can still be thought of as minor arc, when
replacing $N$ by $N^{1-o(1)}$.

\section{A brief overview of the concepts around nilsequences}
Let $G$ be a connected, simply connected, $k$-step nilpotent Lie group,
and let $\Gamma$ be a discrete co-compact subgroup.
Then $G/\Gamma$ is called a $k$-step nilmanifold.
A filtration $G_{\bullet}$ of $G$ is a sequence of subgroups
$$G=G_0=G_1 \geq G_2 \geq \dots \geq G_d \geq G_{d+1} =
\{\mathrm{id}_G\}$$
such that for any $d\geq i,j \geq 0$ the commutator group $[G_i,G_j]$ is
a subgroup of $G_{i+j}$.
The filtration is said to have degree $d$, if $G_{d+1}$ is the first
element in the sequence that is trivial.
By definition, a nilpotent group always has a filtration.  

The quantitative analysis carried out in \cite{green-tao-polynomialorbits}
relies on the existence of a certain type of basis, a Mal'cev basis, for
the Lie algebra $\mathfrak g$ of G.
Adapted to any filtration, there exists a Mal'cev basis for $\mathfrak g$
that parametrises via the exponential map both the groups in the
filtration and the uniform subgroup $\Gamma$ in a very natural way. 
For each such basis $\mathcal X$, Green and Tao introduce a metric
$d_{\mathcal X}$ for $G$ and its quotient $G/\Gamma$ in
\cite[Def.2.2]{green-tao-polynomialorbits}, which then allows them to
define Lipschitz functions on $G/\Gamma$, and also allows them to
introduce a notion of slowly varying (or smooth) sequences 
$(\eps(n))_{n\in \Z}$ that take values in $G$.
Despite the fact that any of the statements on nilsequences require a
fixed choice of Mal'cev basis $\mathcal X$ and corresponding metric
$d_{\mathcal X}$, we will not need to directly work with any of the
specific properties of either of these objects: they will only implicitly
be present through the results from \cite{green-tao-polynomialorbits} we
build on.
For this reason, we content ourselves to refer to \cite[\S2
and App.A]{green-tao-polynomialorbits} for background and exact
definitions.

\begin{definition}[Polynomial sequence; Def.1.8
\cite{green-tao-polynomialorbits}]
 Let $g:\Z \to G$ be a $G$-valued sequence, and define the discrete
derivative $\partial_h g(n):= g(n+h)g(n)^{-1}$ for each $h \in \Z$.
Then $g$ is a \emph{polynomial sequence} with coefficients in
$G_{\bullet}$, when for every $i \in \{0,\dots,d+1\}$, and every choice of
$h_1,\dots, h_i \in \Z$ all $i$-th derivatives satisfy
$\partial_{h_i} \dots \partial_{h_1}g(n) \in G_i$.
We write $\mathrm{poly}(\Z,G_{\bullet})$ for all polynomial sequences
adapted to $G_{\bullet}$ and say they are of degree $d$, where $d$ is
the degree of the filtration.
\end{definition}
Two facts about polynomial sequences are of particular importance. 
The first is a theorem of Lazard: \emph{$\mathrm{poly}(\Z,G_{\bullet})$
forms a group}; see \cite[\S6]{green-tao-polynomialorbits} for a proof
and the reference to the original work.
The second important property is a more explicit description of
polynomial sequences. It is shown in \cite{green-tao-polynomialorbits}
(cf.~\S6 and the remarks following Def.1.8) that every polynomial
sequence can be written in the form $g(n)=a_1^{p_1(n)}\dots a_k^{p_k(n)}$,
where $k$ is some integer, $a_1, \dots, a_k \in G$, and $p_1, \dots,
p_k:\Z \to \Z$ are polynomials. 
Observe that, if the sequence $g_i$ defined by $g_i(n)=a_i^{p_i(n)}$
belongs to $\mathrm{poly}(\Z,G_{\bullet})$, then the assertion that the
discrete derivatives of order $d+1$ all equal $\mathrm{id}_G$ directly
translates to $\deg(p_i) \leq d$.
In general the degree of the polynomial sequence $g$ is much larger than
the degrees of the polynomial exponents $p_1, \dots, p_k$ that appear in
the above mentioned representation. 

\begin{definition}[Horizontal torus]
Write $\pi: G \to (G/\Gamma)_{\mathrm{ab}}:=G/([G,G]\Gamma)$ for the
canonical projection of $G$ on the abelianisation of $G/\Gamma$.
$(G/\Gamma)_{\mathrm{ab}}$ is called the \emph{horizontal torus} of $G$.
\end{definition}

We will extensively work with horizontal characters $\eta:G \to \R/\Z$.
These are additive homomorphisms that annihilate $\Gamma$.
Note that when $g$ has degree $d$, that is, when $g$ has coefficients in
a filtration of degree $d$, then the projection $\eta \circ g$ can be
written as an ordinary polynomial of degree at most $d$ taking values in
$\R/\Z$.
 
\cite[Def.2.6]{green-tao-polynomialorbits} defines the notion of the
modulus $|\eta|$ of a horizontal character.
All that is important to us, is that 
$\|\eta\|_{\mathrm{Lip}} \ll |\eta|$.

\section{Reduction from nilmanifolds to the abelian setting}
\label{reduction_to_torus}

In this section we provide the tool for passing from a general
nilmanifold to the abelian setting of the horizontal torus.
We caution, however, that by far the largest amount of the real work
behind these results is hidden in the application of 
\cite[Thm 1.16]{green-tao-polynomialorbits}, while the converse
statements we prove are fairly straightforward.

Integral to all what follows are the two quantitative notions of
equidistribution that were introduced in 
\cite[Def. 1.2]{green-tao-polynomialorbits}:

\begin{definition}[Quantitative equidistribution,
\cite{green-tao-polynomialorbits}]
Let $G/\Gamma$ be a nilmanifold endowed with Haar measure
and let $\delta_1, \delta_2 \in (0,1)$ be parameters. 
A finite sequence $(g(n)\Gamma)_{n\leq N}$ is said to be
$\delta_1$-\emph{equidistributed} if
$$\Big| \E_{n \in [N]} F(g(n)\Gamma) - \int_{G/\Gamma} F\Big|
\leq \delta_1 \|F\|_{\mathrm{Lip}}$$
for all Lipschitz functions $F:G/\Gamma \to \C$
with
$$\|F\|_{\mathrm{Lip}}
:= \|F\|_{\infty} 
 + \sup_{x,y \in G/\Gamma, x\not=y}
   \frac{|F(x)-F(y)|}{d_{G/\Gamma}(x,y)}~.
$$
$(g(n)\Gamma)_{n\leq N}$ is said to be \emph{totally
$\delta_2$-equidistributed} if
$$\Big| \E_{n \in P} F(g(n)\Gamma) - \int_{G/\Gamma} F\Big|
\leq \delta_2 \|F\|_{\mathrm{Lip}}$$
for all Lipschitz functions $F$ as above and all arithmetic progressions
$P\subseteq [N]$ of length $|P| \geq \delta_2 N$.
\end{definition}
For polynomial sequences these two notions of equidistribution are
equivalent in the sense that every totally $\delta_2$-equidistributed
sequence is $\delta_2$-equidistributed, and every
$\delta_1$-equidistributed sequence is totally
$\delta_2(\delta_1)$-equidistributed, where
$\delta_1^{A} \leq \delta_2(\delta_1) \leq \delta_1$
for some $A$ only depending on the degree of the sequence, and
the dimension and step of the nilmanifold. 
(As this observation will not be used later on, a proof is omitted.)

We set out by recalling the quantitative version of Weyl's inequality from
\cite{green-tao-polynomialorbits}, and the notion of smoothness norms
in terms of which this inequality is phrased.

Any polynomial $g:\Z \to \R/\Z$ of degree $\leq d$ has an expansion of the
form
$$g(n) = \alpha_0 + \alpha_1 \binom{n}{1} + \dots + \alpha_d
\binom{n}{d} ~. $$
The \emph{smoothness norm} of $g$ is defined by
$$\|g\|_{C^{\infty}[N]}:= \sup_{1\leq j \leq d} N^j
\|\alpha_j\|_{\R/\Z}~.$$
This norm was introduced in \cite[Def. 2.7]{green-tao-polynomialorbits} as
a measure of slow variation of polynomial sequences on tori. 
Indeed, 
\begin{align} \label{variation-bound}
\|g(n)-g(n-1)\|_{\R/\Z} \ll_d \|g\|_{C^{\infty}[N]}/N 
\end{align}
holds.
For us it will be more convenient to work with the coefficients of the
ordinary representation of $g$.
When 
$g(n)=\beta_d n^d + \beta_{d-1}n^{d-1} + \dots + \beta_0$, then
(cf.~\cite[Lemma 3.2]{green-tao-nilmobius}) there is $q\geq 1$ with
$q=O_d(1)$ such that 
\begin{align} \label{smoothness_and_betas}
\|q \beta_j\|_{\R/\Z} \ll N^{-j} \|g\|_{C^{\infty}[N]}
\end{align}
for $j=1, \dots, d$.
This follows by expressing each $\beta_j$ as a linear combination of
$\alpha_i$.
The coefficients appearing are bounded by $O_d(1)$.

In the other direction we can show
\begin{align} \label{smoothness_and_betas_2}
 \|g\|_{C^{\infty}[N]} 
\ll \sup_{1 \leq j \leq d} N^{j}
  \|j! \beta_j \|_{\R/\Z}~.
\end{align}
Indeed, $j! \beta_j$ is a linear combination of $\alpha_i$,
$i \geq j$, where the coefficient of $\alpha_j$ is $1$ and all other
coefficients are $O_d(1)$.
Let $j_0$ be the maximal index for which 
$\|g\|_{C^{\infty}[N]} = N^{j_0}\|\alpha_{j_0}\|$.
Then $N^{i} \|\alpha_{i}\| < N^{j_0}\|\alpha_{j_0}\|$
for all $i>j_0$.
Thus $\|\alpha_{i}\| < N^{j_0 - i}\|\alpha_{j_0}\|$.
Then 
$\|j_0! \beta_{j_0} \|
= \| \alpha_{j_0} \| (1 + O_d(N^{-1}))$,
which proves the result.

Part (a) of the following is Green and Tao's Proposition 4.3 from
\cite{green-tao-polynomialorbits}.
While the latter is quite a deep result, its converse, which we prove as
part (b), is rather straightforward.
\begin{proposition}[Weyl]\label{weyl}
$(a)\quad$ Suppose that $g: \Z \to \R$ is a polynomial of degree $d$, and
let $0<\delta < 1/2$. If $(g(n) \Mod{\Z})_{n\in [N]}$ is not
$\delta$-equidistributed in $\R/\Z$, then there is an integer $k$, 
$1\leq k \ll \delta^{-O_d(1)}$ such that 
$\|kg\|_{C^{\infty}[N]} \ll \delta^{-O_d(1)}$.

\noindent $(b)\quad$ 
Suppose that the parameter $\delta=\delta(N) \in (0,1)$ satisfies 
$\delta^{-t} \ll_t N$ for all $t \in \N$.
Further, suppose there are positive integers $k_1, \dots, k_d$ satisfying 
$k_j \ll \delta^{-2^{d-j}}$ such
that 
$$\|k_j \alpha_j\|_{\R/\Z} \leq \delta^{-2^{d-j}}/N^j~.$$
Then, provided $N$ is large enough, there is some positive integer
$A = O_d(1)$ such that $(g(n) \Mod{\Z})_{n\in [N]}$ is not totally
$\delta^{A}$-equidistributed in $\R/\Z$. 
\end{proposition}
\begin{remarks}
 \emph{(1)} The precise choice of exponents in the bounds
$\delta^{-2^{d-j}}$ is not important to this result, but we will later
make use of the fact that
this way
$k_d k_{d-1} \dots k_{d-j} \ll \delta^{-2^{j+1}+1}$.

\emph{(2)} In part $(b)$, the conditions 
$\|k_j \alpha_j\|_{\R/\Z} \leq \delta^{-2^{d-j}}/N^j$
can be replaced by conditions of the form
$\|k_j \beta_j\|_{\R/\Z} \leq \delta^{-2^{d-j}}/N^j$ as they imply
$\|k_j j! \beta_j\|_{\R/\Z} \ll \delta^{-2^{d-j}}/N^j$.
\end{remarks}
\begin{proof} All that is left is to prove part $(b)$. 
Put $k=\mathrm{lcm}(k_1,\dots, k_d)$. 
Then, by the assumption on $\delta$, 
$$\|k \alpha_j\|_{\R/\Z} \leq \delta^{-A'}/N^j = o(1)$$
for some $A'=O_d(1)$ for each $j \in [d]$. 
Consider the sequence
$$(g(kn))_{n\in[\frac{N}{k}\delta^{2 A'}]}~.$$
By \eqref{variation-bound}, each $g(kn)$ in the range satisfies
$\|g(k) - g(kn)\|\ll \delta^{A'}$.
Thus $e\circ g= \exp(2\pi i g(\cdot))$ is almost constant on this range
and we obtain for $N$
sufficiently large
$$\left| \E_{n\in[\frac{N}{k}\delta^{2A'}]} e(g(kn))
-\int_{\R/\Z} e(x)~dx \right| 
\geq 1 - \left(\frac{2\pi\delta^{-A'} \delta^{2A'} N}{N} \right)^2 
\gg \delta^{A'}\|e\|_{\mathrm{Lip}}~,$$
that is, $(g(n) \Mod{\Z})_{n \leq N}$ is not totally
$\delta^{2A'}/k=\delta^{O_d(1)}$-equidistributed.
\end{proof}

The equidistribution of nilsequences is related to the equidistribution
of certain polynomial sequences via the following projection theorem.

\begin{proposition}[Green-Tao `Quantitative Leibman theorem']
\label{prop-projection}
Let $m,d,N$ be positive integers, and let $\delta \in (0,1/2)$ be a
parameter.
Let $G/\Gamma$ be an $m$-dimensional nilmanifold together with a
filtration $G_{\bullet}$ of degree $d$ and a $\delta^{-1}$-rational
Mal'cev basis adapted to this filtration. 
Suppose that $g \in \mathrm{poly}(\Z,G_{\bullet})$.
Then there are positive constants $B$ and $B'$, only depending on $m$ and
$d$, such that the following holds.
If $(g(n)\Gamma)_{n\leq N}$ is not totally $\delta$-equidistributed in
$G/\Gamma$, then there is a non-trivial horizontal character 
$\eta$ of modulus $|\eta|\ll \delta^{-O_{m,d}(1)}$ such that 
$(\eta \circ g(n))_{n\leq N}$ is not totally $\delta^{B}$-equidistributed
in $\R/\Z$.

Conversely, if there is a non-trivial horizontal character $\eta$ of
modulus $|\eta| \ll \delta^{-1}$ such that $(\eta \circ g(n))_{n\leq N}$
fails to be totally $\delta$-equidistributed in $\R/\Z$, then
$(g(n)\Gamma)_{n\leq N}$ is not totally $\delta^{B'}$-equidistributed in
$G/\Gamma$.
\end{proposition}
\begin{proof} 
If $(g(n)\Gamma)_{n\leq N}$ is not totally $\delta$-equidistributed, then
there is a progression  $P = \{ p_0, p_0 + q, \dots, p_0 + \ell q \}$
of length at least $\delta N$ such that the sequence 
$(g(n)\Gamma)_{n \in P}$ fails to be $\delta$-equidistributed.
Define $g'\in \mathrm{poly}(\Z,G_{\bullet})$ by $g'(n):=g(qn + p_0)$.
Then \cite[Thm 2.9]{green-tao-polynomialorbits} implies
that there is a non-trivial horizontal character $\eta: G \to \R/\Z$ of
modulus $|\eta| \ll \delta^{-O_{m,d}(1)}$ such that
$\|\eta \circ g' \|_{C^{\infty}[\delta N]} \ll \delta^{-O_{m,d}(1)}$.
By Proposition \ref{weyl}(b), this implies that $\eta \circ g$ fails to be
totally $\delta^{B}$-equidistributed for some $B=O_{m,d}(1)$.

In the other direction, if there is a non-trivial horizontal $\eta$ of
modulus bounded by $\delta^{-1}$ such that $(\eta \circ g(n))_{n\leq N}$ 
fails to be totally $\delta$-equidistributed, then we again find a
progression  $P = \{ p_0, p_0 + q, \dots, p_0 + \ell q \}$
of length at least $\delta N$ such that the sequence 
$(\eta \circ g(n)\Gamma)_{n \in P}$ fails to be $\delta$-equidistributed.
By Proposition \ref{weyl}(a) we have
\begin{align} \label{slow-var}
\| \eta \circ g(p_0 + jq) - \eta \circ g(p_0 + (j-1)q )\|_{\R/\Z} \ll
\delta^{-O_{m,d}(1)}/N
\end{align}
for all $j \in \{1, \dots, \ell \}$. 
Since $\eta$ is an additive character on a compact group, we have
$\int_{(G/\Gamma)_{\mathrm{ab}}} e(\eta (x)) ~dx = 0$. 
Consider the subprogression 
$P' = \{ p_0, p_0 + q, \dots, p_0 + \ell'q \} \subset P$, 
where $\ell' = \delta^{B'} N$, with $B'=O_{m,d}(1)$ large enough so that
\eqref{slow-var} guarantees
$$
\| \eta \circ g (p_0) - \eta \circ g(p_0+jq)
\|_{\R/\Z} 
\leq \frac{|P'|}{N \delta^{O_{m,d}(1)}}
\leq \frac{1}{4 \pi}$$ 
for all $j, 0 \leq j \leq \ell'$.
This implies
\begin{align*}
\left| \E_{n \in P'} e(\eta \circ g(n)) -
\int_{G/\Gamma} e(\eta(x))~dx  \right| 
= \left| \E_{n \in P'} e(\eta \circ g(n)) \right|
> \frac{1}{2}~,
\end{align*}
using the fact that $\Re (e(x)) = \cos (2 \pi x) \geq 1-(2\pi x)^2 >
\frac{1}{2}$ for
$x\leq \frac{1}{4}$. 

Since $\|e \circ \eta\|_{\mathrm{Lip}(G/\Gamma)} \ll 
\|e\|_{\mathrm{Lip}(\R/\Z)} \|\eta\|_{\mathrm{Lip}(G/\Gamma)}
\ll \delta^{-O_{m,d}(1)}$, where the bound on the Lipschitz constant
of $\eta$ comes from the bound on the modulus
(cf.~\cite[Def.2.6]{green-tao-polynomialorbits}) of the
character, we may in fact choose $B'=O_{m,d}(1)$ large enough to ensure
that also $$\frac{1}{2} > \delta^{B'} \|e \circ \eta\|_{\mathrm{Lip}}$$
holds.
Thus, $(g(n))_{n \leq N}$ is not totally
$\delta^{B'}$-equidistributed in $G/\Gamma$.
\end{proof}

\section{Equidistribution of polynomial subsequences via Weyl's
inequality} \label{poly_subs_section}

With the help of the quantitative Leibman theorem (Proposition
\ref{prop-projection}), which reduces questions about the equidistribution
of polynomial nilsequences to questions about the equidistribution of
polynomials taking values in $\R/\Z$, we analyse in this section the
distribution of polynomial subsequences of polynomial orbits.

The first result states that on the torus polynomial subsequences of
$\delta$-equidistributed sequences are equidistributed too. 
Before stating this proposition properly we give an informal description
of its contents here.
A polynomial $g:\Z \to \R/\Z$ is equidistributed if and only if one of its
coefficients is irrational.
Quantitative equidistribution is an assertion on whether or not there is
a Lipschitz function $F:\R/\Z \to \C$ for which  
$|\E_{n\leq N}F(g(n)) - \int_{\R/\Z} F|$ fails to be small.
Approximating the Lipschitz function $F$ by a Fourier series, one
sees that studying this quantity is equivalent to studying the exponential
sums $\E_{n\leq N} e(\omega g(n))$ for certain rational $\omega$. 
The latter is naturally approached by Weyl's inequality which then shows
that the quantitative equidistribution of $g$ is an assertion about
whether or not there is a coefficient of $g$ that is not close to a
rational with small denominator.
This rational approximation property is preserved when we consider
compositions $g \circ P$ of $g$ with an integral polynomial $P$ whose
leading coefficient is not too large. 
To see this we only need to consider the case where $g$ has a `highly
irrational' coefficient. 
Take the largest-index coefficient of $g$ which is `highly
irrational' and call it $\beta_{i_0}$. Then we may check that the
largest-index coefficient of $g \circ P$ which arises from the highly
irrational coefficient $\beta_{i_0}$ of $g$ is still considerably
irrational. (Some bounds on the lower coefficients of $P$ are needed in
order to avoid cancellation.)
 
\begin{proposition}[Equidistribution of polynomial subsequences:
Abelian case]\label{abelian equid subseqs}
 Suppose that $g: \Z \to \R$ is a polynomial of degree $d$ and that 
$P(n)=\sum_{i=0}^{d'} \gamma_i n^i$ is a polynomial with integer
coefficients of degree $d'$ such that the leading coefficient
$\gamma_{d'}$ is bounded by $L_0$, while all other coefficients satisfy
the inequality $\gamma_i \leq N^{(d'-i)/d'}$. 
Let $0<\delta < 1/2$ and suppose $\delta^{-t} \ll_t N$ for all $t \in \N$.
Then there is some integer $A = O_d(1)$ such that when 
$(g(n)\Mod{\Z})_{n\in [N]}$ is totally $\delta$-equidistributed
and when $L_0 \leq \delta^{-1/A}$, then 
$(g\circ P (n) \Mod{\Z})_{n\in [N^{1/d'}]}$ is totally
$\delta^{1/O_{d,d'}(1)}$-equidistributed.
\end{proposition}

\begin{proof}
Since $g$ is totally $\delta$-equidistributed, Proposition \ref{weyl}(b)
implies that there is an integer $A'=O_d(1)$ such that no $d$-tuple of
positive integers $k_1,\dots,k_d$ satisfies simultaneously $k_j \ll
\delta^{-2^{d-j}/A'}$ and
$\|k_j \beta_j\| \ll \delta^{-2^{d-j}/A'}/N^j$ for all $j=1,\dots,d$.
We deduce that there is some index $i_0$ among them such that
$\|k_{i_0} \beta_{i_0}\| \ll \delta^{-2^{d-i_0}/A'}N^{-i_0}$ 
does not hold for any $k_{i_0} \ll \delta^{-2^{d-i_0}/A'}$.
Suppose $i_0$ is maximal with this property.
Then for all $\ell$ with $i_0< \ell \leq d$ we find 
$\kappa_{\ell} \ll \delta^{-2^{d-\ell}/A'}$ such that 
\begin{align}\label{kappa-ell}
\|\kappa_{\ell} \beta_{\ell} \| 
\leq \delta^{-2^{d-\ell}/A'} N^{-\ell} ~.
\end{align}
For any $j \in \{1, \dots, d\}$, considering the $j$th term of
\begin{align*}
 \sum_{j=0}^d \beta_j (P(n))^{j} = g \circ P (n) 
\end{align*}
we have
\begin{align*}
  \beta_j (P(n))^j 
= \beta_j (\gamma_{d'})^j 
  n^{j d'}
  + \beta_j Q_j(n)~,
\end{align*}
where $Q_j(n)$ is a polynomial of degree $\leq j d' - 1$ such that the
coefficient of $n^i$ for any $i$ is bounded by 
$O_{d,d'}(N^{j}N^{-i/d'} \delta^{-j/A})$
since 
$$(P(n))^j
=\Big(\sum_{t=1}^{d'} \gamma_t n^t\Big)^j 
= \sum_{(t_1, \dots, t_j) \in [d']^j} 
  \gamma_{t_1} \dots \gamma_{t_j}
  n^{t_1+ \dots + t_j}$$
and 
$$\gamma_{t_1} \dots \gamma_{t_j} \leq L_0^j N^{j-(t_1+\dots+t_j)/d'}~.$$ 
Define $\sigma_i$, $i=0, \dots, dd'$, to be the following coefficients
$$ \sum_{i=0}^{dd'} \sigma_i n^i = g \circ P (n) = \sum_{j=0}^d \beta_j
(P(n))^{j}~.$$
Comparing coefficients, each $\sigma_i$ may be written as a linear
combination of $\beta_j$ with $j\geq i/d'$;
$\sigma_{jd'}$ is the $\sigma$-coefficient of largest index whose
representation in terms of $\beta$'s contains $\beta_j$, which
appears with coefficient $(\gamma_{d'})^j$ in the representation. 

Next, we aim to show that there is $A''=O_{d,d'}(1)$ such that every
choice of $k_1, \dots, k_{dd'}$ with $k_j \leq \delta^{-2^{dd'-j}/A''}$
for each $j \in \{1,\dots, dd'\}$ contains some $k_{j_0}$ such that
$$ \|k_{j_0} \sigma_{j_0} \| 
> \delta^{-2^{dd'-j_0}/A''} N^{-{j_0}/d'}~.$$
This, when applied with $k_j=qk$ for any $k \leq \delta^{-1/A''}$, would
in view of \eqref{smoothness_and_betas} imply 
$\|k_j g \circ P\|_{C^{\infty}[N^{1/d'}]} \gg \delta^{-1/A''}$,
from which the result follows by Proposition \ref{weyl}(a).

We will show that we can pick $j_0=i_0 d'$.
Thus, suppose for contradiction that 
\begin{align}\label{for-contr.}
\|k_{i_0d'} \sigma_{i_0d'} \| 
\leq \delta^{-2^{dd'-i_0d'}/A''} N^{-i_0}
\end{align}
holds for some $k_{i_0d'} \leq \delta^{-2^{dd'-i_0d'}/A''}$.
Note that 
$$k_{i_0 d'} \sigma_{i_0 d'} 
= k_{i_0 d'} (\gamma_{d'})^{i_0} \beta_{i_0} + 
\sum_{\ell > i_0} k_{i_0 d'}C_{\ell} \beta_{\ell}~,
$$
where the $C_{\ell}$ are integers of order 
$O_{d,d'}(N^{\ell - i_0}\delta^{-d/A})$ as can be deduced from the
equation
$$C_{\ell} 
= \sum_{\substack{t_1, \dots, t_{\ell} \in [d'] \\
        t_1 +\dots+ t_{\ell}=i_0d'}} 
\gamma_{t_1} \dots \gamma_{t_{\ell}}~.$$
We wish to discard all the terms with $\ell > i_0$ in the above expression
for $k_{i_0 d'} \sigma_{i_0 d'}$ in order to deduce that $\beta_{i_0}$ is
well approximable by rationals which will hopefully lead us to the sought
for contradiction. 
Thus, in view of \eqref{kappa-ell}, we multiply the above expression for
$k_{i_0 d'} \sigma_{i_0 d'}$ by $\kappa:=\prod_{\ell>i_0}\kappa_{\ell}$.
Inequality \eqref{for-contr.} yields
$$ \| \kappa k_{i_0d'} \sigma_{i_0d'}   \| 
\ll \delta^{-2^{dd'-i_0d'}/A''}
    \delta^{-(2^0+\dots+2^{d-i_0-1})/A'} 
  N^{-i_0} 
= \delta^{-2^{dd'-i_0d'}/A''}
  \delta^{-(2^{d-i_0}-1)/A'} 
  N^{-i_0}~.
$$
Writing 
$\bar{k}= k_{i_0d'}(\gamma_{d'})^{i_0}\prod_{\ell>i_0}\kappa_{\ell}$,
we have
$$k_{i_0d'} \sigma_{i_0d'} \prod_{\ell>i_0}\kappa_{\ell}
= \bar k \beta_{i_0} 
+ k_{i_0d'} 
  \sum_{\ell>i_0} \kappa C_{\ell} \beta_{\ell} ~,$$
where in view of \eqref{kappa-ell} and the bound on the $C_{\ell}$
\begin{align*}
\|k_{i_0d'} \kappa C_{\ell} \beta_{\ell} \| 
&\ll \delta^{-2^{dd'-i_0d'}/A''} 
    \delta^{-(2^0+\dots+2^{d-i_0-1})/A'}
    N^{-\ell} N^{\ell - i_0} \delta^{-d/A} \\
&= \delta^{-2^{dd'-i_0d'}/A''} \delta^{-(2^{d-i_0}-1)/A'} \delta^{-d/A}
   N^{- i_0}~.
\end{align*}
Recalling that $\delta^{-t} \ll_t N$ for all $t \in \N$, this upper bound
is seen to be $o(1)$.
Together with the bound on $ \| \kappa k_{i_0d'} \sigma_{i_0d'} \|$ this
allows us to employ the triangle inequality provided $N$ is large enough
that no wrap-around issues can occur.
In particular, this allows us to deduce that
$$\|\bar k \beta_{i_0} \| 
\ll \delta^{-2^{dd'-i_0d'}/A''} \delta^{-(2^{d-i_0}-1)/A'} \delta^{-d/A}
   N^{-i_0}~.$$
Choosing $A''= 2^{dd'+1}A'$ and $A=2dA'$ 
(to ensure that $L_0^{j} \leq \delta^{-d/A} \leq \delta^{-1/(2A')}$)
this translates to
$$\|\bar k \beta_{i_0} \| 
\ll \delta^{-2^{d-i_0}/A'} N^{-i_0}~,$$
while we obtain the following bound on $\bar k$
$$\bar k \leq \delta^{-2^{dd'-i_0d'}/A''}
L_0^{i_{0}}\delta^{-2^{d-i_0-1}/A'} \leq \delta^{-2^{d-i_0}/A'}~.
$$
Hence, we obtained a contradiction to the rational non-approximability
properties of $\beta_{i_0}$.
\end{proof}

Next, we slightly extend this result. 
Consider the binary quadratic form $f(x,y)$ for fixed $y$ and its
restriction to subprogressions modulo $q$ in the $x$ variable:
$$f(qx+r,y) 
= aq^2 x^2  + x(2aqr + bqy) + (ar^2 + bry + cy^2)~.$$
This defines a quadratic polynomial 
$P(x):= \gamma_2 x^2 + \gamma_1 x + \gamma_0:= f(qx+r,y)$
in $x$.
Being interested in $(x,y)$ such that $f(x,y) \leq N$, we may suppose that
$y \ll N^{1/2}$. 
Further assume that $q$ is $k$-smooth (we will be interested in
the case $q=\W$) and satisfies $q \ll N^{o(1)}$. 
Then the coefficients of this quadratic polynomial in $x$ have the
following properties.
$\gamma_2$ is $k$-smooth, 
and $\gamma_2 \ll q^2 N^{(2-2)/2}$,
$\gamma_1 \ll q^1N^{(2-1)/2}$,
$\gamma_0 \ll q^0N^{(2-0)/2}$.

The proposition below is tailored to address polynomials with these
specific properties.
\begin{proposition}\label{modified_abelian_equid}
Let $0 < \delta < 1/2$ and suppose $\delta^{-t} \ll_t N$ for all $t \in
\N$.
Let $k$ be a positive integer and suppose that the polynomial sequence
$g:\Z \to \R$, $g(n) = \sum_{j=0}^d \beta_j n^j$ has the property that
for every $k$-smooth integer $q$, $q \leq N^{o(1)}$, for every choice of 
$k_1, \dots, k_d$ with $0 < k_j \leq \delta^{-2^{d-j}}$ for
$j=1,\dots,d$, and for every sufficiently large $N$, we have
$$\sup_{1 \leq j \leq d} \|q^j k_j \beta_j \| \delta^{2^{d-j}} (N/q)^{j} 
\geq 1~.$$
Then, if $P(n)=\sum_{i=0}^{d'} \gamma_i n^i$ is an integer-coefficient
polynomial of degree $d'$ whose leading coefficient is a \emph{$k$-smooth}
integer satisfying $\gamma_{d'} < N^{o(1)}$, while all other coefficients
satisfy the inequality $\gamma_i \leq N^{(d'-i)/d'} {\gamma_{d'}}^{i/d'}$,
we obtain a conclusion similar to the one in the previous proposition:

Then there is a $k$-smooth number $\tilde q$, $\tilde q \ll N^{o(1)}$,
such that each of the sequences 
$(g\circ P (\tilde q n + r) (\mathrm{mod}\Z))_{n \in
[(N/\gamma_{d'}\tilde q^{d'})^{1/d'}]}$ for $r \in [\tilde q]$ is totally
$\delta^{1/O_{d,d'}(1)}$-equidistributed, provided $N$ is large enough.
\end{proposition}

\begin{remark}
 The unconventional form of the inapproximability conditions imposed on
the $\beta_i$ comes out of our choice of major and minor arcs; cf.
Proposition \ref{modified_general_equid} and the next section.
\end{remark}

\begin{proof}
 Consider $q:= \gamma_{d'}$, and let, as in the previous proof,
$i_0$ be the maximal index for which 
$$\|q^{i_0} k_{i_0} \beta_{i_0} \| \delta^{2^{d-{i_0}}} (N/q)^{i_0}
\gg 1
$$
for all $k_{i_0} \leq \delta^{-2^{d-i_0}}$.
Thus, for $\ell > i_{0}$ there are 
$\kappa_{\ell} \leq \delta^{-2^{d-\ell}}$ such that
\begin{align} \label{kappa--l>i_0}
\|q^{\ell} \kappa_{\ell} \beta_{\ell} \| 
\ll \delta^{-2^{d-\ell}}(N/q)^{-\ell}~. 
\end{align}
We wish to employ this information to proceed as in the previous proof,
that is, we wish to assume for contradiction that all coefficients of $g
\circ P$ are close to rationals. In particular this would apply to the
$(i_0d')$-th coefficient.
Writing that coefficient as a linear combination of $\beta$'s we would
then like to deduce that $\beta_{i_0}$ has to be close to a rational,
which produces a contradiction. 
Unfortunately, the above information is not quite sufficient for our
purposes yet: we require similar bounds on 
$\|\kappa_{\ell} \beta_{\ell}\|$ instead of on 
$ \|q^{\ell} \kappa_{\ell} \beta_{\ell} \|$.
To work around this, we pass to higher powers $q^t$ of $q$, aiming to find
a small $t$ and an index $i_t$ such that 
$\|q^{t d i_t} k_{i_t} \beta_{i_t} \| \gg \delta^{-2^{d-{i_t}}}
(N/q^{dt})^{-i_t}$, while
$\|q^{t \ell} k_{\ell} \beta_{\ell} \| \ll \delta^{-2^{d-{\ell}}}
(N/q^t)^{-\ell}$ for $\ell > i_t$.
The gap between $q^{t d i_t}$ and $q^{t \ell}$ (for $\ell > i_t$)
introduced by the extra factor $d$ will be sufficient to analyse $g \circ
P$ on subprogressions modulo $q^{td/d'}$.

Returning to the proof, note that \eqref{kappa--l>i_0} implies 
$$\|q^{t\ell} \kappa_{\ell} \beta_{\ell} \|
\leq q^{(t-1)\ell}\|q^{\ell} \kappa_{\ell} \beta_{\ell} \|
\ll \delta^{-2^{d-\ell}}(N/q^t)^{-\ell} $$
for $\ell > i_{0}$ and for all positive integers $t$.
By assumption on the rationality properties of the $\beta_j$,
$j=1, \dots, d$, there is an index $i_1$, which by the previous
observation necessarily satisfies $i_1 \leq i_0$, such that
$$\|q^{2 i_1} k_{i_1} \beta_{i_1} \| \delta^{2^{d-{i_1}}}
(N/q^{2})^{i_1}
\gg 1
$$
for all $k_{i_1} \leq \delta^{-2^{d-i_1}}$.

Proceeding like this, we obtain a decreasing sequence 
$i_0 \geq i_1 \geq i_2 \geq \dots$ of positive integers such that for
every $j$ the following two families of inequalities hold:
$$\|q^{(j+1) i_j} k_{i_j} \beta_{i_j} \| 
\gg \delta^{-2^{d-{i_j}}} (N/q^{(j+1)})^{-i_j}
$$
for all $k_{i_j} \leq \delta^{-2^{d-i_j}}$,
and for every $\ell$ with $d \geq \ell > i_j$ there is 
$\kappa_{\ell} \leq \delta^{-2^{d-\ell}}$ such that
\begin{align*}
\|q^{(j+1)\ell} \kappa_{\ell} \beta_{\ell} \| 
\ll \delta^{-2^{d-\ell}}(N/q^{(j+1)})^{-\ell}~. 
\end{align*}
By positivity of the indices $i_j$, there is $t=O_{d,d'}(1)$ such that
$i_{t-1}=i_t=i_{t'}$
for all $t < t' \leq tdd'$. 
Setting $\tau=tdd'$, we therefore have
\begin{align}\label{beta_i_t rationality}
\|q^{(1+\tau)i_t} k_{i_t} \beta_{i_t} \| \delta^{2^{d-{i_t}}}
(N/q^{1+\tau})^{i_t} \gg 1
\end{align}
for all $k_{i_t} \leq \delta^{-2^{d-i_t}}$,
while we find for every $\ell > i_t$ a positive integer 
$\kappa_{\ell} \leq \delta^{-2^{d-\ell}}$ such that
\begin{align} \label{kappa--l>i_t}
\|q^{t\ell} \kappa_{\ell} \beta_{\ell} \| \ll
\delta^{-2^{d-\ell}}
(N/q^{t})^{-\ell}~.
\end{align}
Now recall that $d'= \deg(P)$ and consider the sequence 
$g(P(q^{\tau/d'}n+r))_{n \in [(N/q^{\tau+1})^{1/d'}]}$ 
for an arbitrary $r \in [q^{\tau/d'}]$.
Defining coefficients $\sigma_i$ by 
$$
\sum_{i=0}^{dd'} \sigma_{i}n^i
= g(P(q^{\tau/d'}n+r))
= \sum_{j=0}^d \beta_j(P(q^{\tau/d'}n+r))^j~,$$
we have
\begin{align}\label{sigma_[i_t d']}
\sigma_{i_td'} 
= \beta_{i_t} (\gamma_{d'} q^{\tau})^{i_t} 
+ \sum_{\ell > i_t} \beta_{\ell} q^{\tau i_t } C_{\ell}~,
\end{align}
with integer coefficients $C_{\ell}$.
We need a bound on $C_{\ell}$ and proceed to show that
$C_{\ell} = O(N^{\ell - i_{t}}q^{i_t})$.
Expanding out products yields
\begin{align*}
&(P(nq^{\tau/d'} + r))^{\ell} 
 =\Big(\sum_{j=1}^{d'} \gamma_j (q^{\tau/d'}n+r)^j\Big)^{\ell}  \\
&= \sum_{(j_1, \dots, j_{\ell}) \in [d']^{\ell}} 
  \gamma_{j_1} \dots \gamma_{j_{\ell}}
  \sum_{\substack{(u_1, \dots, u_{\ell}) \leq \\(j_1, \dots, j_{\ell})}}
  \binom{j_1}{u_1} \dots \binom{j_{\ell}}{u_{\ell}}
  (q^{\tau/d'}n)^{u_1+ \dots + u_{\ell}}
  r^{(j_1 - u_1) + \dots + (j_{\ell} - u_{\ell})}~.
\end{align*}
Consider any term involving 
$(q^{\tau/d'} n)^{i_td'}=q^{\tau i_t}n^{i_td'}$.

If $j_1 +  \dots + j_{\ell} = u_1+ \dots + u_{\ell} = i_t d'$, 
then the coefficient of $q^{\tau i_t}n^{i_td'}$ is
$\gamma_{j_1} \dots \gamma_{j_{\ell}} \leq N^{\ell - i_t}q^{i_t}$.
If $j_1 +  \dots + j_{\ell} > u_1+ \dots + u_{\ell} = i_t d'$, then the
coefficient of $(q^t n)^{i_td'}$ 
is bounded by $O_{d,d'}(\gamma_{j_1} \dots \gamma_{j_{\ell}} r^{\ell d'}) 
= O_{d,d'}(N^{\ell - i_t -(1/d')}q^{\ell}r^{\ell d'}) 
= O(N^{\ell - i_t})$,
since $r < q^{\tau/d'} \ll N^{o(1)}$.
Thus in total, $C_{\ell} = O(N^{\ell - i_{t}}q^{i_t})$.

We return to analysing the rational approximations of the individual
terms of \eqref{sigma_[i_t d']}. 
Notice that $\tau \geq t \ell$ for all $\ell \in [d]$. 
Thus \eqref{kappa--l>i_t} guarantees for $\ell > i_t$ the existence of
$\kappa_{\ell} \leq \delta^{-2^{d-\ell}}$ such that 
\begin{align}\label{ell>i_t terms}
  \|\beta_{\ell} \kappa_{\ell} q^{\tau i_t} C_{\ell}\| 
  \ll \delta^{-2^{d-\ell}}
  N^{-\ell}
  q^{\tau i_t}
  C_{\ell}
  \ll \delta^{-2^{d-\ell}}
  N^{-i_t}
  (q^{\tau+1})^{i_t}
\end{align}
holds.

We are finally in the position to show that there is $A=O_{d,d'}(1)$ such
that $$g(P(q^{\tau/d'}n + r))_{n\in[(N/q^{\tau+1})^{1/d'}]} $$ is totally
$\delta^{1/A}$-equidistributed.
More precisely, we show that there is $A'=O_{d,d'}(1)$ such that for
every $k_{i_td'} \leq \delta^{-2^{dd'-i_td'}/A'}$
$$\| k_{i_td'} \sigma_{i_td'} \| > (N/q^{\tau+1})^{-i_t}
\delta^{-2^{d'(d-i_t)}/A'}$$
holds true. From here the result follows from Proposition \ref{weyl}(a).
Suppose for contradiction that
$$\| k_{i_td'} \sigma_{i_td'} \| \leq (N/q^{\tau+1})^{-i_t}
\delta^{-2^{d'(d-i_t)}/A'}$$
for some $k_{i_td'} \leq \delta^{-2^{dd'-i_td'}/A'}$.
Let $\kappa:=\kappa_d \dots \kappa_{i_t+1}$ (or $\kappa=1$ when the
product is empty), then, since $\kappa_{\ell} \leq \delta^{-2^{d-\ell}}$,
$$\| \kappa k_{i_td'} \sigma_{i_td'} \| 
\ll (N/q^{\tau+1})^{-i_t} \delta^{-2^{d'(d-i_t)}/A'}
\delta^{-(2^{d-i_t}-1)}~.$$
Considering the summands in \eqref{sigma_[i_t d']}, the bounds
\eqref{ell>i_t terms} imply
$$\|\beta_{\ell} \kappa q^{\tau i_t} C_{\ell}\| 
  \ll \delta^{-(2^{d-i_t}-1)} \delta^{-2^{d'(d-i_t)}/A'}
  N^{-i_t} (q^{\tau+1})^{i_t}~.$$
Appealing to the assumptions that both $q$ and $\delta^{-1}$ are bounded
by $N^{o(1)}$, the above is seen to equal $O(N^{-1 + o(1)})= o(1)$ since 
$i_t \geq 1$.
Thus, provided $N$ is large enough, no wrap-around issues appear
when examining the circle norm $\|k_{i_td'} \kappa \sigma_{i_td'}\|$ 
and we find the following statement on rational approximation of
$\beta_{i_t}$
\begin{align*}
\| k_{i_td'} \kappa \beta_{i_t} (\gamma_{d'}q^{\tau})^{i_t} \|
&= \| k_{i_td'} \kappa \beta_{i_t} (q^{\tau+1})^{i_t} \| \\
&= \| k_{i_td'} \kappa \sigma_{i_td'}
  -\sum_{\ell>i_t} \beta_{\ell} \kappa q^{\tau i_t} C_{\ell}\| \\
&\leq |k_{i_td'} \kappa \sigma_{i_td'}| 
 + \sum_{\ell>i_t} |\beta_{\ell} \kappa q^{\tau i_t} C_{\ell}| \\
& \ll \delta^{-(2^{d-i_t}-1)} \delta^{-2^{d'(d-i_t)}/A'}
  N^{-i_t} (q^{\tau+1})^{i_t}~.
\end{align*}
Choosing $A'=2^{d'(d-i_t)}$, this shows that there is $\bar k$, namely 
$\bar k= k_{i_td'} \kappa a^{i_t}$, bounded by $\delta^{-2^{d-i_t}}$ such
that
$$
\| \bar k \beta_{i_t} q^{\tau+1} \| \ll 
\delta^{-2^{d-i_t}} (N/q^{\tau+1})^{-i_t}~,$$
contradicting \eqref{beta_i_t rationality}.
\end{proof}

Combining either of the previous two results with the quantitative Leibman
theorem, the general case of the equidistribution theorem for subsequences
follows. 
\begin{proposition}[Equidistribution of polynomial subsequences]
\label{equid_subsequences}
Let $N,d,d'$ be positive integers, and let $L_0$ and $\delta \in (0,1/2)$
be parameters, and suppose that $\delta^{-t} \ll_t N$ for all $t \in \N$.
Let $g \in \mathrm{poly}(\Z,G_{\bullet})$ be a polynomial sequence of
degree $d$ and suppose that the finite orbit $(g(n)\Gamma)_{n\in [N]}$
is totally $\delta$-equidistributed in $G/\Gamma$. 
Let $P:\Z \to \Z$ be an integer-coefficient polynomial of degree $d'$
whose coefficients are bounded by $L_0$.
Then there is some $A = O_{d,d'}(1)$ such that whenever $L_0^{A} <
\delta$, then the polynomial subsequence 
$((g\circ P)(n))_{n\in [(N/\gamma_{d'})^{1/d'}]}$ is totally
$\delta^{1/O_{d,d'}(1)}$-equidistributed on $G/\Gamma$. 
\end{proposition}
\begin{proof}
We first pass to the abelian setting: 
by Proposition \ref{prop-projection}, there are constants
$A,A'=O_{m,d}(1)$ such that every sequence $(\eta \circ g(n))_{n\in [N]}$
for a horizontal character $\eta$ of modulus at most $\delta^{-A}$ is
totally $\delta^{1/A'}$-equidistributed.
Applying Proposition \ref{abelian equid subseqs}, we deduce that for each
such character $\eta$ the sequence
$(\eta \circ g \circ P(n))_{n \in [(N/\gamma_{d'})^{1/d'}]}$
is totally $\delta^{1/O_{d,d',m}(1)}$-equidistributed in $\R/\Z$.
An application of the other direction of Proposition \ref{prop-projection}
then allows us to return to $G/\Gamma$ and deduce the
stated equidistribution property of 
$(g \circ P(n))_{n \in [(N/\gamma_{d'})^{1/d'}]}$ in $G/\Gamma$.
\end{proof}

Similarly, Proposition \ref{modified_abelian_equid} results in an
assertion for polynomial orbits on general nilsequences:

\begin{proposition}\label{modified_general_equid}
Let $N,d,d',k$ be positive integers, 
and let $\delta \in (0,1/2)$ be such that 
$\delta^{-t} \ll_t N$ for all $t \in \N$.
Let $g \in \mathrm{poly}(\Z,G_{\bullet})$ be a polynomial sequence of
degree $d$ and suppose that for every $k$-smooth number $q$, 
$q \ll N^{o(1)}$, the sequence $(g(qn)\Gamma)_{n\in[N/q]}$ is totally
$\delta$-equidistributed in $G/\Gamma$.

Suppose further that $P:\Z \to \Z$ is an integer-coefficient polynomial of
degree $d'$ as in Proposition \ref{modified_abelian_equid}. 
That is, if $P(n)=\sum_{i=0}^{d'} \gamma_i n^i$, then
$\gamma_{d'}$ is a $k$-smooth integer with $\gamma_{d'} < N^{o(1)}$,
while all other coefficients satisfy the inequality 
$\gamma_i \leq N^{(d'-i)/d'} {\gamma_{d'}}^{i/d'}$.

Then there is a $k$-smooth number $\tilde q$, $\tilde q \ll N^{o(1)}$,
such that each of the sequences 
$(g\circ P (\tilde q n + r) \Gamma)_{n\in [(N/\gamma_{d'} \tilde
q^{d'})^{1/d'}]}$ for $r \in [\tilde q]$ is totally
$\delta^{1/O_{d,d'}(1)}$-equidistributed in $G/\Gamma$, 
provided $N$ is large enough.
\end{proposition}

\begin{proof}
Let $\eta:G/\Gamma \to \R/\Z$ be an arbitrary non-trivial horizontal
character of modulus bounded by $\delta^{-O_{m,d}(1)}$ and suppose that
$\eta \circ g$ has the polynomial representation
$\eta \circ g (n) = \sum_{j=0}^d \beta_j n^j$
in $\R/\Z$.
Let $q$, $q \leq N^{o(1)}$, be $k$-smooth and consider the sequence 
$$(\eta \circ g(qn)\Gamma)_{n \in [N/q]}~.$$
By the equidistribution assumption on the subsequences of $g$, by
Proposition \ref{prop-projection} and by Proposition \ref{weyl}(b), there
is an integer $B=O_{d}(1)$ such that for every choice of 
$k_1, \dots, k_d$ with $0 < k_j \leq \delta^{-2^{d-j}/B}$ for
$j=1,\dots,d$, and for every sufficiently large $N$, we have
$$\sup_{1 \leq j \leq d} \|q^j k_j \beta_j \| \delta^{2^{d-j}/B}
(N/q)^{j} \geq 1~.$$
Thus, with $\delta^{1/B}$ in place of $\delta$, the conditions of
Proposition \ref{modified_abelian_equid} are satisfied and
hence there is $\tilde q \ll N^{o(1)}$ such that for every 
$r \in [\tilde q]$
the sequence
$$(\eta \circ g \circ P (\tilde q n + r) \Mod{\Z})_{ 
n \in [(N/\gamma_{d'} \tilde q)^{1/d'}]}$$ 
is totally $\delta^{1/O_{d,d'}(1)}$-equidistributed in $\R/\Z$, 
provided $N$ is large enough.
An application of Proposition \ref{prop-projection} to get back to
$G/\Gamma$ gives the result.
\end{proof}

\section{The factorisation into minor and major arcs}\label{factorisation}

In view of the previous section, a `minor arc sequences' 
$g \in \mathrm{poly}(\Z,G_{\bullet})$ 
should satisfy the conditions of Proposition \ref{modified_general_equid}
in order to guarantee its applicability.
That is, given $k \in \N$, $\delta = \delta(N) \in (0,1/2)$, and
$R\ll N^{o(1)}$, the sequence $g$ should have the property that for every
$k$-smooth number $q \leq R$ the finite sequence 
$(g(q n)\Gamma)_{n \in [N/q] }$ is $\delta$-equidistributed in
$G/\Gamma$.

In this section we will achieve a factorisation of an arbitrary polynomial
sequence $g$ into a product $\eps g' \gamma$, where $\eps$ is slowly
varying (`smooth'), $\gamma$ is periodic with a $k$-smooth common
difference, and $g'$ has the `minor arc property' described above.
We will ensure that $g'$ satisfies a slightly stronger version of this:
when we restrict $g'$ to subprogressions on which $\gamma$ is constant
and on which $\eps$ is almost constant, then the restricted sequence
still enjoys the `minor arc property'.

This factorisation will be obtained by iteration of the Green-Tao
factorisation theorem \cite[Thm. 1.19]{green-tao-polynomialorbits}
employing its dimension reduction as a guarantee for termination of the
iteration.
Before we state the factorisation theorem, we recall the notion of
smoothness of sequences.

\begin{definition}[$(M,N)$-smooth sequence,
\cite{green-tao-polynomialorbits} Def.1.18]
Let $G/\Gamma$ be a nilmanifold with $Q$-rational Mal'cev basis 
$\mathcal X$ and metric $d=d_{\mathcal X}$.
Let $(\eps(n))_{n\in\Z}$ be a sequence in $G$, and let $M,N \geq 1$.
Then $\eps$ is said to be $(M,N)$-smooth if both 
$d(\eps(n),\mathrm{id}_G)\leq M$ and $d(\eps(n),\eps(n-1))\leq M/N$ are
satisfied for all $n \in [N]$. 
\end{definition}

In the later iteration of the Green-Tao factorisation theorem we will
encounter a product of smooth sequences, which needs to be shown to be
smooth itself.
Notice therefore that, when $(\eps(n))_{n \in \Z}$ is $(M,N)$-smooth and
when $(\eps'(n))_{n \in \Z}$ is $(M,N/q)$-smooth, then the triangle
inequality and right-invariance of the metric $d$ yield 
$$d(\eps(qn+j)\eps'(n),\mathrm{id}_G) \leq
d(\eps(qn+j),\mathrm{id}_G)+d(\eps'(n),\mathrm{id}_G) \leq 2M$$
for all $n \in [N/q]$.
Employing also the approximate left-invariance of $d$ (see
\cite[Lemma A.5]{green-tao-polynomialorbits}), we obtain
$$d(\eps(qn+j)\eps'(n),\eps(q(n-1)+j)\eps'(n-1))
\leq 2 q Q^{O(1)}M/N~.$$
Thus,  $(\eps(qn+j)\eps'(n))_{n\in \Z}$ is $(2Q^{O(1)}M,N/q)$-smooth.

The tool to split into major and minor arcs is the following Green-Tao
factorisation theorem. 

\begin{theorem}[Green-Tao, Thm 1.19 \cite{green-tao-polynomialorbits}]
 Let $m,d \geq 0$, and let $Q_0, N \geq 1$ and $A > 0$ be real numbers.
Suppose that $G/\Gamma$ is an $m$-dimensional nilmanifold together with a
filtration $G_{\bullet}$ of degree $d$. 
Suppose that $\mathcal X$ is a $Q_0$-rational Mal'cev basis $\mathcal X$
adapted to $G_{\bullet}$ and that $g \in \mathrm{poly}(\Z,G_{\bullet})$. 
Then there is an integer $Q$ with 
$Q_0 \leq Q \ll Q_0^{O_{A,m,d}(1)}$, 
a rational subgroup $G' \subseteq G$, a Mal'cev basis $\mathcal X'$ for
$G'/\Gamma'$ in which each element is a $Q$-rational combination of the
elements of $\mathcal X$, and a decomposition $g = \eps g' \gamma$ into
polynomial sequences 
$\varepsilon, g', \gamma \in \mathrm{poly}(\Z,G_{\bullet})$ 
with the following properties:
\begin{enumerate}
\item $\varepsilon : \Z \to G$ is  $(Q,N)$-smooth;
\item $g' : \Z \to G'$ takes values in $G'$, and the finite
sequence $(g'(n)\Gamma')_{n \in [N]}$ is $1/Q^A$-equidistributed
in $G'/\Gamma'$, using the metric $d_{\mathcal X'}$ on $G'/\Gamma'$;
\item $\gamma: \Z \to G$ is $Q$-rational, and 
$(\gamma(n)\Gamma)_{n \in \Z}$ is periodic with period at most $Q$.
\end{enumerate}
\end{theorem} 

The proof of our modified factorisation theorem will proceed via an
iterative application of the theorem stated above.
Our next aim is to prove an auxiliary lemma which will guarantee
that the iteration process stops after finitely many steps.
The way this goal is attained is to ensure that every time we refine
our splitting of $[N]$ into subprogressions the polynomial sequence $g$
we try to factorise fails to be totally equidistributed (with some
parameter) on \emph{each} of the new subprogressions.
This way an application of the factorisation theorem on any new
subprogression yields a \emph{lower} dimensional rational subgroup.

\begin{lemma}\label{residue-independence}
Let $G/\Gamma$ be an $m$-dimensional nilmanifold and let 
$g \in \mathrm{poly}(G_{\bullet}, \Z)$ be a polynomial sequence of
degree $d$. Let $\delta \in (0,1/2)$ be such that 
$\delta^{-t} \ll_t N$ for all $t \in \N$. 
Further let $a \leq \delta^{-1}$ be an integer, and $b \in [a]$.
Suppose that
$(g(q(an + b))\Gamma)_{n \in [N/q]}$ fails to be $\delta$-equidistributed
in $G/\Gamma$ for some $q \ll N^{o(1)}$.
Then there is some $B=O_{m,d}(1)$ such that each of the sequences
$(g(n(aq)^d + r) \Gamma)_{n \in [N/q^d]}$ for 
$r \in [(aq)^d]$ fails to be $\delta^{B}$-equidistributed in $G/\Gamma$. 
\end{lemma}

\begin{proof}
By Proposition \ref{prop-projection} and Proposition \ref{weyl}(a),
there is a non-trivial horizontal character $\eta$ of modulus bounded by
$\delta^{-O_{m,d}(1)}$ such that the function $h:\Z \to \R/\Z$ defined by
$h(n):= \eta \circ g (q(an +b))$ satisfies
$\|h\|_{C^{\infty}[N/aq]} \ll \delta^{-O_{m,d}(1)}$.
Let $\eta \circ g (n) = \sum_{j=0}^d \beta_j n^j$ and 
$\eta \circ g (q(an+b)) = \sum_{j=0}^d \sigma_j n^j$
be polynomial representations in $\R/\Z$.
Then
$$ \sup_{1 \leq j \leq d} \| \sigma_j \| (N/aq)^{j} \ll
\delta^{-O_{m,d}(1)}~.$$
Since
\begin{align} \label{sigma_j--lemma}
\sigma_j 
= \beta_j(aq)^j 
+ \sum_{\ell > j} \binom{\ell}{j} \beta_{\ell}(aq)^j(bq)^{\ell-j}~,
\end{align}
we find, using a downwards induction starting with $j=d$, that
$$ 
\| \beta_j (aq)^{d} \|  \ll \delta^{-O_{m,d}(1)} N^{-j}(aq)^{d} = o(1)~.
$$
Indeed, for $j=d$ the assertion is immediate. 
Suppose now it holds for $j \in \{j_0 + 1, \dots , d\}$ for some
$j_0 \geq 1$. We proceed to check the case where $j=j_0$ by analysing
\eqref{sigma_j--lemma} for $j=j_0$, multiplied through by
$t=(aq)^{d-{j_0}}$. 
Observe that for all positive integers $t$ and for all 
$i\in\{1, \dots,d\}$
$$ 
\|  \sigma_i t \| 
\ll t \|  \sigma_i \| 
\ll t  N^{-i} (aq)^{i} \delta^{-O_{m,d}(1)} ~.
$$
By the assumptions on $\delta, a$ and $q$, this bound is $o(1)$ when
$t=(aq)^{d-{j_0}}$ and $i = j_0$.
Similarly, we have by induction hypothesis for 
$\ell \in \{j_0 + 1, \dots, d\}$ and all $t$
$$ 
\|  \beta_{\ell} (aq)^{d} t \| 
\ll t \|  \beta_{\ell} (aq)^{d} \| 
\ll t  N^{-\ell} (aq)^{d} \delta^{-O_{d,m}(1)}~,
$$
which certainly is $o(N^{-j_0}(aq)^{d} \delta^{-O_{d,m}(1)})$ if we set 
$t=(qb)^{\ell - j_0}$.
This allows us to apply the triangle inequality to split up 
$\| \sigma_{j_0} (aq)^{d-{j_0}}\|$ in the manner of \eqref{sigma_j--lemma}
to deduce the assertion for $j_0$.

Next, pick $r \in [(aq)^d]$ and define 
$\tilde \sigma_0, \dots, \tilde\sigma_d$ such that
$\eta \circ g((aq)^d n + r) = \sum_{j=0}^d \tilde \sigma_j n^j$, thus
$$\tilde \sigma_j
= \sum_{\ell=j}^d \binom{\ell}{j} r^{\ell - j}(aq)^{jd} \beta_{\ell}~.$$
Since $jd\geq d $ for all $j \in \{1, \dots, d\}$, we have
for each of the summands
\begin{align*}
\left\| \binom{\ell}{j} r^{\ell - j}(aq)^{jd} \beta_{\ell} \right\| 
&\ll \binom{\ell}{j} r^{\ell - j}(aq)^{(j-1)d} \|(aq)^{d}\beta_{\ell}\| \\
&\ll_d r^{\ell - j}(aq)^{jd} N^{-\ell} \delta^{-O_{d,m}(1)}
\ll_d (aq)^{\ell d} N^{-\ell} \delta^{-O_{d,m}(1)}~.
\end{align*}
By the assumptions on $\delta$ and $q$, this bound equals $o(1)$ and hence
we can apply the triangle inequality to split up $\| \tilde \sigma_j \|$:
$$
\| \tilde \sigma_j \| 
\ll_d \sum_{\ell = j}^d (aq)^{\ell d} N^{-\ell} \delta^{-O_{d,m}(1)}
\ll_d (N/(aq)^d)^{-j} \delta^{-O_{d,m}(1)}~.
$$
By Proposition \ref{weyl}(b) and Proposition \ref{prop-projection}, this
implies the result.
\end{proof}

Now we finally turn to the modified factorisation theorem which
gives the correct type of minor arcs.

\begin{theorem}[Modified factorisation theorem]
\label{modified-factorisation}
Let $m,d,N,A \geq 1$ be integers, and let $k, Q_0, R \geq 1$ be
integer parameters.
Suppose that $G/\Gamma$ is an $m$-dimensional nilmanifold together with a
filtration $G_{\bullet}$ of degree $d$. 
Suppose that $\mathcal X$ is a $Q_0$-rational Mal'cev basis $\mathcal X$
adapted to $G_{\bullet}$ and that $g \in \mathrm{poly}(\Z,G_{\bullet})$.
Suppose further that $Q_0 \ll \log k$ and 
$k,R=O(N^{o(1)})$. 
Then there is an integer $Q$ with 
$Q_0 \leq Q \ll Q_0^{O_{A,m,d}(1)}$,
and a partition of $[N]$ into at most $R^{dm}$ disjoint subprogressions
$P$, each of length at least $N/R^{dm}$ and each of $k$-smooth
common difference bounded by $R^{dm}$ such that the restriction of
$(g(n))_{n \in P}$ to any of the progression $P$ can be factorised as
follows.

There is a rational subgroup $G' \leq G$, depending on $P$, and 
a Mal'cev basis $\mathcal X'$ for $G'/\Gamma'$ such that 
every element of $\mathcal X'$ is a $Q$-rational combination of elements
from $\mathcal X$ (that is, each coefficient is rational of height
bounded by $Q$).
Suppose $P=\{n \equiv r \Mod{q}\}$, then we have a factorisation
$$g(q n + r) = \eps_P(n) g'_P(n) \gamma_P(n)~,$$
where $\eps_P, g'_P, \gamma_P$ are polynomial
sequences from $\mathrm{poly}(\Z, G_{\bullet})$ with the properties
\begin{enumerate}
 \item $\eps_P: \Z \to G$ is $(Q,N/q)$-smooth;
 \item $g'_P: \Z \to G'$ takes values in $G'$ and for each  $k$-smooth
number $\tilde q \leq R$ the finite sequence 
$(g'_P(\tilde q n)\Gamma')_{n \leq N/(q \tilde q)}$ is totally
$Q^{-A}$-equidistributed in $G'/\Gamma'$;
 \item $\gamma_P: \Z \to G$ is $Q$-rational and 
$(\gamma_i(n)\Gamma)_{n \in \Z}$ is periodic with a $k$-smooth period
which is bounded by $R^{md}Q$.
\end{enumerate}
\end{theorem}

\begin{proof}
We may suppose that $g$ does not satisfy (2), that is, there is some
$k$-smooth integer $q_1 \leq R$  and $b_1 < a_1 \leq Q_0^A$ such that 
$(g(q_1(a_1 n + b_1))\Gamma)_{n \leq N/q_1}$ fails to be 
$Q_0^{-A}$-equidistributed.
Writing $z_1:=(a_1q_1)^d$, Lemma \ref{residue-independence} implies
that each of the sequences $(g(z_1 n + r_1)\Gamma)_{n \leq N/z_1}$ with
$r_1 \in [z_1]$ fails to be $Q_0^{-AA'}$-equidistributed for some
$A'=O_{m,d}(1)$.
Now, we run through all $r_1 \in [z_1]$ in turn.

Applying the factorisation theorem in its original form to any of these
sequences yields some $Q_1 \ll Q_0^{O(A,m,d)}$, a proper $Q_1$-rational
subgroup $G_1 < G$ of dimension strictly smaller than $m$, and a
factorisation 
$$g(z_1 n + r_1) = \eps_{r_1}(n) g'_{r_1}(n)\gamma_{r_1}(n)$$
where the finite sequence 
$(g'_{r_1}(n) \Gamma_1)_{n\leq N/z_1}$ is totally
$Q_1^{-A}$-equidistributed in  $$G_1/\Gamma_1:=G_1/(\Gamma \cap G_1)~.$$

If $g'_{r_1}$ is $Q_1^{-A}$-equidistributed on every subprogression 
$\{n \equiv b_2 \Mod{a_2q_2} \}$ of $k$-smooth common difference $a_2q_2$,
where $b_2<a_2<Q_1^{A}$ and $q_2<R$, then we stop (and turn to the next
choice of $r_1$). 
Otherwise, invoking Lemma \ref{residue-independence} again, there is a
$k$-smooth integer $a_2q_2$ as above such that with $z_2:=(a_2q_2)^d$ the
finite sequence
$(g_{r_1,r_2}(n))_{n \leq N/(z_1z_2)}$ defined by 
$g_{r_1,r_2}(n) := g'_{r_1}(z_2 n + r_2)$ is not
$Q_1^{-A}$-equidistributed for any $r_2 \in [z_2]$. 
We proceed as before.

This process yields a tree of operations which has height at most $m =
\dim G$, since each time the factorisation theorem is applied a new
sequence $g'_{r_1, \dots, r_i}$ is found that takes values in some
strictly lower dimensional submanifold 
$G_i= G_i(r_1, \dots, r_i)$ of $G_{i-1}(r_1, \dots, r_{i-1})$. 
Thus, we can apply the factorisation theorem at most $m$ times in a row
before the manifold involved has dimension $0$. 

The tree we run through starts with $g$, which has $z_1$
neighbours $g_{r_1}$, one for each $r_1 \in [z_1]$. 
Each $g_{r_1}$ has $z_2=z_2(r_1,r_2)$ neighbours $g_{r_1,r_2}$, one for
each $r_2 \in [z_2]$, etc..

As a result, we obtain a decomposition of the range $[N]$ into at most
$R^{2dm}$ subprogressions of the form 
\begin{align*}
P
& = \{z_1(z_2(z_3( \dots (z_t m + r_t)
     \dots) + r_3) + r_2) + r_1 
    ~:~ m \leq N/(z_1z_2 \dots z_t) \} \\
& = \{ z_1z_2 \dots z_t m + r 
    ~:~ m \leq N/(z_1z_2 \dots z_t)\} ~,
\end{align*}
for some $r$, and where each $z_i$ depends on $r_1, \dots,r_{i-1}$.
The common difference of such a progression $P$ is
$k$-smooth and bounded by $R^{2dm}$. 
Thus, $P$ has length at least $N/R^{2dm}=N^{1-o_m(1)}$.
The iteration process furthermore yields a factorisation of 
$g_{r_1, \dots, r_t}$, which is the restriction of $g$ to $P$:
$$ g_{r_1,\dots,r_t} (m)
=g(z_1z_2 \dots z_t m + r) 
= \tilde \eps_{r_1,\dots,r_t}(m)
  g'_{t}(m)
  \tilde \gamma_{r_1,\dots,r_t}(m)~,
$$ 
where
$$\tilde \eps_{r_1,\dots,r_t}(m)
= \eps_{r_1}( z_2 \dots z_t m + \tilde r_2) \dots
  \eps_{r_1,\dots,r_{t-1}}(z_t m + \tilde r_t)
  \eps_{r_1,\dots,r_t}(m)
$$
for certain integers $\tilde r_2, \tilde r_3, \dots, \tilde r_t$,
and
$$ \tilde \gamma_{r_1,\dots,r_t}(m)
  =
  \gamma_{r_1,\dots,r_t}(m)
  \gamma_{r_1,\dots,r_{t-1}} (z_t m + \tilde r_t) \dots
  \gamma_{r_1}(z_2 \dots z_t m + \tilde r_2)~.
$$
In view of the remarks following the definition of smoothness of
sequences, the factor
$\tilde \eps_{r_1,\dots,r_t}(m)$
is a $(Q_0^{O_{A,d,m}(1)},N/ (z_1 \dots z_t) )$-smooth sequence.
Further, the periodic sequences
$\tilde \gamma_{r_1,\dots,r_t}(m)$
are easily seen to have a $Q_0^{O_{A,d,m}(1)}$-smooth, i.e.~$k$-smooth,
period.
\end{proof}

\section{Reduction to the case of minor arc nilsequences}

With the help of the modified factorisation theorem,
Theorem \ref{modified-factorisation}, we will show that the
general non-correlation estimate follows from the special case of
non-correlation with `minor arc nilsequences' that enjoy property (ii) of
the modified factorisation theorem.

The general case is the following proposition.
\begin{proposition}\label{non-corr}
Let $G/\Gamma$ be a nilmanifold of dimension $m \geq 1$, let
$G_{\bullet}$ be a filtration of $G$ of degree $d \geq 1$, and let
$g \in poly(\Z,G_{\bullet})$ be a polynomial sequence. Suppose that
$G/\Gamma$ has a $Q$-rational Mal'cev basis $\mathcal X$
for some $Q \geq 2$, defining a metric $d_{\mathcal X}$ on $G/\Gamma$.
Suppose that $F: G/\Gamma \to [-1,1]$ is a Lipschitz function. Then we
have for $M_0= \log \log \log N$ and $N'= \lfloor N/\W \rfloor$
$$|\E_{n\in[N']} (r'_{f,\beta}(n) - 1) F(g(n))\Gamma| 
\ll_{m,d,\gamma,A} 
  Q^{O_{m,d,\gamma,A}(1)} (1+ \|F\|) M_0^{-A}$$
for any $A>0$ and $N \geq 2$.
\end{proposition}

Similarly as in \S2 of \cite{green-tao-nilmobius}, we will deduce this
result from the following special case involving only
`minor arc nilsequences'. 

\begin{proposition}[Non-correlation, equidistributed case]
\label{non-corr,equid.case}
Let $N>0$ be a large integer and let $\delta$, $k$ and $R$ be
parameters such that $\delta \in (0,1/2)$, $\delta^{-t}\ll_t N'$ for all 
$t\in \N$, $R \ll N^{o(1)}$ and $k=w(N)$. 
 Suppose that $(G/\Gamma, d_{\mathcal X})$ is an $m$-dimensional
nilmanifold with some filtration $G_{\bullet}$ of degree $d$ and suppose
that $g \in \mathrm{poly}(\Z,G_{\bullet})$. 
Suppose further that for every $k$-smooth number $\tilde q \leq R$ the
finite sequence $(g(\tilde q n)\Gamma)_{n\in [N'/\tilde q]}$ is
$\delta$-equidistributed in $G/\Gamma$. 
Then for every Lipschitz function $F:G/\Gamma \to \R$ satisfying
$\int_{G/\Gamma}F=0$ and for every $k$-smooth number $q \ll N^{o(1)}$ and
every $r \in [q]$, we have
$$|\E_{n\in[N']} (r'_{f,\beta}(qn+r)-1) F(g(n)\Gamma)| \ll
\delta^{c} \|F\|$$
for some $c$ such that $c^{-1}=O_{m,d}(1)$.
\end{proposition}

\begin{proof}[Proof of Proposition \ref{non-corr} assuming Proposition 
\ref{non-corr,equid.case}]
Observe that $N'=N^{1-o(1)}$. 
We may assume that $Q \leq M_0$, thus $Q \leq M_0= \log w(N)$.
The modified factorisation theorem can now be applied to the sequence
$(g(n)\Gamma)_{n \leq N}$ with the following parameters: $k=w(N)$,
$Q_0=\log w(N)$, $R= N^{o(1)}$.
This yields a partition of $[N']$ into at most $R^{2md}$ progressions of
$w(N)$-smooth common differences.
By the triangle inequality, it suffices to show that
$$|\E_{n\in P } (r'_{f,\beta}(n) - 1) F(g(n)\Gamma)|
\ll_{m,d,\gamma,A} 
   Q^{O_{m,d,\gamma,A}(1)} (1+ \|F\|) M_0^{-A}$$
for every progression $P$ in the partition.

For each of these progressions $P=:\{q_P n + r_P\}$, the modified
factorisation theorem provides us with a factorisation of the restriction
of $g$ to $P$: 
$$g(q_P n + r_P)=: g_P(n) = \eps_P(n) g'_P(n) \gamma_P(n)~,$$
where $\eps_P, g'_P, \gamma_P$ satisfy $(i)$, $(ii)$ and $(iii)$ from
Theorem \ref{modified-factorisation}.
Proceeding as in \cite[\S2]{green-tao-nilmobius} 
(see \emph{loc.~cit.}~for full details), we split each $P$ into
subprogressions $P=P_1 \cup \dots \cup P_t$ in such a way that
\begin{itemize}
 \item $\gamma_P(n)$ is constant on each progression, 
say $\gamma_P(n)=\gamma_j$ for $n \in P_j$, and
\item $\eps_P(n)$ is almost constant: to be precise, the $P_j$ are such
that $|n-n'| \leq N'/(qQ^{B})$ for some $B=O(1)$ and all $n,n' \in P_j$
which implies $d(\eps_P(n),\eps_P(n')) \leq Q^{-B+1}$ by smoothness of
$\eps_P$.
\end{itemize}
From each $P_j$, we choose a fixed element, say $n_j$. 
Then the Lipschitz property of $F$, right-invariance of the
metric, and smoothness of $\eps_P$ imply that for every $n \in P_j$
\begin{align*}
 |F(\eps_P(n) g'_P(n) \gamma(n) \Gamma)
 -F(\eps_P(n_j) g'_P(n) \gamma_j \Gamma)| \leq Q^{-B/2}~,
\end{align*}
provided $B$ was chosen large enough.
Hence it suffices to show that
$$|\E_{n\in P_j } 
(r'_{f,\beta}(n) - 1) 
F(\eps_P(n_j) \gamma_j (\gamma_j^{-1} g'_P(n) \gamma_j) \Gamma)| 
\ll_{m,d,\gamma,A} 
   Q^{O_{m,d,\gamma,A}(1)} (1+ \|F\|) M_0^{-A}~.$$
The aim is now to apply Proposition \ref{non-corr,equid.case} to
$g_j:\Z \to \gamma_j^{-1} G \gamma_j=:H_j$, 
$$g_j(n) := \gamma_j^{-1} g'_P(q_{P_j}n + r_{P_j}) \gamma_j~.$$
Property $(ii)$ of the modified factorisation theorem was set up so as
to ensure that $g_j$ still enjoys the `minor arc property' (on
$H_j/(\Gamma \cap H_j)$ rather than $G/\Gamma$, of course).
Note that the Lipschitz constant of $F_j: H_j/(\Gamma \cap H_j) \to \C$, 
$F_j(x(\Gamma \cap H_j)):= F(\eps_P(n_j) \gamma_j \Gamma)$ is
bounded by $M^{} \|F\|$ by \cite[Lemma A.16]{green-tao-polynomialorbits}.
Since $P_j$ has a $w(N)$-smooth common difference and length at least 
$N^{1-o(1)}$, Proposition \ref{major_arc_estimate} implies that
$(r_{f,\beta} - 1)$ does not correlate with any function 
$n \mapsto c1_{P_j}(n)$, where $c$ is a constant.
Hence we can subtract off the mean value of $F_j$ and reduce to the
assumption $\int_{H_j/\Lambda_j} F_j = 0$.

All remaining technical details work exactly as in
\cite[\S2 and App.B]{green-tao-nilmobius}, so we have chosen, given
their technical complexity, to omit them here.
\end{proof}

\section{Completion of the non-correlation estimate}
\label{conclusion-section}
We complete the proof of Proposition \ref{non-corr,equid.case} and
therefore the analysis of correlation of $r'_{f,\beta}$ with nilsequences.
Recall the conditions of Proposition \ref{non-corr,equid.case}.
In particular, we are given a polynomial sequence
$(g(n)\Gamma)_{n\in[N']}$ such that for every
$w(N)$-smooth number $\tilde q \leq R$ the finite sequence 
$(g(\tilde q n)\Gamma)_{n\in [N'/\tilde q]}$ is $\delta$-equidistributed
in $G/\Gamma$. The parameter $\delta$ satisfies the condition 
$\delta^{-t} \ll_t N'$,  which will allow us later to apply Proposition
\ref{modified_general_equid}.
We are required to show that for every Lipschitz function $F:G/\Gamma
\to \R$ satisfying $\int_{G/\Gamma}F=0$, for every $w(N)$-smooth number
$q \ll N^{o(1)}$, and for every $r \in [q]$, we have
$$|\E_{n\in[N'/q]} (r'_{f,\beta}(qn+r)-1) F(g(n)\Gamma)| \ll
\delta^{1/O_{m,d}(1)} \|F\|~.$$
By $\delta$-equidistribution of 
$(g(n)\Gamma)_{n\in [N]}$ and since $\int F=0$, it
suffices to show that
$$|\E_{n\in[N'/q]} r'_{f,\beta}(qn+r) F(g(n)\Gamma)| \ll
\delta^{1/O_{m,d}(1)} \|F\|~.$$

\begin{figure}
\centering
\includegraphics{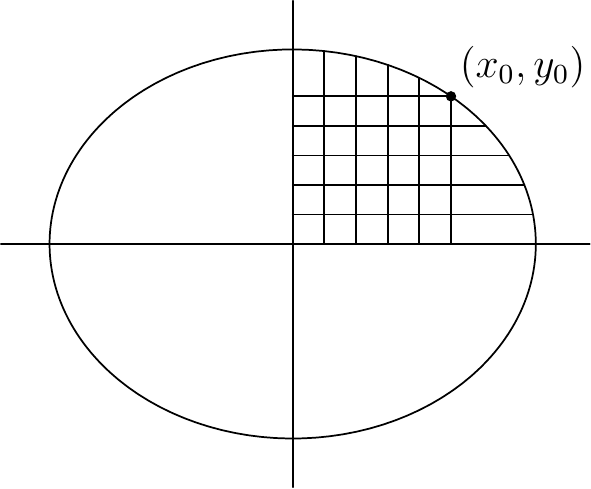}

%
%
%
%
%
%
%
%
%
%
%
%
%
%

\caption{Schematic of how the summation is split: we sum along
horizontal lines ($y$ fixed), along vertical lines ($x$ fixed), and
another time over the segments of the horizontal lines that are contained
in the `box', that is, over the double-counted segments.}
\label{picture}
\end{figure}

We may suppose that $f = \<a,b,c\>$ has reduced form, that is 
$|b|\leq a \leq c$.
Writing $X(f, N) := \{(x,y): f(x,y) \leq N \}$, our aim is to decompose
the binary sequence 
$$\{g((ax^2 + bxy +cy^2)\Gamma)\}_{(x,y) \in X(f,N)}$$
into a sum of polynomial subsequences
$(g'(P(n) \Gamma))_{n \leq (N')^{1/\deg(P)}}$ of some equidistributed
sequence $(g'(n) \Gamma)_{n \leq N'}$.
In order to do so, let $(x_0,y_0) \in \R^2$ be the point
on the ellipse $f(x,y)=ax^2 + bxy +cy^2 = N$ that satisfies 
$x_0 = y_0 \sim N^{1/2}$.
Since $f$ has reduced form, both $a x_0^2 \leq N$ and 
$c y_0^2 \leq N$ hold. 
With respect to $(x_0,y_0)$, the summation over $(x,y) \in X(f,N)$ now
splits into three parts (cf.~Figure \ref{picture}) such that on each
part one of the variables $x$ and $y$ may be fixed, while the free
variable will range over an interval of length at least 
$x_0 \sim N^{1/2}$.
This decomposition yields
\begin{align}
\nonumber
&  
\frac{2\pi}{\sqrt{-D}}
~
\Big| \sum_{n \leq (N'-r)/q} r'_{f,\beta}(qn+r) F(g(n) \Gamma) \Big|
\\
\nonumber
\leq& \quad \bigg(\frac{\rho_{f,\beta}(\W)}{\W}\bigg)^{-1}
   \sum_{y \leq y_0} \bigg| \sum_{ x : f(x,y) \leq N} 
    1_{f(x,y) \equiv \W r + \beta \Mod{\W q}}
    F\Big(g\Big(\frac{f(x,y) -\beta -\W r}{\W q}\Big)\Gamma\Big) \bigg| \\
\nonumber
&+ \bigg(\frac{\rho_{f,\beta}(\W)}{\W}\bigg)^{-1}
    \sum_{x \leq x_0} \bigg| \sum_{ y : f(x,y) \leq N} 
    1_{f(x,y) \equiv \W r + \beta \Mod{\W q}}
    F\Big(g\Big(\frac{f(x,y) -\beta -\W r}{\W q}\Big)\Gamma\Big) \bigg| \\
\label{3-summands}
&+ \bigg(\frac{\rho_{f,\beta}(\W)}{\W}\bigg)^{-1}
    \sum_{y \leq y_0} \bigg| \sum_{ x \leq x_0 } 
    1_{f(x,y) \equiv \W r + \beta \Mod{\W q}}
    F\Big(g\Big(\frac{f(x,y) - \beta - \W r}{\W q}\Big)\Gamma\Big) \bigg| 
  ~.
\end{align}

To remove the congruence condition 
$f(x,y) \equiv \W r + \beta \Mod{\W q}$ in
this explicit form, we consider the set $S(q\W,\W r + \beta)$ of all
solutions $(x',y') \in [q\W]^2$ to the congruence $f(x',y') \equiv \W r +
\beta \Mod{q\W}$. 
By Corollary \ref{lifted-densities} the density of these
solutions for a $w(N)$-smooth integer $q$, $r \in [q]$ and $\beta \in \A$
satisfies
\begin{align*}
\frac{\rho_{f,\W r + \beta}(q\W)}{q\W} 
= \frac{\rho_{f,\beta}(\W)}{\W}~. 
\end{align*}
To simplify the notation, define $\beta':= \W r + \beta$ and
$q':= \W q$.
Considering any of the three parts of our summation above, we may continue
this as follows 
\begin{align*}
  &~ \bigg(\frac{\rho_{f,\beta}(\W)}{\W}\bigg)^{-1}
    \sum_{y \leq y_0} \bigg| \sum_{ x : f(x,y) \leq N} 
    1_{f(x,y) \equiv \beta' \Mod{q'}}
    F\Big(g\Big(\frac{f(x,y) -\beta'}{q'}\Big)\Gamma\Big) \bigg|
\\ 
&= q\W ~\E_{(x',y') \in S(q',\beta')}
    \sum_{y: q' y + y' \leq y_0} 
   \bigg| \sum_{ \substack{ x : f(q' x + x', q' y + y') \\ \leq N}}
    F\Big(g\Big(\frac{f(q' x + x', q' y + y') - \beta'}{q'}
      \Big) \Gamma \Big) \bigg|~.   
\end{align*}
Observe that
\begin{align*}
\frac{f(\W q x + x',\W q y + y') - \W r - \beta }{q\W} 
= \W q a x^2 + b' x + c'~,
\end{align*}
for some $b',c'$ depending on $y,y',x',b,c,q$ and $\W$, is a polynomial
that satisfies the conditions of Proposition \ref{modified_general_equid}.
Thus, setting $P(x):= \W q a x^2 + b' x + c'$, we are considering the
polynomial subsequence $(g \circ P (n) \Gamma)_{n \leq (N'/qq')^{1/2}}$
of $(g(n)\Gamma)_{n \leq N'/q}$.
By Proposition \ref{modified_general_equid} there is for each $P$ a
$w(N)$-smooth integer $\tilde q \ll N^{o(1)}$ such that
for every $\tilde r \in [\tilde q]$ the sequence
$$ g(P(\tilde q x + \tilde r))_{x \leq N^{1/2}/ (\W q \tilde q)} $$
is totally $\delta^{1/O_{d}(1)}$-equidistributed.
Splitting the summation into subprogressions modulo $\tilde q$, we have
via the triangle inequality
\begin{align*}
&q\W ~\E_{(x',y') \in S(q',\beta')}
    \sum_{y: q' y + y' \leq y_0}
\bigg| \sum_{ \substack{ x : f(q' x + x', q' y + y') \\ \leq N}}
    F (g \circ P(x) \Gamma )~ \bigg| \\
&\leq q\W  ~\E_{(x',y') \in S(q',\beta')}
    \sum_{y: q' y + y' \leq y_0} \sum_{\tilde r}
 \bigg| 
 \sum_{ \substack{ x : f(q'(\tilde q x + \tilde r)+x',q'y +y')\\ \leq N}}
 F (g \circ P(\tilde q x + \tilde r) \Gamma )~ \bigg|  \\
&\ll q\W  ~\E_{(x',y') \in S(q',\beta')}
    \sum_{y: q' y + y' \leq y_0}
    \tilde q ~\frac{N^{1/2}}{\W q \tilde q}~ \delta^{1/O_d(1)}
\|F\|_{\mathrm{Lip}} \\
&\ll \delta^{1/O_d(1)} \frac{N'}{q} \|F\|_{\mathrm{Lip}}~.
\end{align*}
As these arguments also apply to the two remaining parts of the sum
\eqref{3-summands} this completes the proof of Proposition
\ref{non-corr,equid.case} and also the proof of the main theorem.

\subsection*{Acknowledgements}
I should like to thank my PhD supervisor Ben Green for suggesting the
problems studied in this paper and for many valuable discussions and
advice.
I am also very grateful to Tim Browning for insightful comments and
suggestions, and to Tom Sanders for helpful conversations.

\providecommand{\bysame}{\leavevmode\hbox to3em{\hrulefill}\thinspace}


\begin{thebibliography}{99}

\bibitem{browning-breteche} R.~de la Bret{\`e}che and T.~D.~Browning,
\emph{Binary linear forms as sums of two squares,}
Compos.~Math. {\bf 144} (2008), no.~6, 1375--1402.

\bibitem{browning-munshi} T.D.~Browning and R.~Munshi,
\emph{Rational points on singular intersections of quadrics},
{arXiv:1108.1902}, 2011.

\bibitem{cook} R.~J.~Cook, 
{\em Simultaneous quadratic equations,\/} 
{J. London Math. Soc.} {\bf s2-4 (2)} (1971), {319--326}.

\bibitem{cox} D.~A.~Cox,
{\em Primes of the form $x^2 + ny^2$,\/}
Pure and Applied Mathematics, Wiley, 1989.

\bibitem{erdos} P.~Erd\H{o}s,
\emph{On the sum $\sum_{k=1}^x d(f(k))$,} 
{J. London Math. Soc.} {\bf 27} (1952), no.~1, {7--15}.

\bibitem{green-tao-longprimeaps} B.~J.~Green and T.~C.~Tao,
\emph{The primes contain arbitrarily long arithmetic progressions,}
{Annals of Math.} {\bf 167 } (2008), No.~2, 481--547.

\bibitem{green-tao-quadraticmöbius} \bysame, 
\emph{Quadratic uniformity of the M\"obius function,}
Annales de l'Institut Fourier (Grenoble) {\bf 58} (2008), no.~6,
1863--1935.

\bibitem{green-tao-linearprimes} \bysame, 
\emph{Linear equations in primes,} 
{Annals of Math.,} {\bf 171 } (2010), No.~3 , 1753--1850.

\bibitem{green-tao-nilmobius} \bysame, 
\emph{The M\"obius function is strongly orthogonal to nilsequences,}
{Annals of Math.} {\bf 175} (2012), No.~2, 541--566.

\bibitem{green-tao-polynomialorbits} \bysame, 
\emph{The quantitative behaviour of polynomial orbits on nilmanifolds,}
{Annals of Math.} {\bf 175} (2012), No.~2, 465--540.

\bibitem{gtz} B.~J.~Green, T.~C.~Tao and T.~Ziegler
\emph{An inverse theorem for the Gowers $U^k[N]$ norm},
{Annals of Math.,} to appear.
Preprint available at 
{arXiv:1009.3998v2}.

\bibitem{heath-brown} D.~R.~Heath-Brown,
Linear Relations Amongst Sums of two Squares. 
{\em Number theory and algebraic geometry, \/} 133--176, London Math. Soc.
Lecture Note Ser. {\bf 303}, Cambridge Univ. Press, Cambridge, 2003.

\bibitem{henriot} K.~Henriot,
\emph{Nair--Tenenbaum bounds uniform with respect to the discriminant,}
Math.~Proc.~Camb. Phil.~Soc., {\bf 152} (2012), 405--424.

\bibitem{hua} L-K.~Hua,
\emph{Introduction to Number Theory,}
{Springer-Verlag,}
{Berlin, Heidelberg,}
{1982}.


\bibitem{iwaniec} H.~Iwaniec, 
\emph{The half dimensional sieve,}
{Acta Arith.},
{\bf 29} {(1976)},{no.1}, {69--95}.

\bibitem{IK} H.~Iwaniec and E.~Kowalski, 
{\em Analytic Number Theory, \/}
Colloquium Publications, vol. {\bf 53}, American Mathematical
Society, Providence, RI, 2004.

\bibitem{landau} E.~Landau,
{\em \"Uber die Einteilung der positiven ganzen Zahlen in vier Klassen
nach der Mindestzahl der zu ihrer additiven Zusammensetzung
erforderlichen Quadrate, \/}
Archiv der Mathematik und Physik, ser.3, vol. 13 (1908), 305--312.

\bibitem{m-divisorfunction} L.~Matthiesen,
\emph{Correlations of the divisor function}, 
{Proc.~London Math.~Soc.,} 
{\bf 104} (2012), 827--858.
Preprint available at 
{arXiv:1011.0019}.

\bibitem{MV}  H.~L.~Montgomery and R.~C.~Vaughan, 
\emph{Multiplicative Number Theory I, Classical Theory},
Cambridge Studies in Advanced Mathematics {\bf 97}, Cambridge University
Press, 2006.

\bibitem{prachar} K.~Prachar,
{\em \"Uber Zahlen der Form $a^2+b^2$ in einer arithmetischen Progression,
\/}
Math. Nachr. 10 (1953), 51--54.

\bibitem{rose} H.~E.~Rose,
\emph{A Course in Number Theory,}
second ed.,
Oxford science publications, Clarendon Press, Oxford, 1994.

\bibitem{stewart} C.~L.~Stewart, 
\emph{On the number of solutions of polynomial congruences and Thue
equations}, 
J.~Amer.~Math.~Soc. {\bf 4} (1991), 793--835.

\end{thebibliography}
\end{document}